%Paper: q-alg/9503022
%From: David Kazhdan <kazhdan@abel.math.harvard.edu>
%Date: Fri, 31 Mar 1995 14:02:42 -0500

\input vanilla.sty
\scaletype{\magstep1}
\overfullrule=0pt
\input amssym.def
\input amssym

\def\ident{1\!\!1}
\def\llongrightarrow{\mskip-\medmuskip \mkern5mu
   \mathbin{\vrule height 2.7pt depth -2.3pt width 30pt
   \mathrel{\mkern-5mu}\longrightarrow} \penalty900 \mkern2mu
\mskip-\medmuskip}

\def\llongrightarrow{\mskip-\medmuskip \mkern5mu
   \mathbin{\vrule height 2.7pt depth -2.3pt width 30pt
   \mathrel{\mkern-5mu}\longrightarrow} \penalty900 \mkern2mu
\mskip-\medmuskip}

\baselineskip = 18pt

%this file is qkzint.tex  Introduction

\newpage

\baselineskip=12pt

\centerline{\bf Representations of  Quantum Affine Algebras}
\vskip .2in

$$
{\matrix \format \l &\qquad\qquad\qquad  \l \\
\text{D.\ Kazhdan}\\
\text{Dept.\ of\ Mathematics}\\
\text{Harvard\ University}\\
\text{Cambridge,\ MA\ 02138}
\endmatrix}\qquad\qquad
{\matrix \format \l &\qquad\qquad\qquad  \l \\
\text{Y.\ Soibelman}\\
\text{Dept.\ of\ Mathematics}\\
\text{Kansas\ State\ University}\\
\text{Manhattan, KS 66506}
\endmatrix}
$$
\vskip .3in

\baselineskip=18pt

\centerline{\bf Introduction}
\vskip .2in

Let  $\overline{\frak G}$  be a finite dimensional simple Lie algebra and
${\frak G}$  the corresponding affine Ka$\check{\text c}$-Moody algebra.  The
notion of the
fusion in the category ${\Cal O}$  of representations of affine
Ka$\check{\text c}$-Moody algebras ${\frak G}$  was introduced ten years ago
by physicists in the framework of Conformal Field Theory.  This notion was
developed in a number of mathematical papers (see, for example, [TUY]), where
the notion of fusion is rigorously defined and [D1] where the relation between
the fusion and quantum groups in the quasiclassical region was established).
This line of development was extended in [KL] to a construction on equivalence
between the ``fusion'' category for an arbitrary negative charge and the
category of representations of the  corresponding quantum group (for
simply-laced affine Ka$\check{\text c}$-Moody algebras).

It is natural to try to define the notion of fusion for the category ${\Cal O}$
for affine quantum groups.  One can hope that it might be useful for a deformed
CFT (see [FR]).  The usual construction does not work since it is
based on the existence of a subalgebra  $\Gamma \subset {\frak G} \times
{\frak G} \times {\frak G}$  which is a central extension of the Lie algebra
$\Gamma({\Bbb P}^1 - \{ 0,1,\infty\}, \overline{\frak G})$.  Unfortunately
such subalgebra  $\Gamma$  does not have a natural generalization for the
quantum case (see [D3] where the problem was stated).

In the remarkable paper [FR],\ I. Frenkel and N. Reshetikhin found a way to
describe the fusion between finite-dimensional representations and
representations from the category ${\Cal O}$ for affine quantum groups.
In the present paper we reconsider this problem from the point of view of [KL].
Our main result is the construction of the quasi-associativity constraints
(see Sect.3 and 5). Main phenomena which should be mentioned in connection
with the problem are the appearance of elliptic curves instead of genus zero
curves and {\bf Z}-sheaves instead of bundles with flat connections. From a
general point of view this reflects the fact that the categories we are
considering do not carry monoidal structure.  Nevertheless those
categories of representations of affine quantum algebras carry some other
interesting structures discussed in Sect.1 in the general situation.
These structures explain the categorical meaning of the so-called quantum
Knizhnik-Zamolodchikov equations introduced in [FR].
Among other results we can mention meromorphicity of the quantum $R$-matrix
for any two finite-dimensional representations (see Sect. 4) which has also
general categorical origin.

$Acknowledgements$. It is a pleasure for us to thank P. Deligne, V. Drinfeld,
P.Etingof and N.Reshetikhin for the useful discussions on the topic.  Y.S. is
grateful to the Department of Mathematics of Harvard University, the Institute
for Advanced Study and  the Max-Planck Institute f\"ur Mathematik for the
support and hospitality during various stages of this work.

%this file is qkz1a.tex	  [contains 1.1*]

\newpage

\pageno=3

\centerline{\bf \S1.\ Monoidal categories}
\vskip .3in

\subheading{1.1 Definitions and basic properties}

\proclaim{{\bf 1.1.1.}\ Definition}\ A monoidal category ${\Cal C}$  is a
triple ${\Cal C}
= (\widetilde{\Cal C}, \otimes,a)$  where  $\widetilde{\Cal C}$  is a category,
$\otimes: \widetilde{\Cal C} \times \widetilde{\Cal C} \to \widetilde{\Cal C}$
is a functor and  $a$  is the natural transformation between the functors
$\otimes(\otimes \times id)$  and $\otimes(id \times \otimes)$ from
$\widetilde{\Cal C} \times \widetilde{\Cal C} \times \widetilde{\Cal C}$ to
$\widetilde{\Cal C}$, $a = \{ a_{X,Y,Z}: (X \otimes Y) \otimes Z
\longrightarrow X \otimes (Y \otimes Z)\}$, $X,Y,Z$  in $\widetilde{\Cal C}$
such that all  $a_{X,Y,Z}$  are isomorphisms, the pentagon axiom is satisfied,
and there exists an object $U$ in $\widetilde{\Cal C}$  and an isomorphism
$u: U \otimes U {\buildrel \sim\over\longrightarrow} U$  such that the functors
$X \to X \otimes U$  and $X \to U \otimes X$  are autoequivalences of
$\widetilde{\Cal C}$  (see {\rm [DM]}).  We  call such a pair $(U,u)$  the
identity
object of  ${\Cal C}$.  (It is clear that the identity object
is determined uniquely up to a unique isomorphism.)

Let ${\Cal C} = (\widetilde{\Cal C}, \otimes, {\Cal A})$  be a monoidal
category, and $(U,u)$  the identity object in ${\Cal C}$.
\endproclaim

\proclaim{Lemma}\ a)  For any object $X$ in $\widetilde{\Cal C}$ there exist
unique isomorphisms  $r_X: X {\buildrel \sim\over\longrightarrow} X~\otimes~U$,
$\ell_X: X \to U \otimes X$  such that
$$
a_{X,U,U} \circ (r_X \otimes \text{id}_U) = id_X \otimes u\quad \text{and}\quad
\text{id}_U \otimes
\ell_X = a_{U,U,X} \circ (u \otimes \text{id}_X).
$$

b)  $r_U = u$.

c)  For any morphism  $\alpha: Y \to X$  in $\widetilde{\Cal C}$  the diagrams
$$
{\matrix \format \l &\quad \c &\quad \l \\
Y & {\buildrel \alpha\over\longrightarrow} & X \\
\big\downarrow r_Y && \big\downarrow r_X \\
Y \otimes U & {\buildrel {\alpha \otimes 1}\over\longrightarrow} & X \otimes U
\endmatrix}
\qquad\qquad
{\matrix \format \l &\quad \c &\quad \l \\
Y & {\buildrel \alpha\over\longrightarrow} & X \\
\big\downarrow \ell_Y && \big\downarrow \ell_X \\
U \otimes Y & {\buildrel {\alpha \otimes 1}\over\longrightarrow} & U \otimes X
\endmatrix}
$$
are commutative.
\endproclaim

\demo{Proof}  Well known.  (See [DM].)
\enddemo

\proclaim{{\bf 1.1.2}\ Definition}  We say that a monoidal category  ${\Cal C}
=
(\widetilde{\Cal C}, \otimes, a)$
is strict if for any $X,Y,Z$  in  $\widetilde{\Cal C}$\ we have
$(X \otimes Y) \otimes Z = X \otimes (Y \otimes Z)$, $a_{X,Y,Z} = \text{id}$,
there exists an object $\ident$ in $\widetilde{\Cal C}$  such that
$\ident \otimes X = X \otimes \ident = X$  for all $X$ in $\widetilde{\Cal C}$,
and for any  $X,Y$ in $\widetilde{\Cal C}$ the composition
$$
X \otimes Y = (X \otimes \ident) \otimes Y = X \otimes (\ident \otimes Y) =
X \otimes Y
$$
is equal to $\text{id}_{X\otimes Y}$.
\endproclaim

Mac Lane's theorem says that any monoidal category is equivalent to a strict
one
(see [M, Chapter 7]).   To simplify formulas
we will from now on assume that all monoidal categories are strict.
We will also assume that the category $\widetilde{\Cal C}$
is abelian, contains inductive limits and that $\otimes$  is an additive exact
functor.  Let
$F_{\Cal C}\ {\buildrel {\text{def}}\over =}\ End_{\widetilde{\Cal
C}}(\ident)$.  It
is easy to see that  $F_{\Cal C}$  is a commutative ring and for any
$X,Y$  in $\widetilde{\Cal C}$  the group  Hom$_{\widetilde{\Cal C}}(X,Y)$  has
a natural structure of an $F_{\Cal C}$-module.

{}From now on we will denote the category $\widetilde{\Cal C}$  simply as
${\Cal C}$.  We hope that this will not create any confusion.

\subheading{1.1.3}  Assume until the end of this section that  $F_{\Cal C}$  is
a field.  Any  $F_{\Cal C}$-vector space $R$ can be written as an inductive
limit of finite dimensional spaces  $R = {\dsize\lim_\rightarrow} R_i$.  We
define
$R_{\Cal C}\ {\buildrel{\text{def}}\over =}\ {\dsize\lim_\rightarrow} R_i
\otimes \ident$.
It is clear that the object $R_{\Cal C}$  in ${\Cal C}$  does not depend on a
choice of a presentation $R = {\dsize\lim_\rightarrow} R_i$.

We denote by Vec$_{\Cal C}$  the category of  $F_{\Cal C}$-vector spaces and
for
any object $X$ in ${\Cal C}$  we denote by $L_X$  the functor from
Vec$_{\Cal C}$ to itself such that
$$
L_X(R)\ {\buildrel{\text{def}}\over =}\ Hom_{\Cal C}(X, R_{\Cal C})\quad
\text{for\ all}\ R \in \text{Vec}_{\Cal C}.
$$
It is easy to see that the functor $L_X$  satisfies the conditions of Theorem
5.6.3 in [M].  Therefore the functor $L_X$  is representable.

\proclaim{Definition}  We denote by $\langle X\rangle$ the object in
Vec$_{\Cal C}$  which represents the functor $L_X$.
\endproclaim

As follows from the definition of $\langle X\rangle$, for any $X$ in
${\Cal C}$  we have a natural isomorphism
$\langle X\rangle^\vee \cong Hom_{\Cal C}(X,\ident)$  where
$L^\vee\ {\buildrel{\text{def}}\over =}\ Hom(L,F_{\Cal C})$  for any $L$ in
Vec$_{\Cal C}$.

\subheading{1.1.4}\ Let  ${\Cal C}$
be a strict monoidal category.  We denote by  $\overline{\Cal C}$  the set of
isomorphism classes of objects in  ${\Cal C}$.  For any  $X$  in
$\widetilde{\Cal C}$  we denote by  $[X] \in \overline{\Cal C}$  the
equivalence
class of $X$.

\proclaim{Definition}    We say that an object
$X$ in ${\Cal C}$ is rigid if there exists an object $Y$ in
${\Cal C}$ and a pair of morphisms
$$
i_X: \ident \to X \otimes Y \qquad\qquad  e_X: Y \otimes X \to \ident
$$
such that the compositions
$$
\aligned
X &= \ident \otimes X {\buildrel {i_X \otimes \text{id}}\over\longrightarrow}
X \otimes Y \otimes X {\buildrel {\text{id}\otimes e_X}\over\longrightarrow}
X\\
Y &= Y \otimes \ident  {\buildrel {\text{id} \otimes i_X} \over\longrightarrow}
Y \otimes X \otimes Y {\buildrel {e_X\otimes \text{id}}\over\longrightarrow} Y
\endaligned
$$
are equal to $\text{id}_X$ and $\text{id}_Y$  correspondingly.  In this case we
say that $Y$ is dual to $X$.
\endproclaim

It is easy to see that such a triple $(Y, i_X, e_X)$  if it exists, is unique
up to a unique isomorphism.   For any rigid $X$ in ${\Cal C}$  we
denote by $[X]^* \in \overline{\Cal C}$  the isomorphism class of objects $Y$
as in Definition 1.1.4.

\proclaim{{\bf 1.1.5}\ Definition}\ We say that the category ${\Cal C}$  is
rigid if
all its objects are rigid and for any $[Y]$ in $\overline{\Cal C}$  there
exists
$X$ in ${\Cal C}$  such that $[X]^* = [Y]$.
\endproclaim

If  ${\Cal C}$  is a rigid category, then there exists an equivalence
$*: {\Cal C} 	\to {\Cal C}^{op}$ of categories such that
$X^*$  is dual to $X$ for all $X$ in ${\Cal C}$,
$\alpha \in \text{Hom}_{{\Cal C}}(X,Y)$  the morphism
$\alpha^* \in \text{Hom}_{{\Cal C}}(Y^*,X^*)$  is the composition
$$
Y^*\longrightarrow Y^* \otimes \ident
{\buildrel {\text{id}\otimes i_X}\over\longrightarrow}
Y^* \otimes X \otimes X^*
{\buildrel {\text{id}\otimes\alpha\otimes \text{id}}\over\longrightarrow}
Y^* \otimes Y \otimes X^*
{\buildrel {e_Y \otimes \text{id}}\over\longrightarrow}
\ident \otimes X^* = X^*.
$$
Moreover, such an equivalence is unique up to a unique isomorphism.

We will fix such an equivalence between the categories
${\Cal C}$ and  ${\Cal C}^{op}$.  Then the functor $X \longrightarrow  X^{**}$
is an auto-equivalence
of ${\Cal C}$.

\subheading{1.1.6.\ Example}  Let  $(H,m,\Delta,i,\varepsilon,S)$  be a Hopf
algebra over  a ring $F$,  ${\Cal C}_H$  be the category of $H$-modules
$X = (\rho_X,\underline{X})$ and  $\otimes: {\Cal C}_H \times {\Cal C}_H \to
{\Cal C}_H$  be
the tensor product over $F$.  That is,  $X \otimes Y = ((\rho_X \otimes \rho_Y)
\circ \Delta, \underline{X} \otimes \underline{Y})$, and $\ident\
{\buildrel {\text{def}}\over =}\ (F,\epsilon)$.  Then $({\Cal
C}_H,\otimes,\ident)$
is a strict monoidal category.  In this case, $F_{\Cal C} = F$.

Let  $H_0 \subset H$  be the kernel of the counit $\varepsilon$.  For any
$X = (\rho_X,\underline{X})$  in ${\Cal C}_H$  we define
${X}^H_{(0)}\ {\buildrel {\text{def}}\over =}\ \rho_X(H_0)\underline{X}$.
Often we will write  $X_{(0)}$  instead of  $X^H_{(0)}$.  In
the case when $F$ is a field, we can identify  $\langle X\rangle$  (see
Definition 1.1.3)
with the quotient $\langle X\rangle\ {\buildrel {\text{def}}\over =}\
\underline{X}/{X}_{(0)}$.

Assume that $F$ is a field.  An object  $X = (\rho_X,\underline{X})$
is rigid if dim$_F\underline{X}
< \infty$.  In this case, $X^* = (\rho_{X^*}, \underline{X}^*)$ where
$\underline{X}^*\ {\buildrel {\text{def}}\over =}\ \underline{X}^\vee(=
\text{Hom}_F(\underline{X},F))$,
and $\rho_{X^*}(h) = (\rho_X(S(h)))^*$, where for any
$\alpha \in \text{End}\ \underline{X}$  we denote by   $\alpha^\vee$  the
endomorphism of  $\underline{X}^\vee$  dual to $\alpha$.  In this case, the
morphism  $e_X: X^* \otimes X \to \ident$ is induced by the natural pairing
$\underline{X}^* \otimes \underline{X} \to F$.  This pairing defines a
canonical isomorphism of the linear space $\underline{X}^{**}$ with
$\underline{X}$  and the action  $\rho_{X^{**}}$  of $H$ on
$\underline{X}^{**} = \underline{X}$  is given by the rule
$\rho_{X^{**}}(h) = \rho(S^2(h))$.

\proclaim{{\bf 1.1.7}\ Proposition}  For any $X$ in ${\Cal C}$  and
$\varphi \in Hom(X,X)$  the diagrams
$$
{\matrix \format \c &\quad \c &\quad \c \\
\ident & {\buildrel{i_X}\over\llongrightarrow} & X \otimes X^* \\
\vspace{1\jot}
\bigg\downarrow i_X && \bigg\downarrow \varphi\otimes id \\
\vspace{1\jot}
X \otimes X^* &  {\buildrel {id \otimes \varphi^*} \over\llongrightarrow}
& X \otimes X^*
\endmatrix}\qquad\qquad\qquad
{\matrix \format \c &\quad \c &\quad \c \\
X^*\otimes X & {\buildrel{\varphi^*\otimes id}\over\llongrightarrow} &
X^* \otimes X \\
\vspace{1\jot}
\bigg\downarrow 1\otimes\varphi && \bigg\downarrow e_X \\
\vspace{1\jot}
X^* \otimes X &  {\buildrel {e_X} \over\llongrightarrow}
& \ident
\endmatrix}
$$
are commutative.
\endproclaim

\demo{Proof}  We prove the commutativity of the first diagram.  The proof
of the commutativity of the second diagram is completely analogous.

Let
$$
a\ {\buildrel {\text{def}}\over =}\  (\varphi \otimes id_{X^*}) \circ
i_X,\ \ b\ {\buildrel {\text{def}}\over =}\  (id_{X} \otimes \varphi^*) \circ
i_X \in Hom_{{\Cal C}}(\ident, X \otimes X^*).
$$
By the definition of  $i_X$ and  $e_X$,   $a$ is equal to the
composition
$$
\ident\ {\buildrel{i_X}\over\longrightarrow}\ X \otimes X^*\
{\buildrel {\varphi \otimes id_{X^*}}\over\llongrightarrow}\ X \otimes X^*\
{\buildrel {i_X\otimes id_{X\otimes X^*}}\over\llongrightarrow}\
X \otimes X^* \otimes X \otimes X^*\
{\buildrel {id_X \otimes e_X \otimes id_{X^*}}\over\llongrightarrow}\
X \otimes X^*.
$$
But the composition  $(i_X \otimes id_{X\otimes X^*}) \circ
(\varphi \otimes id_{X^*})$  is equal to the composition

\noindent
$(id_{X\otimes X^*} \otimes \varphi \otimes id_{X^*}) \circ
(i_X \otimes id_{X\otimes X^*})$.  Since the compositions

\noindent
$\ident {\dsize\buildrel {i_X}\over\longrightarrow}\ X \otimes X^*\
{\dsize \buildrel {id_{X\otimes X^*} \otimes i_X}\over\llongrightarrow}\
X \otimes X^* \otimes~X~\otimes~X^*$ and
$\ident {\dsize\buildrel {i_X}\over\longrightarrow}\ X \otimes X^*\
{\dsize \buildrel {i_X \otimes id_{X\otimes X^*}}\over\llongrightarrow}\
X \otimes X^* \otimes X \otimes X^*$  are equal (and both coincide with
the composition $\ident = \ident \otimes \ident\
{\dsize\buildrel {i_X \otimes i_X}\over\longrightarrow}\ X \otimes X^* \otimes
X \otimes X^*)$  we see that  $\alpha$ is equal to the composition
$$
\ident {\dsize\buildrel {i_X}\over\longrightarrow}\ X \otimes X^*\
{\dsize \buildrel {id_{X\otimes X^*} \otimes i_X}\over\llongrightarrow}\
X \otimes X^* \otimes X \otimes X^*\
{\buildrel{id_{X\otimes X^*}\otimes\varphi\otimes
id_{X^*}}\over\llongrightarrow}\
X \otimes X^* \otimes X \otimes X^*\
{\buildrel {id_X \otimes e_X \otimes id_{X^*}}\over\llongrightarrow}\
X \otimes X^*.
$$
But the composition of the last three arrows is equal to $b = id_X \otimes
\varphi^*$. Lemma 1.1.7 is proved.
\enddemo

\subheading{1.1.8}   For any  $X,Y$ in ${\Cal C}$  the morphisms  $i,e$
defined as the compositions
$$
\aligned
&i: \ident {\buildrel {i_X}\over\longrightarrow} X \otimes X^* =
X \otimes \ident \otimes X^*
{\buildrel {\text{id}\otimes i_Y \otimes \text{id}}\over\llongrightarrow}
(X \otimes Y) \otimes (Y^* \otimes X^*)\\
& e: (Y^* \otimes X^*) \otimes (X \otimes Y)
{\buildrel {\text{id}\otimes e_X \otimes \text{id}}\over\llongrightarrow}
Y^* \otimes \ident \otimes Y = Y^* \otimes Y
{\buildrel {e_Y}\over\longrightarrow} \ident
\endaligned
$$
satisfy the conditions of Definition 1.1.4.  Therefore they define isomorphisms
of  $(X \otimes Y)^*$  with  $Y^* \otimes X^*$ such that $i = i_{X \otimes Y}$
and $e = e_{X \otimes Y}$.  We will
freely use this identification and will therefore identify the second dual
$(X \otimes Y)^{**}$  with  $X^{**} \otimes Y^{**}$.

\demo{Remark}  As follows from ([DM], 1.17) the rigidity of ${\Cal C}$
implies the semisimplicity of $\ident$.  We will always assume that $\ident$
is simple.  In this case  $F_{\Cal C}$  is a field.
\enddemo

\proclaim{{\bf 1.1.9.}\ Definition} Let  ${\Cal C}$  be a monoidal category,
${\Cal D} \subset {\Cal C}$  a full monoidal subcategory such that any
object of ${\Cal D}$  is rigid.

For any  $V,W$  in ${\Cal C}$  and $X,Y$  in ${\Cal D}$
we define the maps  $\alpha_V^{X}:  \langle V \otimes X\rangle^\vee \to
\langle X^{**} \otimes V\rangle^\vee$  and
$$
\varphi^{X,Y}_{V,W}: \text{Hom}_{{\Cal C}}(V \otimes X, Y \otimes W)
\to \text{Hom}_{{\Cal C}}(Y^* \otimes V, W \otimes X^*)
$$
as the compositions
$$
\aligned
\alpha^X_V: \langle V\otimes X\rangle^\vee &= Hom_{{\Cal C}}(V \otimes
X,\ident)\
\longrightarrow\ Hom_{{\Cal C}}(X^{**} \otimes V \otimes X \otimes
X^*, X^{**} \otimes X^*)\\
&\qquad {\buildrel{id_{X^{**}\otimes V}\otimes i_X \otimes
e_{X^*}}\over\llongrightarrow}\ Hom_{{\Cal C}}(X^{**} \otimes V,\ident) \\
&= \langle X^{**} \otimes V\rangle^\vee.
\endaligned
$$

$$
\varphi^{X,Y}_{V,W}(a):  Y^* \otimes V\
{\buildrel {id \otimes i_X}\over\longrightarrow}\ Y^* \otimes V \otimes X
\otimes X^*\ {\buildrel {id \otimes a \otimes id}\over\longrightarrow}\
Y^* \otimes Y \otimes W \otimes X^*\ {\buildrel {e_Y}\over\longrightarrow}\
W \otimes X^*
$$
for all $a \in \text{Hom}_{{\Cal C}}(V \otimes X, Y \otimes W)$.
\endproclaim

\demo{Remark}  We will write  $\varphi_V^{X,Y}$  instead of
$\varphi^{X,Y}_{V,V}$.
\enddemo

\proclaim{Lemma}  If  $F_{\Cal C}$  is a field, then there exist
$F_{\Cal C}$-linear maps $\beta^V_{X^{**}}: \langle X^{**} \otimes V\rangle
\to \langle V \otimes X\rangle$  such that $\alpha^X_V =
(\beta^V_{X^{**}})^\vee$.
\endproclaim

\demo{Proof}  For any $R$ in Vec$_{\Cal C}$  we can define a $F_{\Cal
C}$-linear
map  $\alpha^X_V(R):  Hom(V \otimes X, R_{\Cal C}) \to Hom(X^{**} \otimes V,
R_{\Cal C})$  exactly in the same way as we have defined the map
$\alpha^X_V$.   Then maps $\alpha^X_V(R)$  define a morphism from the functor
$L_{V \otimes X}$  to the functor  $L_{X^{**} \otimes V}$.  The corresponding
morphism
between the representing objects is  $\beta^V_X$.
\enddemo

\subheading{1.1.10} For any $V$ in ${\Cal C}$  we denote by $E_V$  the ring of
endomorphisms of $V$ and by  $E_V^{op}$  the opposite ring.  If
$V$ is rigid, the map $f \to f^*$  defines an isomorphism of $E_V$
with  $E^{op}_{V^*}$.

\proclaim{Lemma}  a)\ The linear maps  $\alpha_V^X$, $\beta^V_{X^{**}}$  and
$\varphi_{V,W}^{X,Y}$  are isomorphisms.

b)  For any $V$ in ${\Cal C}$, $X,Y$ in
${\Cal D}$  the linear map  $\alpha_V^{(X\otimes Y)}$  is equal
to the composition  $\alpha_{Y^{**}\otimes V}^{X} \circ
\alpha_{V\otimes X}^{Y}$  and the map  $\beta^V_{(X\otimes Y)^{**}}$  is
equal to the composition  $\beta_{Y^{**}}^{V\otimes X} \circ
\beta_{X^{**}}^{Y^{**}\otimes V}$.

c)  For any $f_X \in E_X,\ f_Y \in E_Y$, $f_W \in E_W$,$f_V\in E_V$  and
$a \in \text{Hom}_{{\Cal C}}(V \otimes X, Y\otimes W)$  we have
$$
\varphi^{X,Y}_{V,W}((f_Y \otimes f_W) \circ a \circ (f_V \otimes f_X)) =
(f_W \otimes f^*_X) \circ \varphi^{X,Y}_{V,W}(a) \circ (f^*_Y \otimes
f_V).
$$
\endproclaim

\demo{Proof}  a)  Consider the map from  $\langle X^{**} \otimes V\rangle^\vee$
to  $\langle V \otimes X\rangle^\vee$ defined as the composition
$$
\aligned
\langle X^{**} \otimes V\rangle^\vee &=
\text{Hom}_{{\Cal C}}(X^{**} \otimes V, \ident) \longrightarrow
\text{Hom}_{{\Cal C}}(X^* \otimes X^{**} \otimes
V \otimes X, X^* \otimes X)\\
&\qquad {\buildrel {i_{X^*} \otimes id_{V\otimes X}
\otimes e_X}\over\llongrightarrow}\
\text{Hom}_{{\Cal C}}(V \otimes X, \ident)  \\
&= \langle V \otimes X\rangle^\vee.
\endaligned
$$
As follows immediately from Definition 1.1.4, this map is the inverse to
$\alpha_V^{X}$.  The construction of the inverse to  $\varphi_{V,W}^{X,Y}$  is
completely analogous.

b)  Follows immediately from the definitions.

c)  We have to show that for any $a \in Hom_{{\Cal C}}(V \otimes X,
X \otimes V)$  the following equalities hold:

\itemitem{i)}  $\varphi^{X,Y}_{V,W}(a \circ (f_V \otimes id_X)) =
\varphi^{X,Y}_{V,W} (a)(id_{Y^{*}} \otimes f_V)$  for any  $f_V \in E_V$

\itemitem{ii)}  $\varphi_{V,W}^{X,Y}(a \circ (id_W \otimes f_X)) = (id_W
\otimes
f^*_X) \circ \varphi_{V,W}^{X,Y}(a)$ for all  $f_X \in E_X$.

\itemitem{iii)} $\varphi^{X,Y}_{V,W}((id_Y\otimes f_W) \circ a) =
(f_W \otimes id_{X^*}) \circ \varphi^{X,Y}_{V,W}(a)$

\itemitem{iv)}  $\varphi^{X,Y}_{V,W}((f_Y \otimes id_W) \circ a) =
\varphi^{X,Y}_{V,W}(a) \circ (f^*_Y \otimes id_V)$.

\medpagebreak

The proofs of i) and iii) are straightforward.  We will show how to prove ii).
The proof of (iv) is completely analogous.

By the definition
the map  $\varphi_{V,W}^{X,Y}(a \circ (id_V \otimes f_X))$  is equal to the
 composition
$$
\aligned
& Y^* \otimes V\ {\buildrel {\ident \otimes i_X}\over\longrightarrow}\
Y^* \otimes V \otimes X \otimes X^*\
{\buildrel {id\otimes f_X \otimes id}\over\llongrightarrow}\ Y^* \otimes V
\otimes X \otimes X^* \\
&\quad {\buildrel {id \otimes a \otimes id}\over\llongrightarrow}\
Y^* \otimes Y \otimes W \otimes X^* \
{\buildrel {e_Y \otimes id}\over\llongrightarrow}\ W \otimes X^*.
\endaligned
$$
As follows from Lemma 1.1.7 a) this composition is equal to the composition
$$
\aligned
& Y^* \otimes V\ {\buildrel {\ident \otimes i_X}\over\longrightarrow}\
Y^* \otimes V \otimes X \otimes X^*\
{\buildrel {id\otimes id_X\otimes f^*_X}\over\llongrightarrow}\ Y^* \otimes V
\otimes X \otimes X^* \\
&\quad {\buildrel {id \otimes a \otimes id}\over\llongrightarrow}\
Y^* \otimes Y \otimes W \otimes X^*
{\buildrel {e_Y \otimes id}\over\longrightarrow}\ W \otimes X^*.
\endaligned
$$
But the last composition is equal to $(id_W \otimes f^*_X) \circ
\varphi_{V,W}^{X,Y}(a)$.

Lemma 1.1.10 is proved.
\enddemo

\subheading{1.1.11}  Let $H$ be a Hopf algebra, ${\Cal C} = {\Cal C}_H$,
$X = (\rho_X,\underline{X}), Y = (\rho_V,\underline{V})$ be $H$-modules
such that  dim$_{\Bbb C} \underline{X} < \infty$.

Let $P^X_V: \underline{X} \otimes \underline{V} \to \underline{V} \otimes
\underline{X}$  be the permutation  $P^X_V(x\otimes v) = v \otimes x$.  We
can consider  $P^X_V$  as a linear map from  $\underline{X}^{**} \otimes
\underline{V}$  to  $\underline{V} \otimes \underline{X}$.

\proclaim{{\bf 1.1.12}\ Lemma}  $P^X_V$  maps the subspace
$({X}^{**} \otimes {V})_{(0)}$  into
$({V} \otimes {X})_{(0)}$  and induces
the linear map  from  $\langle X^{**} \otimes V\rangle$  to
$\langle V \otimes X\rangle$ equal to   $\beta^V_{X^{**}}$.
\endproclaim

\demo{Proof}  Let  $(P^X_V)^*: (\underline{V} \otimes \underline{X})^* \to
(\underline{X}^{**} \otimes \underline{V})^*$  be the linear map dual to
$P^X_V$.  It is sufficient to show that for any  $\lambda \in
\langle V \otimes X\rangle^\vee \subset (\underline{V} \otimes
\underline{X})^*$
we have  $(P^X_V)^*(\lambda) = (\alpha_V^{X})(\lambda)$.
\enddemo

Consider first the case when  $V = X^*$  and  $\lambda = e_X$.  It is easy to
check that  $\alpha^{X}_{X^*}(e_X)~=~e_{X^*}$  and the validity of Lemma 1.1.12
follows from definition of  $e_X$.

For any $\lambda \in \langle V \otimes X\rangle^\vee =
\text{Hom}_{{\Cal C}_H}(V \otimes X, \ident)$  we have $\lambda =
e_X \circ (\widetilde{\lambda} \otimes id_X)$ where

\noindent
$\widetilde{\lambda} \in \text{Hom}_{{\Cal C}_H}(V,X^*)$  and
Lemma 1.1.12 follows from the functoriality of  $\alpha^{X}_V$  and $P^X_V$.
Lemma 1.1.12 is proved.

\proclaim{{\bf 1.1.13}\ Definition}  A strict endomorphism ${\Cal T}$ of ${\Cal
C}$
is a functor from the category ${\Cal C}$  to itself such that

\itemitem{a)}  ${\Cal T}(\ident) = \ident$,

\itemitem{b)}  ${\Cal T}(X \otimes Y) = {\Cal T}(X) \otimes {\Cal T}(Y)$  for
all  $X,Y$  in ${\Cal C}$  and those identifications are compatible
with morphisms in ${\Cal C}$, and

\itemitem{c)}  for any rigid object $X$ in ${\Cal C}$  we have
${\Cal T}(X^*) = ({\Cal T}(X))^*$,

\itemitem{d)}  ${\Cal T}(i_X) = i_{{\Cal T}(X)}, {\Cal T}(e_X)
= e_{{\Cal T}(X)}$.
\endproclaim

We say that  ${\Cal T}$  is a strict automorphism if it is an equivalence of
categories.

If  ${\Cal T}$  is a strict automorphism of ${\Cal C}$, then for any $X$ in
${\Cal C}$,\ ${\Cal T}$  defines an isomorphism from Hom$_{\Cal C}(X,\ident) =
\langle X\rangle^\vee$  to  Hom$_{\Cal C}({\Cal T}(X),\ident) =
\langle {\Cal T}(X)\rangle^\vee$.  It is easy to see that this isomorphism
$\langle X\rangle^\vee {\buildrel {\sim}\over\longrightarrow}\
\langle {\Cal T}(X)\rangle^\vee$  comes from an isomorphism
$\langle X\rangle {\buildrel {\sim}\over\longrightarrow}\
\langle {\Cal T}(X)\rangle$.

\demo{{\bf 1.1.14}\ Example}  For any rigid monoidal category ${\Cal D}$ the
functor
$X \to X^{**}$  is a strict automorphism.  We denote the inverse to this
automorphism by  ${\Cal T}_*$.  It
is clear that any strict automorphism ${\Cal T}$  of ${\Cal D}$  commutes with
${\Cal T}_*$.
\enddemo

\demo{{\bf 1.1.15}\ Remark}  We can interpret the linear map $\beta^V_X$
defined
in Lemma 1.1.9 as the linear map from $\langle X \otimes V\rangle$
to $\langle V \otimes {\Cal T}_*(X)\rangle$.
\enddemo

\subheading{1.1.16}  Let  ${\Cal C}$  be a monoidal category, ${\Cal D}$
and ${\Cal O}$  full subcategories  of ${\Cal C}$ such that ${\Cal D}$  is
rigid and  ${\Cal T}$  a strict automorphism of ${\Cal D}$.

\proclaim{Definition}   A right ${\Cal T}$-braiding of ${\Cal D}$ on ${\Cal O}$
 is
a functorial isomorphism
$$
s_{V,X}: V \otimes X \to {\Cal T}(X) \otimes V,
$$
defined for all $V$ in ${\Cal O}$  and $X$ in ${\Cal D}$  such that
$s_{V,\ident} = id_V$  and for all $X,Y$  in ${\Cal D}$  we have
$$
s_{V,X\otimes Y} = (id_{{\Cal T}(X)} \otimes s_{V,Y}) \circ (s_{V,X} \otimes
id_Y).
$$
\endproclaim

\proclaim{{\bf 1.1.17}\ Proposition}  For any  $X$ in ${\Cal D}$  and $V$ in
${\Cal O}$  the diagram
$$
\matrix \format \c &\quad \c &\quad \c \\
{\Cal T}(X)^* \otimes V \otimes X &
{\buildrel{id_{{\Cal T}(X)} \otimes s_{V,X}}\over\llongrightarrow} &
{\Cal T}(X^*) \otimes {\Cal T}(X) \otimes V \\
\vspace{1\jot}
\Bigg\downarrow s^{-1}_{V,X^*} \otimes id_X && \Bigg\downarrow e_{{\Cal T}(X)}
\\
\vspace{1\jot}
V \otimes X^* \otimes X{\phantom{T(X)}}  & {\buildrel
{e_X}\over\llongrightarrow}
& V{\phantom{T(X)}}
\endmatrix
$$
is commutative.
\endproclaim

\demo{Proof}  As follows from the functoriality of  $s_{V,X}$  and the equality
$s_{V,\ident} = id$, the diagram
$$
\matrix \format \r &\quad \c &\quad \l \\
V \otimes (X^* \otimes X) & {\buildrel {s_{V,X^* \otimes
X}}\over\llongrightarrow}
& (X^* \otimes X) \otimes V \\
\vspace{1\jot}
\bigg\downarrow e_X && e_X\bigg\downarrow \\
\vspace{1\jot}
V& = & V
\endmatrix
$$
is commutative.  On the other hand, we know that
$$
s_{V,X^* \otimes X} = (id_{{\Cal T}(X^*)} \otimes s_{V,X}) \circ
(s_{V,X^*} \otimes id_X).
$$
Therefore the diagram
$$
\matrix \format \c &\quad \c &\quad \c \\
{\Cal T}(X)^* \otimes V \otimes X &
{\buildrel{id_{{\Cal T}(X^*)} \otimes s_{V,X}}\over\llongrightarrow} &
{\Cal T}(X^*) \otimes {\Cal T}(X) \otimes V \\
\vspace{1\jot}
\Bigg\uparrow s_{V,X^*} \otimes id_X && \Bigg\downarrow e_{{\Cal T}(X)} \\
\vspace{1\jot}
V \otimes X^* \otimes X{\phantom{T(X)}}  & {\buildrel
{e_X}\over\llongrightarrow}
& V{\phantom{T(X)}}
\endmatrix
$$
is commutative.  Proposition 1.1.17 is proved.
\enddemo

\proclaim{{\bf 1.1.18}\ Proposition}  For any $V$ in ${\Cal O}$  and any
$X$ in ${\Cal D}$
$$
\varphi^{X,{\Cal T}(X)}_V (s_{V,X}) = s^{-1}_{V,X^*}.
$$
\endproclaim

\demo{Proof}  By the definition the map $\varphi_V^{X,{\Cal T}(X)}(s_{V,X}) \in
Hom_{{\Cal C}}({\Cal T}(X)^* \otimes V, V \otimes X^*)$  is defined as
the composition
$$
{\Cal T}(X)^* \otimes V\ {\buildrel {id\otimes i_X}\over\llongrightarrow}
{\Cal T}(X)^* \otimes V \otimes X \otimes X^*
{\buildrel{s_{V,X}}\over\longrightarrow}
{\Cal T}(X)^* \otimes {\Cal T}(X) \otimes V \otimes X^*\
{\buildrel {id \otimes e_{{\Cal T}(X)} \otimes id}\over\llongrightarrow}
V \otimes X^*.
$$
As follows from Lemma 1.1.17, this  composition is equal to the composition
$$
{\Cal T}(X)^* \otimes V\ {\buildrel {id\otimes i_X}\over\llongrightarrow}
{\Cal T}(X)^* \otimes V \otimes X \otimes X^*
{\buildrel{s^{-1}_{V,X^*}\otimes id}\over\llongrightarrow}
V \otimes X^* \otimes X \otimes X^*\
{\buildrel {id_V \otimes e_X \otimes id_{X^*}}\over\llongrightarrow}
V \otimes X^*.
$$
As follows from the definition of rigidity, the composition of the last two
morphisms is equal to  $s^{-1}_{V,X^*}$.  Proposition 1.1.18 is proved.
\enddemo

\proclaim{{\bf 1.1.19}\ Definition}  Let ${\Cal C}$
be a strict monoidal category.  A braiding $s$  on ${\Cal C}$  is
a functorial system of isomorphisms
$s_{X,Y} \in \text{Isom}(X \otimes Y, Y \otimes X)$  for
$X, Y$ in ${\Cal C}$	 such that

\itemitem{a)}  for any $X, Y, Z$ in ${\Cal C}$  we have
$$
\aligned
s_{X\otimes Y,Z} &= (s_{X,Z} \otimes id_Y) \circ (id_X \otimes s_{Y,Z}), \\
s_{X,Y\otimes Z} &= (id_X \otimes s_{X,Z}) \circ (s_{X,Y} \otimes id_Z);
\endaligned
$$

\itemitem{b)}  $s_{X,\ident} = s_{\ident,X} = id_X$ for all $X$ in
${\Cal C}$.
\endproclaim
\vskip .3in

%this file is qkz1b.tex	  [contains 1.2* and 1.3*]

\subheading{1.2 $KZ$-data}

\proclaim{{\bf 1.2.1}\ Definition}  Let  ${\Cal D}$  be a rigid monoidal
category.
A weak braiding $b = ({\Cal D}^{(2)}, s)$ on ${\Cal D}$  is

\item{1)}\ A choice of a subset  $\overline{\Cal D}^{(2)}$  of
$\overline{\Cal D} \times \overline{\Cal D}$  such that

\itemitem{a)}\ 	For any  $X,Y,Z$ in ${\Cal D}$  such that
$([X],[Y]), ([X],[Z])$ are in  $\overline{\Cal D}^{(2)}$  the pair
\newline
$([X],[Y \otimes Z])$ is in $\overline{\Cal D}^{(2)}$.

\itemitem{b)}\ For any $X,Y,Z$ in ${\Cal D}$  such that  $([X],[Z])$  and
$([Y],[Z])$  in  $\overline{\Cal D}^{(2)}$  the pair
\newline
$([X \otimes Y], [Z])$  is in $\overline{\Cal D}^{(2)}$.

\itemitem{c)}\ $(\overline{X},\ident) \in \overline{\Cal D}^{(2)}$ and
$(\ident,\overline{X}) \in  \overline{\Cal D}^{(2)}$,
for all $\overline{X} \in \overline{\Cal D}$.

\item{2)}\ A functorial isomorphism $s_{X,Y}$  between the restrictions of
functors
\newline
$(X,Y) \to X~\otimes~Y, (X,Y) \to (Y \otimes X)$  on
${\Cal D}^{(2)}$  such that

\itemitem{a)} $s_{\ident,X} = s_{X,\ident} = \text{id}$  for all $X \in
{\Cal D}$

\itemitem{b)} $s_{X,Y\otimes Z} = (1_Y \otimes s_{X,Z}) \circ (s_{X,Y} \otimes
1_Z)$  for all $X,Y,Z \in {\Cal D}$  such that
\newline
$(X,Y), (X,Z) \in {\Cal D}^{(2)}$

\itemitem{c)} $s_{X\otimes Y,Z} = (s_{X,Z} \otimes 1_Y) \circ (1_X \otimes
s_{Y,Z})$  for all $X,Y,Z \in {\Cal D}$  such that
\newline
$(X,Z), (Y,Z) \in {\Cal D}^{(2)}$,
\newline
where we denote by  ${\Cal D}^{(2)}$  the
full subcategory of  ${\Cal D} \times {\Cal D}$  of pairs $(X,Y)$  such that
$([X],[Y]) \in \overline{\Cal D}^{(2)}$.
\endproclaim

\proclaim{{\bf 1.2.2}\ Definition}  A $KZ$-data consists of a monoidal category
${\Cal C}$, its full subcategories ${\Cal D}, {\Cal C}^\pm$,
strict automorphisms\   ${\Cal T}_\pm$ of ${\Cal D}$,
right ${\Cal T}_\pm$ braidings
$s_\pm$  of ${\Cal D}$  on  ${\Cal C}^\pm$,
and a weak braiding  $({\Cal D}^{(2)},s)$  compatible with the
automorphisms  ${\Cal T}_\pm$.
\endproclaim

\subheading{1.2.3}  We say that  $KZ$-data $({\Cal C}, {\Cal D}, {\Cal C}^\pm,
s_\pm, {\Cal T}_\pm, {\Cal D}^{(2)},s)$  is rigid if ${\Cal D}$  is
rigid.  In this case we denote by ${\Cal T}$  the automorphism of ${\Cal D}$
which is the composition  ${\Cal T} = {\Cal T}_- {\Cal T}_* {\Cal T}_+$.

\proclaim{{\bf 1.2.4}\ Definition}  Let ${\Cal K} = ({\Cal C}, {\Cal D}, {\Cal
C}^\pm,
s_\pm, {\Cal T}_\pm, {\Cal D}^{(2)},s)$  be a rigid  $KZ$-data.

\item{a)}  We say that a pair $(X,Y) \in
{\Cal D} \times {\Cal D}$  is ${\Cal T}$-generic
(or simply generic) if for any
$r \in {\Bbb Z}$, $({\Cal T}^r(X),Y) \in {\Cal D}^{(2)}$.

\item{b)}  For any $n \in {\Bbb Z}$  we denote by  $S_n = S^{\Cal T}_n$  the
set of
$n$-tuples $(X_1,\ldots,X_n) \in {\Cal D}^n$  such that all the pairs
$(X_i,X_j)$, $1 \le i \ne j \le n$  are generic.

\item{c)}  For any $i$,\ $1 \le i \le n-1$  we denote by  $p_i: S^{\Cal T}_n
\to
S^{\Cal T}_{n-1}$  the map
$$
(X_1,\ldots,X_n) \to (X_1,\ldots,X_{i-1}, X_i \otimes X_{i+1},
X_{i+2},\ldots,X_n).
$$

\item{d)} We define an action of the group  ${\Bbb Z}^n$  on $S^{\Cal T}_n$  by
the rule
$$
{\Cal E}_i(X_1,\ldots,X_n) = (X_1,\ldots,X_{i-1}, {\Cal T}(X_i),
X_{i+1},\ldots,X_n),
$$
where   $\{ {\Cal E}_i\},\ 1 \le i \le n$,  is the standard
set of generators of the group  ${\Bbb Z}^n$.

\item{e)}  For any  pair  $V \in {\Cal C}^+, W \in {\Cal C}^-$  and a point
$x = (X_1,\ldots,X_n) \in S_n$  we define a vector space
${\Cal F}^{(n)}_{V,W}(x)\ {\dsize\buildrel {\dsize\text{def}}\over =}\
\langle V \otimes X_1 \otimes \cdots \otimes X_n \otimes W\rangle$.  We will
consider  ${\Cal F}^{(n)}_{V,W}$  as a set-theoretical vector bundle over
$S_n$.
\endproclaim

It is clear that for any  $i$,\ $1 \le i \le n-1$  we have a
canonical isomorphism
$$
{\Cal F}^{(n)}_{V,W}\ {\buildrel \sim\over\longrightarrow}\
p^*_i({\Cal F}^{(n-1)}_{V,W}).
$$
\enddemo

\proclaim{{\bf 1.2.5}\ Definition}  We denote by  $\theta_n: S_n \to S_n$  an
automorphism
given by the rule
$$
\theta_n(X_1,\ldots,X_n) = (X_2,\ldots,X_n, {\Cal T}(X_1))
$$
and by   $\widehat{\theta}_n$  its lifting  $\widehat{\theta}_n:
{\Cal F}^{(n)}_{V,W} \to \theta^*_n({\Cal F}^{(n)}_{V,W})$  to the bundle
defined as the composition
$$
\aligned
& \widehat{\theta}_n: {\Cal F}^{(n)}_{V,W}(x) = \langle V \otimes X_1 \otimes
\cdots \otimes X_n \otimes W)\rangle = \langle V \otimes X_1 \otimes Y \otimes
W\rangle \to \\
&{\buildrel {s^+_{V,X_1}}\over\longrightarrow} \langle {\Cal T}^+(X_1) \otimes
V
\otimes Y \otimes W\rangle\
{\buildrel{\beta_{{\Cal T}_+(X_1)}^{V\otimes Y \otimes
W}}\over\longrightarrow}\
\langle V \otimes Y \otimes W \otimes {\Cal T}_* {\Cal T}_+ (X_1) \rangle \to
\\
& {\buildrel {s^-_{W,{\Cal T}_*{\Cal T}_+(X_1)}}\over\longrightarrow}\ \langle
V
\otimes Y \otimes {\Cal T}(X_1) \otimes W\rangle =
{\Cal F}^{(n)}_{V,W}(\theta_n(x))
\endaligned
$$
for all $x = (X_1,\ldots,X_n) \in S_n$  where  $Y = X_2 \otimes \cdots \otimes
X_n$  and the isomorphism

\noindent
$\beta^U_Z: \langle Z \otimes U\rangle\
{\dsize\buildrel{\sim}\over\longrightarrow}\ \langle U \otimes {\Cal
T}_*(Z)\rangle$
is defined in 1.1.15.
\endproclaim

It is clear that  $p_{n-1} \circ \theta^2_n = \theta_{n-1} \circ p_{n-1}$.
Since  ${\Cal F}^{(n)}_{V,W} = p^*_{n-1}({\Cal F}^{(n-1)}_{V,W})$  we can
consider
$\widehat{\theta}^2_n$  and  $p^*_{n-1}(\widehat{\theta}_{n-1})$  as
automorphisms
of ${\Cal F}^{(n)}_{V,W}$  lifting the transformation  $\theta^2_n$ of $S_n$.
Analogously for any $i$,\ $1 \le i < n-1$  we can consider
$p^*_i(\widehat{\theta}_{n-1})$  as an automorphism of  ${\Cal F}^{(n)}_{V,W}$
over  $\theta_n$.

\proclaim{{\bf 1.2.6}\ Proposition}

\itemitem{a)}\ $\widehat{\theta}_n$  commutes with endomorphisms of
$Y = X_1 \otimes \cdots \otimes X_{n-1}$  and for any endomorphism  $a$ of
$X_n$  we have  $(\theta_n \circ a) = {\Cal T}(a) \circ \theta_n$.

\itemitem{b)}\ $\widehat{\theta}^2_n = p^*_{n-1}(\widehat{\theta}_{n-1})$.

\itemitem{c)}\ $\widehat{\theta}_n = p^*_i(\widehat{\theta}_{n-1})$  for all
$i$,\ $1 \le i < n-1$.
\endproclaim

\demo{Proof}  Follows immediately from the definition and the properties of
$s_\pm$ and $s$.
\enddemo

\proclaim{{\bf 1.2.7}\ Definition}  For any $i$,\ $1 \le i \le n$  we denote by
$\delta_i^{(n)}: {\Cal F}^{(n)}_{V,W}(x) \to {\Cal F}^{(n)}_{V,W}({\Cal
E}_i(x))$
a lifting of  ${\Cal E}_i$  to ${\Cal F}^{(n)}_{V,W}$  defined by a composition
$$
\aligned
& \delta_i^{(n)}: {\Cal F}^{(n)}_{V,W}(x) = \langle V \otimes X^+ \otimes
X_i \otimes X^- \otimes W\rangle \longrightarrow\\
&{\buildrel {s^{-1}_{X_i,X^+}}\over\longrightarrow}\ \langle V \otimes X_i
\otimes X^+ \otimes X^- \otimes W\rangle\
{\buildrel {\widehat{\theta}_n}\over\longrightarrow}\
\langle V \otimes X^+ \otimes X^-  \otimes {\Cal T}(X_i) \otimes W\rangle\\
& {\buildrel {s_{X^-,{\Cal T}(X_i)}}\over\longrightarrow}\
\langle V \otimes X^+ \otimes {\Cal T}(X_i) \otimes X^- \otimes W\rangle =
{\Cal F}^{(n)}_{V,W}({\Cal E}_i(x))
\endaligned
$$
for all $x = (X_1,\ldots,X_n) \in S_n$  where  $X^+ = X_1 \otimes \cdots
\otimes X_{i-1}, X^- = X_{i+1} \otimes \cdots \otimes X_n$.
\endproclaim

\proclaim{{\bf 1.2.8}\ Proposition}

\itemitem{a)}  $\delta_{i+1}^{(n)} \delta_{i}^{(n)} =
p^*_i(\delta_{i}^{(n-1)})$

\itemitem{b)}  $\delta_i^{(n)} \delta_{i+1}^{(n)} = p^*_i(\delta_i^{(n-1)})$

\itemitem{c)}  If $j < i$, then  $\delta_{i}^{(n)} =
p^*_j(\delta_{i-1}^{(n-1)})$

\itemitem{d)}  If $j > i+1$, then $\delta_{i}^{(n)} =
p^*_j(\delta_{i}^{(n-1)})$.
\endproclaim

\demo{Proof}  a)  For simplicity we consider the case  $i = 1$.  The proof in
the
general case is completely analogous.  Let  $x = (V_1,\ldots,V_n) \in S_n$.  We
define  $Y = X_3 \otimes \cdots \otimes X_n$.  To prove a) we consider the
composition
$$
\delta_2^{(n)}\delta_1^{(n)}(x):  \langle V \otimes X_1 \otimes X_2
\otimes Y \otimes W\rangle \to \langle V \otimes Y \otimes {\Cal T}(X_1)
\otimes {\Cal T}(X_2) \otimes Y \otimes W\rangle.
$$
Using the equality  $s_{X_2 \otimes Y,X_1} = (s_{X_2,X_1} \otimes id_Y)
\circ (id_{X_2} \otimes s_{Y,X_1})$   we can write this map as the
composition
$$
\aligned
& \langle V \otimes X_1 \otimes X_2 \otimes Y \otimes W\rangle\
{\buildrel {\widehat{\theta}_n}\over\longrightarrow}\ \langle V \otimes X_2
\otimes Y \otimes {\Cal T}(X_1) \otimes W\rangle
{\buildrel {s_{Y,{\Cal T}(X_1)}}\over\llongrightarrow}\\
& \langle V \otimes X_2 \otimes {\Cal T}(X_1) \otimes Y \otimes W\rangle\
{\buildrel {\widehat{\theta}_n}\over\longrightarrow}\ \langle V
\otimes {\Cal T}(X_1) \otimes Y \otimes {\Cal T}(X_2) \otimes W\rangle\
{\buildrel {s_{Y, {\Cal T}(X_2)}}\over\llongrightarrow}\\
& \longrightarrow \langle V \otimes {\Cal T}(X_1) \otimes {\Cal T}(X_2)
\otimes Y \otimes W\rangle.
\endaligned
$$
\enddemo

By  Proposition 1.2.6 a)  this composition is equal to the composition
$$
\aligned
& \langle V \otimes X_1 \otimes X_2 \otimes Y \otimes Y \otimes W\rangle\
{\buildrel{\widehat{\theta}_n}\over\longrightarrow}\
\langle V \otimes X_2 \otimes Y \otimes {\Cal T}(X_1) \otimes W\rangle\
{\buildrel{\widehat{\theta}_n}\over\longrightarrow} \\
& \langle V \otimes Y \otimes {\Cal T}(X_1) \otimes {\Cal T}(X_2) \otimes
W\rangle\
{\buildrel{s_{Y,{\Cal T}(X_1)}}\over\llongrightarrow}\
\langle V \otimes {\Cal T}(X_1) \otimes Y \otimes {\Cal T}(X_2) \otimes
W\rangle \longrightarrow \\
& {\buildrel {s_{Y,{\Cal T}(X_2)}}\over\longrightarrow}\ \langle V \otimes
{\Cal T}(X_1) \otimes {\Cal T}(X_2) \otimes Y \otimes W\rangle\
\endaligned
$$

Part a) follows now from Proposition 1.2.6 b) and the definition of weak
braiding.

To prove part b) we consider the composition
$$
\aligned
& \delta^{(n)}_1 \delta^{(n)}_2(x): \langle V \otimes X_1 \otimes X_2 \otimes Y
\otimes W\rangle\ {\buildrel{s^{-1}_{X_2,X_1}}\over\llongrightarrow}\
\langle V \otimes X_2 \otimes X_1 \otimes Y \otimes W\rangle\
{\buildrel{\widehat{\theta}_n}\over\longrightarrow}\\
& \longrightarrow \langle V \otimes X_1 \otimes Y \otimes {\Cal T}(X_2) \otimes
W\rangle\ {\buildrel{s_{Y,{\Cal T}(X_2)}}\over\llongrightarrow}\
\langle V \otimes X_1 \otimes {\Cal T}(X_2) \otimes Y \otimes W\rangle\
{\buildrel {\widehat{\theta}_n}\over\longrightarrow} \\
& \langle V \otimes {\Cal T}(X_2) \otimes Y \otimes {\Cal T}(X_1) \otimes
W\rangle\
{\buildrel {s_{Y,{\Cal T}(X_1)}}\over\llongrightarrow}\ \langle V \otimes
{\Cal T}(X_2) \otimes {\Cal T}(X_1) \otimes Y \otimes W\rangle \longrightarrow
\\
& {\buildrel {s_{{\Cal T}(X_2),{\Cal T}(X_1)}}\over\longrightarrow} \langle V
\otimes {\Cal T} (X_2 \otimes X_1) \otimes Y \otimes W\rangle\
{\buildrel {{\Cal T}(s_{X_2,X_1})}\over\llongrightarrow}\
\langle V \otimes {\Cal T}(X_1 \otimes X_2) \otimes Y \otimes W\rangle.
\endaligned
$$
It follows now from Proposition 1.2.6 a) that  $\delta_1^{(n)}
\delta_2^{(n)}(x) =
p^*_1(\delta_1^{(n-1)})$.

\proclaim{{\bf 1.2.9}\ Proposition}  The isomorphisms $\delta^{(n)}_i:  {\Cal
F}^{(n)}_{V,W}
\longrightarrow {\Cal E}^*_i({\Cal F}^{(n)}_{V,W})$, $1 \le i \le n$  commute.
\endproclaim

\demo{Proof}  It follows from Proposition 1.2.8 that  $\delta^{(n)}_{i-1}$
and
$\delta^{(n)}_i$  commute.  To prove that $\delta^{(n)}_{i-2}$  and
$\delta^{(n)}_i$  commute we observe that (by Proposition 1.2.8)
$\delta^{(n-2)}_{i-2}$  and
$\delta^{(n-1)}_{i-1}$  commute.  Therefore
$$
p^*_{n-2}(\delta^{(n-1)}_{i-2})p^*_{n-2}(\delta^{(n-1)}_{i-1}) =
p^*_{n-2}(\delta^{(n-1)}_{i-1})p^*_{n-2}(\delta^{(n-1)}_{i-2}).
$$
By Lemma 1.2.8 we can rewrite this equality in the form
$$
\delta^{(n)}_{i-2} \delta^{(n)}_{i-1}\delta^{(n)}_{i}  =
\delta^{(n)}_{i}\delta^{(n)}_{i-2} \delta^{(n)}_{i-1} .
$$
Since the transformation $\delta^{(n)}_{i-1}$  is invertible and,
by the same Proposition 1.2.8 it commutes
with both  $\delta^{(n)}_{i-2}$  and  $\delta^{(n)}_{i}$   we see that
$\delta^{(n)}_{i-2}$  and  $\delta^{(n)}_{i}$  commute.  Analogous arguments
show that  $\delta^{(n)}_{i}, \delta^{(n)}_j$  commute for all  $i,j$,
$1 \le i,j \le n$.  Proposition 1.2.9 is proved.
\enddemo
\vskip .3in

\subheading{1.3.\ A useful formula}

\subheading{1.3.1}\ Let $(H,1,\varepsilon,m,\Delta)$  be a Hopf algebra.
$H_0 {\buildrel {\text{def}}\over =}\ \ker \varepsilon \subset H$  and
$(H \otimes H)_0~\subset~H~\otimes~H$  be the subgroup of linear combinations
of elements of the form
$$
\Delta(x_0)a,\quad x_0 \in H_0,\quad a \in H \otimes H.
$$
Let $S$  be  the antipode of $H$.

\proclaim{{\bf 1.3.2}\ Lemma}  For any $x \in H$  such that
$\Delta(x) = {\dsize\sum^n_{r=0}} x'_r \otimes x''_r$  we have
$$
x \otimes 1 = \sum^n_{r=0} \Delta(x'_r)(1 \otimes S(x''_r)).
$$
\endproclaim

\demo{Proof}  The right side is equal to $(id \otimes m)(id \otimes id \otimes
S)(\Delta \otimes 1)\Delta(x) = (id \otimes m)(id \otimes id \otimes S)
(1 \otimes \Delta)\Delta(x)$.  By the definition of the antipode this
expression is equal to $(id \otimes\epsilon)\Delta(x) = x\otimes 1$.
Lemma 1.3.2 is proved.
\enddemo

\proclaim{{\bf 1.3.3}\ Lemma}\  For any  $x \in H$  we have $x \otimes 1 - 1
\otimes S(x) \in (H \otimes H)_0$.
\endproclaim

\demo{Proof}  Since  $(\epsilon \otimes 1)\Delta(x) = 1 \otimes x$  we can
write
$\Delta(x)$  in the form
$\Delta(x) = 1 \otimes x + {\dsize\sum^n_{r=1}} x'_r \otimes x''_r$   where
$x'_r \in H_0, x''_r \in H, 1 \le r \le n$.  Then it follows from Lemma 1.3.2
that
$$
x \otimes 1 = \sum^n_{r=1} \Delta(x'_r)(1 \otimes S(x''_r)) + 1 \otimes
S(x).
$$

Lemma 1.3.3 is proved.

%this file is qkz2a.tex	 [contains 2.1*]

\newpage

\pageno=17

\centerline{\bf \S2.\ Quantum affine algebras}
\vskip .3in

In this section we will use freely notations from [L].

\subheading{2.1\ Basic definitions}

\subheading{2.1.1}  Let  $(I,\cdot)$  be an affine irreducible Cartan datum
and  $(\Lambda^\vee, \Lambda, \langle\ ,\ \rangle)$  be a simply connected
root datum of type $(I,\cdot)$  (see [L] 2.1 and 2.2) which is an
affinization of a finite root datum.  Any such datum is either a symmetric
root datum (see [L] 2.1) or is obtained as a quotient of a symmetric one
by a finite group of automorphisms which preserve some special vertex (see
[L] 14.1.5).  We denote by $i_0 \in I$  this special vertex.  As follows from
[L] 14.1.4 such a vertex is defined uniquely up to an automorphism of the root
datum and  $(i_0 \cdot i_0) = 2$.
We define  $\overline{I}\ {\dsize\buildrel {\text{def}}\over =}\ I -
\{ i_0\}$.   (In the terminology of [K] we consider the non-twisted affine
case).

Let  ${\Bbb Z}[I] \to \Lambda$  be the group
homomorphism such that  $i \longmapsto i'$  for all  $i \in I$  (see

\noindent
[L] 2.2.1).  This homomorphism is not injective, it has a kernel
isomorphic to ${\Bbb Z}$.  As follows from [K] 6.2 there exists a generator
${\dsize\sum_{i\in I}} n_i\cdot i$ of this kernel such that  $n_{i_0} = 1$.  We
define the dual Coxeter
number $h^\vee$ as the sum\ $h^\vee\ {\dsize\buildrel{\text{def}}\over =}\
{\dsize\sum_{i\in I}} n_i {\frac{(i\cdot i)}{2}}$ and denote by
$\Lambda^\vee_0 \subset \Lambda^\vee$  the subgroup generated by an element
$\Sigma n_i {\frac{(i\cdot i)}{2}} i$.

Let  $\overline{\Lambda}^\vee \subset \Lambda^\vee$  and $\overline{\Lambda}
\subset \Lambda$ be the subgroups generated by elements  $i$ and  $i'$
correspondingly for  $i \in \overline{I}$.  Since $n_{i_0} = 1$  we have a
direct sum
decomposition  $\Lambda^\vee = \overline{\Lambda}^\vee \oplus \Lambda^\vee_0$
which induces a direct sum decomposition $\Lambda =\ {}'\overline{\Lambda}
\oplus
\Lambda_0$, where  $\Lambda_0\ {\dsize\buildrel{\text{def}}\over =}\
\{ \lambda \in \Lambda| \langle i,\lambda\rangle = 0\ \forall i \in
\overline{I}\}$  and $'\overline{\Lambda} = \overline{\Lambda} \otimes_{\Bbb Z}
{\Bbb Q} \cap \Lambda$.  Then $\overline{\Lambda}$  is a subgroup of finite
index $d$ in $'\overline{\Lambda}$.  As follows from the definition of a root
datum we have  $(i \cdot i) \langle i,\lambda'\rangle = 2(i \cdot \lambda)$
for
all $i \in I, \lambda \in \Lambda^\vee$.  The map  $\lambda \longmapsto
\lambda'$
defines an imbedding  $\overline{\Lambda}^\vee \hookrightarrow\
{}'\overline{\Lambda}$  and there exists unique symmetric bilinear form
$[\ ,\ ]:\  {}'\overline{\Lambda} \times '\overline{\Lambda} \to 1/d {\Bbb Z}$
such that  $[\lambda',\mu'] = (\lambda \cdot \mu)$  for all
$\lambda, \mu \in \overline{\Lambda}^\vee$.

We denote by  $\rho' \in\ {}'\overline{\Lambda}$  the unique element such that
$\langle i,\rho'\rangle = 1$  for all $i \in \overline{I}$.  Then
$[i',\rho'] = {\frac{(i\cdot i)}{2}}$ for all $i \in \overline{I}$.

\demo{Remark}  In the terminology of [K], $\Lambda$ is the weight lattice,
$\Lambda^\vee$  is the coroot lattice,\ ${}'\overline{\Lambda}$  is the weight
lattice of finite-dimensional root datum  $(\overline{I}, \cdot)$.
\enddemo

\subheading{2.1.2}  We fix a number $\widetilde{q} \in {\Bbb C}^*$  such that
$|\widetilde{q}| < 1$  and define  $q\ {\buildrel{\text{def}}\over =}\
\widetilde{q}^{dh^\vee}$,\  $q_i\ {\buildrel{\text{def}}\over =}\ q^{(i\cdot
i)/2}$
for all $i \in I$.

For any $n \in {\Bbb N}, i \in I$  we define $[n]_i\
{\buildrel{\text{def}}\over =}\
{\frac{q^n_i - q^n_i}{q_i-q_i^{-1}}}$  and  $[n]_i! {\buildrel{\text{def}}\over
=}\
\prod^n_{s=1} [s]_i$.

\subheading{2.1.3}  We  denote by
$\check{\bold U}$  the ${\Bbb C}$-algebra generated by elements $E_i,F_i,
K_\mu,
i \in I, \mu \in \Lambda^\vee$  and relations  (a)-(e) below.

\noindent
(a)  $K_0 = 1, K_\mu K_{\mu'} = K_{\mu + \mu'}$, $\mu,\mu' \in \Lambda^\vee$.

\noindent
(b)  $K_\mu E_i = q^{\langle \mu,i'\rangle} E_iK_\mu$, $i \in I, \mu \in
\Lambda^\vee$.

\noindent
(c)  $K_\mu F_i = q^{- \langle \mu,i'\rangle} F_iK_\mu$, $i \in I, \mu \in
\Lambda^\vee$.

\noindent
(d)  $E_iF_j - F_jE_i = \delta_{i,j}\
{\dsize\frac{\widetilde{K}_i-\widetilde{K}_{-i}}{q_i - q_i^{-1}}}$  for
$i,j \in I$, where
$$
\widetilde{K}_{\pm i} = K_{\pm {\frac{(i\cdot i)}{2}} i}
$$

\noindent
(e)  ${\dsize\sum_{p+p' = 1 - 2i \cdot j/(i \cdot i)}} (-1)^{p'} E_i^{(p)}
E_j E_i^{(p')} = 0$\ for all $i \ne j \in I$,
where $E_i^{(p)}\ {\dsize\buildrel{\text{def}}\over =}\ E^p_i/[p]_i!$

\subheading{2.1.4}  Let  $\widetilde{Z}\ {\dsize\buildrel\text{def}\over  =}\
{\dsize\prod_{i\in I}} \widetilde{K}_i^{n_i}$.  Then $\widetilde{Z}$  lies in
the center
of  $\check{\bold U}$  we define ${\bold U} = \check{\bold U}[Z]$  where  $Z$
is defined as a central
element such that  $Z^{dh^\vee} = \widetilde{Z}$.

\subheading{2.1.5}  Let $\Lambda_{{\Bbb C}^*}$ be the tensor product
$\Lambda_{{\Bbb C}^*}\ {\dsize\buildrel{\text{def}}\over =}\
\Lambda \otimes {\Bbb C}^*$  and  $\Lambda\ \hookrightarrow\
\Lambda_{{\Bbb C}^*}$  be the imbedding induced by the imbedding ${\Bbb Z}
\hookrightarrow {\Bbb C}^*: n \longmapsto q^n$.   We extend an imbedding
${\Bbb Z} \hookrightarrow {\Bbb C}^*$  to an imbedding\
${\frac{1}{d}} {\Bbb Z} \hookrightarrow {\Bbb C}^*$  in such a way that\
${\frac{1}{d}} \longmapsto \widetilde{q}^{h^\vee}$.  This imbedding defines an
imbedding  ${\frac{1}{d}} \Lambda \hookrightarrow \Lambda_{{\Bbb C}^*}$.

We denote the group structure on $\Lambda_{{\Bbb C}^*}$
as  $+$.  The pairing $\langle\ ,\ \rangle: \Lambda^\vee \times \Lambda \to
{\Bbb Z}$,  $(\mu, \lambda) \longmapsto \langle \mu,\lambda\rangle$  defines
the pairing
${\frac{1}{d}} \Lambda^\vee \times \Lambda_{{\Bbb C}^*} \to {\Bbb C}^*$  which
we will also
denote by  $(\mu, \lambda) \longmapsto \langle \mu,\lambda\rangle$.

For any representation $(\rho, \underline{V})$  of  ${\bold U}$  and $\lambda
\in
\Lambda_{{\Bbb C}^*}$  we define

$\underline{V}_\lambda = \{ v \in V\big\vert K_\mu v =
\langle \mu,\lambda\rangle v\}$  for all $\mu \in \Lambda^\vee$.

\proclaim{Definition}  We denote by  ${\Cal C}$  the category of
representations
$(\rho, \underline{V})$  of  ${\bold U}$  such that $\underline{V} =
{\dsize\bigoplus_{\lambda\in\Lambda_{{\Bbb C}^*}}} \underline{V}_\lambda$.
For any commutative ring $B$ containing ${\Bbb C}$  we keep the notation
${\Cal C}$  for the  category of  ${\bold U}_B$-modules (where
${\bold U}_B = {\bold U} \otimes_{\Bbb C} B)$  obtained from ${\Cal C}$ by
changing the scalars.
\endproclaim

For any complex-valued function  $F$  on $\Lambda_{{\Bbb C}^*}$ and any
$V = (\rho,\underline{V})$  in ${\Cal C}$  we denote by  $\underline{F}$  the
linear endomorphism of $\underline{V}$  which preserves the direct sum
decomposition  $V = {\dsize\oplus_{\lambda\in\Lambda_{{\Bbb C}^*}}} V_\lambda$
and such that

\noindent
$\underline{F}|_{V_\lambda} = F(\lambda)Id_{V_\lambda}$  for all
$\lambda \in \Lambda_{{\Bbb C}^*}$.

For any  $V = (\rho, \underline{V})$  in ${\Bbb C}$  we define
$\underline{V}^*$
as the direct sum  $\underline{V}^*\ {\dsize\buildrel{\text{def}}\over =}\
{\dsize\oplus_{\lambda\in\Lambda_{{\Bbb C}^*}}}
\text{Hom}(\underline{V}_\lambda,
{\Bbb C})$.  Then  $\underline{V}^*$  has a natural structure of a ${\bold
U}$-module.
We denote this ${\bold U}$-module as $V^*$.  By the definition  $V^*$  belongs
to
${\Cal C}$.

In this paper by an expression ``a ${\bold U}$-module'' we will always
understand
``a ${\bold U}$-module from the category ${\Cal C}$''.

\subheading{2.1.6}  For any  $z \in {\Bbb C}^*$  we denote by  ${\Cal C}_z$
the
full subcategory of  ${\Cal C}$  of representations  $(\rho, \underline{V})$
such that  $\rho(Z) = z Id_V$.  The full subcategory of ${\Cal C}_1$
consisting
of finite-dimensional representations is denoted by  ${\Cal D}$.

\subheading{2.1.7}  Let $\overline{\Lambda}_{{\Bbb C}^*}\
{\dsize\buildrel{\text{def}}\over =}\
\overline{\Lambda} \otimes {\Bbb C}^*$  and
$r: \Lambda_{{\Bbb C}^*} \to {}'\overline{\Lambda}_{{\Bbb C}^*}$  be
the projection induced by the direct sum decomposition $\Lambda =
{}'\overline{\Lambda} \oplus \Lambda_0$.  By the definition $\langle
i,\lambda\rangle
= 0$ for all $i \in \overline{I}$, $\lambda \in \ker r$.   Therefore we
obtain a pairing  ${}'\overline{\Lambda}^\vee  \times \overline{\Lambda}_{{\Bbb
C}^*}
\to {\Bbb C}^*$  which we also denote as  $\langle\ ,\ \rangle$.  It is clear
that this pairing coincides with the restriction of the pairing
$\langle\ ,\ \rangle: {\frac{1}{d}} \Lambda^\vee \times \Lambda_{{\Bbb C}^*}
\to
{\Bbb C}^*$  to  ${}'\overline{\Lambda}^\vee \times \overline{\Lambda}_{{\Bbb
C}^*}.$

The imbedding  $\Lambda \hookrightarrow \Lambda_{{\Bbb C}^*}$  induces the
imbedding  $\overline{\Lambda} \hookrightarrow \overline{\Lambda}_{{\Bbb
C}^*}$.
We denote by  ${\frak S}$  the torus $\overline{\Lambda}_{{\Bbb
C}^*}/\overline{\Lambda}$
and by  $\Pi:  \overline{\Lambda}_{{\Bbb C}^*} \to {\frak S}$  the natural
projection of
$\overline{\Lambda}_{{\Bbb C}^*}$  to ${\frak S}$.

\proclaim{Definition}  a)  For any $s \in {\frak S}$  we denote by
${}^s{\Cal C}$  the full subcategory ${\Cal C}$  of ${\bold U}$-modules
$(\rho,\underline{V})$  such that  $\underline{V}_\lambda = \{ 0\}$  for all
$\lambda \in \Lambda_{{\Bbb C}^*}$  such that $\Pi(r(\lambda)) \ne s$.

b)  For any  $s \in {\frak S}, z \in {\Bbb C}^*$  we define  ${}^s{\Cal C}_z$
as
the intersection of  ${}^s{\Cal C}$  and  ${\Cal C}_z$.
\endproclaim

\demo{Remark}  In this paper, we will almost always assume that our
${\bold U}$-modules are in the subcategory
$$
{}^{[0]}{\Cal C}\ {\buildrel {\text{def}}\over =}\
{\bigcup_{s\in\Pi('\overline{\Lambda})}} {}^s{\Cal C} ,
$$
The reason to consider the more general category ${\Cal C}$  is to have a
possibility to prove results for ${}^s{\Cal C}$, where $s$ is ``generic'' first
and then to ``deform''  ${}^{[0]}{\Cal C}$  to ${}^{-s}{\Cal C}$  for
``generic'' $s$ (see, for example, the proof of Theorem 3.2.2).

\subheading{2.1.8}  Let  ${\bold U}_0 \subset {\bold U}$  be the subalgebra
generated by $Z$ and  $K_\mu$, $\mu \in \Lambda^\vee$ and
${\bold U}^{\bold f}_0 \subset {\bold U}_0$  its subalgebra generated by
$K_\mu, \mu \in \overline{\Lambda}^\vee$.  For any  $\nu \in {\Bbb N}[I]$
we denote by  ${\bold U}^\nu_+$  the subspace of ${\bold U}$ spanned by
elements of the form
$x^+$  where $x \in {\bold f}_\nu$ (see [L] 3.1.1).  We define ${\bold U}^n_+$
as
a direct sum
$$
{\bold U}^n_+\ {\buildrel{\text{def}}\over =}\ {\bigoplus\Sb\\ \nu \in {\Bbb
N}[I]\\
tr \nu = n\endSb} {\bold U}^\nu_+
$$
and denote by  ${\bold U}^{\ge n}_+$  the direct sum
${\dsize\oplus_{n'\ge n}} {\bold U}_+^{n'}$.
We denote    ${\bold U}_+^{\ge 1}$  as  ${\bold U}_+^>$  and define  ${\bold
U}_+ \subset {\bold U}$  to be
the subalgebra generated by  ${\bold U}_0$  and  ${\bold U}_+^>$  and
${\bold U}_+^{\bold f} \subset {\bold U}_+$  as the subalgebra generated by
${\bold U}_0^{\bold f}$  and  ${\bold U}_+^>$.  Then
${\bold U}_+ = {\bold U}_+^{\bold f}[Z,Z^{-1}]$.

\subheading{2.1.9}  We denote by  $\overline{\bold U} \subset {\bold U}$  the
subalgebra generated
by  $K_\mu, \mu \in \overline{\Lambda}^\vee$  and  $E_i,F_i$  for $i \in
\overline{I}$.  Then $\overline{\bold U}$  is the quantum algebra corresponding
to
the root datum  $(\overline{\Lambda}^{\vee},\ {}'\overline{\Lambda})$.

\subheading{2.1.10}  Let $\Delta: {\bold U} \to {\bold U} \to {\bold U} \otimes
{\bold U}$  be the
comultiplication such that
$$
\matrix \format \l & \l &\qquad\qquad \l \\
\Delta(E_i) &= E_i \otimes 1 + \widetilde{K}_i \otimes E_i, & i \in I\\
\vspace{1.5\jot}
\Delta(F_i) &= F_i \otimes \widetilde{K}_{-i} + 1 \otimes F_i, & i \in I\\
\vspace{1.5\jot}
\Delta(K_\mu) &= K_\mu \otimes K_\mu,  & \mu \in \Lambda^\vee\\
\vspace{1.5\jot}
\Delta(Z) &= Z \otimes Z
\endmatrix
$$
(see [L] 23.1.5) .

We define a co-unit $\epsilon$ on ${\bold U}$  as the unique
homomorphism  $\epsilon: {\bold U} \to {\Bbb C}$  such that
$$
\epsilon(E_i) = \epsilon(F_i) = 0\quad i \in I,
$$
$\epsilon(K_\mu) = 1$,\ $\mu \in \Lambda^\vee$  and $\epsilon(Z) = 1$.
Then  ${\bold U}$  is a Hopf algebra where the antipode $S$ is given by
the formulas  $S(E_i) = - \widetilde{K}_{-i} E_i$, $S(F_i) =
- F_i \widetilde{K}_i$, $S(K_\mu) = K_{-\mu}$ and  $S(Z) = Z^{-1}$.

This comultiplication
defines a strict monoidal structure on the category ${\Cal C}$.

\subheading{2.1.11}  Let   $\omega$  be the involution of the algebra
${\bold U}$  such that
$$
\omega(E_i) = F_i,\ \omega(F_i) = E_i,\ i \in I,\ \
\omega(K_\mu) = K_{-\mu},\ \mu \in \Lambda^\vee,\quad \omega(Z) = Z^{-1}
$$
(see [L] 3.1.3).

For any ${\bold U}$-module $(\rho,M)$ we denote by ${}^\omega M$  the
representation $(\widetilde{\rho}, \widetilde{M})$  of
${\bold U}$ on the same space $M$ such that  $\widetilde{\rho}(x) =
\rho(\omega(x))$.

\subheading{2.1.12} Let $F$ be a field containing ${\Bbb C}$, ${\bold U}_F =
{\bold U} \otimes_{\bold C} F$.
For any central invertible $u \in {\bold U}_F$  we define two
automorphisms  $\varphi$  and $\psi$  of the Hopf algebra ${\bold U}_F$
where

\medpagebreak

\noindent
$\varphi_u(E_{i_0}) = u^{dh^\vee}E_{i_0},\ \varphi_u(F_{i_0}) =
u^{-dh^\vee}F_{i_0},\
\varphi_u(E_i) = E_i,\ \varphi_u(F_i) = F_i,\ i \in \overline{I},$\

\noindent
$\varphi_u(K_\mu) = K_\mu,$\qquad  $\mu \in \Lambda^\vee,\ \varphi_u(Z) = Z,$

\noindent
$\psi_u(E_i) =
u^{\frac{d(i\cdot i)}{2}} E_i,\ i \in I,\ \psi_u(F_i) =
u^{-\frac{d(i\cdot i)}{2}} F_i,\ i \in I,\
\psi_u(K_\mu) = K_\mu,\ \mu \in \Lambda^\vee,\ \psi_u(Z) = Z.$

\medpagebreak

\demo{Remark}  The automorphisms $\varphi_u$  and $\psi_u$  of ${\bold U}_F$
define strict automorphisms of the category ${\Cal C}$  which we denote
as  ${\Cal T}_u^\varphi$ and  ${\Cal T}_u$ respectively.
\enddemo

\subheading{2.1.13}  Fix a point $s \in {\frak S}$  and an element
$a \in \Pi^{-1}(s) \subset \overline{\Lambda}_{{\Bbb C}^*}$.    For any central
invertible element $u \in  {\bold U}_F$  we denote by ${\Cal L}^a_u$
the function on $\Pi^{-1}(s)$  with values in ${\bold U}_F$  such that
${\Cal L}^a_u(a) = 1$  and
$$
{\Cal L}^a_u(\overline{\lambda} + i') =
u^{d(i\cdot i)/2}{\Cal L}^a_u(\overline{\lambda})\quad
\text{for\ any}\quad \overline{\lambda} \in \Pi^{-1}(s), i \in \overline{I}.
$$

We extend  ${\Cal L}^a_u$  to  $\Lambda_{{\Bbb C}^*}$  in such a way that
${\Cal L}^a_u(\lambda) = 0$  if  $\lambda \notin \Pi^{-1}(s)$,
$\lambda \in \overline{\Lambda}_{{\Bbb C}^*}$  and ${\Cal L}^a_u$  is constant
on the fibers of $r$.

\proclaim{Definition}  	For any $V = (\rho, \underline{V}) \in
{}^s{\Cal C}$  we denote by  $\underline{\Cal L}^a_u$  a linear endomorphism
of  $\underline{V}$  corresponding to the function  ${\Cal L}^a_u$  (see
2.1.5).
\endproclaim

\demo{Remark}  For any two  $a',a'' \in \Pi^{-1}(s)$  the operators
$\underline{\Cal L}^{a'}_u$  and  $\underline{\Cal L}^{a''}_u$  are
proportional and the coefficient of proportionality does not depend on a
choice of  $V$  in  ${}^s{\Cal C}$.  Therefore we can consider
$\underline{\Cal L}^a_u$  as a section of a line bundle
$\check{\underline{\Cal L}}_u$  on  ${\frak S}$.
\enddemo

\proclaim{{\bf 2.1.14.}\ Proposition}  For any $V = (\rho, \underline{V})$ in
${}^s{\Cal C}$  the map  ${\underline{\Cal L}}_u^a$  defines an isomorphism of
the ${\bold U}_F$-module  ${\Cal T}^\varphi_u(V)$  with the ${\bold
U}_F$-module ${\Cal T}_u(V)$.
\endproclaim

\demo{Proof}  We have to show that for any  $x \in {\bold U}_F$,\ we have
$\underline{\Cal L}^a_u \varphi_u(x) = \psi_u(x) \underline{\Cal L}^a_u$.
This equality is obvious if  $x = K_\mu$,  $\mu \in \Lambda^\vee$  or
$x = Z$.    In the case when $x = E_i$ or $F_i$, $i \in
\overline{I}$  the claim follows immediately from the definition of
${\Cal L}^a_u$  and the equalities  $x^+\underline{V}_\lambda \subset
\underline{V}_{\lambda +\nu}$,\ $x^-\underline{V}_\lambda \subset
\underline{V}_{\lambda-\nu}$  for $x \in {\bold f}_\nu$.  In the
case when $x = E_{i_0}$ or $F_{i_0}$  the claim follows from the definition
of the dual Coxeter number $h^\vee$.   Proposition 2.1.14 is proved.
\enddemo

\proclaim{{\bf 2.1.15}\ Proposition}  For any finite-dimensional
$X = (\rho_X, \underline{X}) \in {\Cal C}$ the identity map
$\underline{X} \longrightarrow 	\underline{X}$  defines an isomorphism
$X^{**} \longrightarrow {\Cal T}_{\widetilde{q}^{-2h^\vee}}(X)$ of ${\bold
U}_F$-modules.
\endproclaim

\demo{Proof}  By the definition  $X^{**} = (\rho_{X^{**}},\underline{X})$
where
$\rho_{X^{**}}(a) = \rho_X(S^2(a))$, $a \in {\bold U}_F$  and
$S: {\bold U}_F \to {\bold U}_F$  is the
antipode.  The claim follows now from the formulas for the antipode $S$ in
2.1.10.
\enddemo
\vskip .3in

%this file is qkz2b.tex	  [contains 2.2*]

\subheading{2.2\ Sugawara operators}

\subheading{2.2.1}  For any $V = (\rho, \underline{V})$  in ${\Cal C}$  and
$n \in {\Bbb N}$  we define a subspace $\underline{V}(n)$ of $\underline{V}$
$$
\underline{V}(n)\ {\buildrel{\text{def}}\over =}\
\left \{  v \in \underline{V} \big\vert av = 0\ \forall a \in
{\bold U}^{\ge n}_+ \right \}.
$$
It is clear that  $\underline{V}(1) \subset \underline{V}(2) \subset \cdots
\subset \underline{V}(n) \subset \cdots$.  We define a subspace
$\underline{V}(\infty)$ of  $\underline{V}$  as the union of all
$\underline{V}(n)$, $n > 0$.

\proclaim{Proposition}  For any  $i \in I$,\  $n \in {\Bbb N}$  we have
$$
(E_i) \underline{V}(n) \subset \underline{V}(n)\quad \text{and}\quad
(F_i) \underline{V}(n) \subset \underline{V}(n+1).
$$
\endproclaim

\demo{Proof}  The first part follows immediately from the definitions and the
second follows from [L] 3.1.6.
\enddemo

\proclaim{Corollary}  $\underline{V}(\infty)$  is a ${\bold U}$-invariant
subspace of $\underline{V}$.
\endproclaim
\vskip .2in

We denote the corresponding ${\bold U}$-module by  $V(\infty)$.  Let  ${\Cal
C}^+$
be the full subcategory of ${\Cal C}$  of modules $V$  such that  $V(\infty)
= V$  and  let  ${\Cal C}^-$  be the subcategory of  ${\Cal C}$  of ${\bold
U}$-modules
$W$  such that  ${}^wW \in {\Cal C}^+$.  For any $W$ in ${\Cal C}^-$  we define
$D(W) = W^*(\infty)$  where  $W^*$  is as in 2.1.5. We define ${}^s{\Cal
C}^\pm$
as\ ${}^s{\Cal C} \cap {\Cal C}^\pm$  and  ${\Cal C}^\pm_z$  as
${\Cal C}^\pm \cap {\Cal C}_z$.

\demo{Remark}  For any ${V}$  in ${\Cal C}$  the submodule
$V(\infty) \subset {V}$  is the maximal subobject of ${V}$
contained in ${\Cal C}^+$.
\enddemo

\subheading{2.2.2}  For any  $V$ in ${\Cal C}^+$  and  $W$  in ${\Cal C}^-$
we define a linear map\quad   Hom$_{\bold U}(V,D(W)) \to$

\noindent
$\text{Hom}_{\bold U}(V \otimes W, {\Bbb C})$
as in [KL] 2.30.  As follows from 1.1.6 we can identify the linear space
Hom$_{\bold U}(V \otimes W, {\Bbb C})$  with  $\langle V \otimes
W\rangle^\vee$.

\proclaim{Lemma}  The map Hom$_{\bold U}(V,D(W)) \to
\langle V \otimes W\rangle^\vee$  is an isomorphism.
\endproclaim

\demo{Proof}  Analogous to the proof of Lemma 2.3.1 in [KL].
\enddemo

\subheading{2.2.3}  Let  $\Omega_{\le p}$  be elements defined in [L] 6.1.1.
Then for any  $M = (\rho, \underline{M})$ in ${\Cal C}^+$  there exists an
operator $\Omega$
on $\underline{M}$  such that  $\Omega(m) = \Omega_{\le p}(m)$  for any
$m \in M$  and all sufficiently big $p$.  We have
$\widetilde{K}_{-i} E_i\Omega = \widetilde{K}_i\Omega E_i,$\
$\Omega F_i = F_i \widetilde{K}_i \Omega \widetilde{K}_i$  and
$\Omega K_\mu = K_\mu\Omega$, $i \in I, \mu \in \Lambda^\vee$.

\subheading{2.2.4}   Fix a point  $s \in {\frak S}$  and an element
$a \in \Pi^{-1}(s)
\subset  \overline{\Lambda}_{{\Bbb C}^*}$.   Let $G^a: \Pi^{-1}(s) \to
{\Bbb C}^*$    be the function such that  $G^a(a) = 1$  and for any
$\overline{\lambda} \in \Pi^{-1}(s)$,\ $i \in \overline{I}$  we have
$$
G^a(\overline{\lambda}) G^a(\overline{\lambda} - i')^{-1} =
\langle i,\overline{\lambda}\rangle^{(i\cdot i)} .
$$
We continue  $G^a$  to  $\Lambda_{{\Bbb C}^*}$  in the same way as we did
with  ${\Cal L}^a_u$  in 2.1.13.
For any  $V = (\rho,\underline{V})$  in ${}^s{\Cal C}$  we denote by
$\underline{G}^a$  the corresponding endomorphism of  $\underline{V}$ (see
(2.1.5)).

\demo{Remark}  As in the case of  $\underline{\Cal L}^a_u$  we can consider
$\underline{G}^a$  as a section of a line bundle on ${\frak S}$.
\enddemo

\proclaim{Definition}  For any $s \in {\frak S}$, $V = (\rho_V,\underline{V})
\in
{}^s{\Cal C}^+$  and  $a \in \Pi^{-1}(s)$  we define linear endomorphisms
$T^a_V$ (or simply $T^a$) of  $\underline{V}$  as composition\
$T^a\ {\buildrel{\text{def}}\over =}\
\underline{\Cal L}^a_{(Z\widetilde{q}^{h^\vee})^{-2}} \Omega \underline{G}^a$.
We will call all of them Sugawara operators.
\endproclaim

\proclaim{Proposition}
$$
T^a \in \text{Hom}_{\bold U}(V, {\Cal
T}_{(Z{\widetilde{q}}^{h^\vee})^{-2}}(V)).
$$
\endproclaim

\demo{Proof}  We have to check that for any  $x \in {\bold U}$  we have
$T^ax = \psi_{(Z{\widetilde q}^{h^\vee})^{-2}}(x) T^a$.  Let

\noindent
$'T^a\ {\dsize\buildrel{\text{def}}\over =}\ \Omega \cdot \underline{G}^a$.
It is sufficient to prove that

\noindent
$\alpha$)\quad $'T^a K_\mu = K_\mu\ {}'T^a,\ 'T^a Z = Z 'T^a$

\noindent
$\beta$)\quad $'T^a E_i = E_i\ {}'T^a$,\ $'T^a F_i = F_i\ {}'T^a$\ for\ $i \in
\overline{I}$

\noindent
$\gamma$)  $'T^a E_{i_0} = (Z{\widetilde q}^{h^\vee})^{-2dh} E_{i_0}\
{}'T^a$,\
$'T^a F_{i_0} = (Z{\widetilde q}^{h^\vee})^{2dh^\vee} F_{i_0}\ {}'T^a$.

\noindent
The equalities $\alpha$) are obviously true.  The proof of equalities
$\beta$) is completely analogous to the proof of Proposition 6.1.7 in [L].
So we give only the proof of equalities  $\gamma$).  Actually we only give a
proof of the equality  $'T^a E_{i_0} = (Z{\widetilde q^{h^\vee}})^{2dh^\vee}
E_{i_0}\ {}'T^a$.  The proof of the second part of  $\gamma$) is completely
analogous.
\enddemo

If  $v \in \underline{V}_\lambda$ then  $E_{i_0}v \in \underline{V}_{\lambda +
i'_0}$.  Let  $\nu\ {\dsize\buildrel {\text{def}}\over =}\
{\dsize\sum_{i\in \overline{I}}} n_i i \in \Lambda$.  Since
$\lambda + i'_0 = \lambda - \nu'$  we have
$$
\aligned
'T^a E_{i_0} v &= G^a(\lambda - \nu') \Omega E_{i_0} v = G^a(\lambda - \nu')
K_{i_0}^{-(i_0\cdot i_0)} E_{i_0} \Omega v = \\
&= G^a(\lambda - \nu') K_{i_0}^{-2} E_{i_0} \Omega v = G^a(\lambda - \nu')
\widetilde{Z}^{-2} {\dsize\prod_{i\in \overline{I}}}
\widetilde{K}_i^{2n_i} E_{i_0} \Omega v = \\
&= Z^{-2dh^\vee} G^a(\lambda -\nu')G^a(\lambda)^{-1}
{\dsize\prod_{i\in \overline{I}}} \widetilde{K}_i^{2n_i}  E_{i_0}\ {}'T^av =
(Z \widetilde{q}^{h^\vee})^{-2dh^\vee} E_{i_0}\ {}'T^av.
\endaligned
$$

Proposition 2.2.4 is proved.

\medpagebreak

\subheading{2.2.5}  In the case when $V = (\rho_V, \underline{V})$ lies in
${}^{[0]}{\Cal C}$  we can give a more explicit formula for the operators
$\underline{\Cal L}_u$  and  $\underline{G}$  on $\underline{V}$.  Let
${\Cal L}_u$  be the function on $'\overline{\Lambda}$  with values in the
center of ${\bold U}$ and let $G$  be the complex-valued function on
$'\overline{\Lambda}$
such that
$$
{\Cal L}_u(\lambda) = u^{d[\lambda,\rho']}\quad
\text{and}\quad G(\lambda) = q^{([\lambda +\rho',
\lambda + \rho'] - [\rho',\rho'])}
$$
for $\lambda \in {}'\overline{\Lambda}$  where the bilinear form
$[\ ,\ ]$  is as 2.1.1.  They define $\underline{\Cal L}_u$  and
$\underline{G}$  in the obvious way.

Let $(\underline{M} \otimes \underline{N})_0 \subset \underline{M} \otimes
\underline{N}$  be the subspace defined as in 1.1.6.

\medpagebreak

\subheading{2.2.6}  For any  $V = (\rho_V, \underline{V})$  in
${}^s{\Cal C}^-$  we define an endomorphism  $\check{T}^a$  of
$\underline{V}$   by the rule
$$
\check{T}^a\ {\buildrel {\text{def}}\over =}\ (T_{{}_{{}^\omega V}}^a)^{-1},
$$
where the module  ${}^\omega V \subset {}^s{\Cal C}^+$  is defined as in
2.1.11.

\demo{Remark}  We can consider  $T^a$  and  $\check{T}^a$  as elements in
appropriate completions of ${\bold U}$.  Then $\check{T}^a = \omega(T^a)$.
\enddemo

\proclaim{{\bf 2.2.7}\ Proposition}  For any  $i \in I$  we have

\noindent
a)\quad  $T^a E_i = (Z\widetilde{q}^{h^\vee})^{-(i\cdot i)d} E_i T^a$,

\noindent
b)\quad  $\check{T}^a E_i = (Z\widetilde{q}^{-h^\vee})^{-(i\cdot i)d}
E_i \check{T}^a,$

\noindent
c)\quad  $T^a F_i = (Z\widetilde{q}^{h^\vee})^{(i\cdot i)d} F_i T^a$,

\noindent
d)\quad  $\check{T}^a F_i = (Z\widetilde{q}^{-h^\vee})^{-(i\cdot i)d} F_i
\check{T}^a,$.
\endproclaim

\demo{Proof}  Parts a) and c) are corollaries of Proposition 2.2.4.  We prove
b).  The proof of d) is completely analogous.  We have
$$
\check{T}^aE_i = \omega[((T^a)^{-1}F_i)] =
\omega[(Z\widetilde{q}^{h^\vee})^{-(i\cdot i)d} F_i(T^a)^{-1}] =
(Z\widetilde{q}^{-h^\vee})^{(i\cdot i)d} E_i \check{T}^a,
$$
where the second equality is a restatements of c).   Proposition 2.2 is proved.
\enddemo

\subheading{2.2.8}  For any  $V = (\rho_V, \underline{V})$ in ${}^s{\Cal C}^+$
and $W = (\rho_W, \underline{W})$  in  ${}^{-s}{\Cal C}^-$  the linear map
$T^a \otimes \check{T}^a$  of  $\underline{V} \otimes \underline{W}$  does not
depend on a choice of  $a \in \Pi^{-1}(s)$  and we denote this operator as
$T \otimes \check{T}$.

\proclaim{Proposition} $(T \otimes \check{T})$  preserves the subspace
$(M \otimes N)_0$ of $\underline{M} \otimes \underline{N}$.
\endproclaim

\demo{Proof}  It is sufficient to show that for any  $m \in M,\ n \in N$,
$i \in I$  and any  $\mu \in \Lambda^\vee$,  we have

\item{a)}  $(T \otimes \check{T})(E_i) (m \otimes n) \in (M \otimes N)_0$,

\item{b)}  $(T \otimes \check{T})(F_i)(m \otimes n) \in (M \otimes N)_0$,

\item{c)}  $(T \otimes \check{T})[(K_\mu)(m \otimes n) - (m \otimes n)] \in
(M \otimes N)$.
\enddemo

\demo{Proof} of a):  By the definition of the action of ${\bold U}$  on
$M \otimes N$  we have
$$
(T \otimes \check{T}) E_i(m \otimes n) = TE_im \otimes \check{T}n +
T\widetilde{K}_im \otimes \check{T} E_in.
$$
\enddemo

Therefore it follows from Proposition 2.2.7 that
$$
\aligned
(T \otimes \check{T}) E_i (m \otimes n)
&= (Z\widetilde{q}^{h^\vee})^{-d(i\cdot i)} E_i Tm \otimes \check{T}n
+ \widetilde{K}_i Tm \otimes (Z\widetilde{q}^{-h^\vee})^{d(i\cdot i)}
E_i \check{T}n  \\
&= \widetilde{q}^{-d(i\cdot i)h^\vee} [\Delta (E_i)(Z^{-d(i\cdot i)} Tm \otimes
\check{T}n) - (Z^{-d(i\cdot i)} A \otimes B - A \otimes Z^{d(i\cdot i)}B) \\
&= \widetilde{q}^{-d(i\cdot i)h^\vee} \Delta (E_i)(Z^{-d(i\cdot i)} Tm \otimes
\check{T}n) - (\Delta(Z^{-d(i\cdot i)}) \\
&\qquad  - \epsilon(Z^{-d(i\cdot i)})) (A \otimes Z^{d(i\cdot i)}B),
\endaligned
$$
where  $A = \widetilde{K}_iTm$  and  $B = \check{T}n$.

The inclusion a) is proved.  The proof of the inclusion is completely analogous
and the proof of c) is obvious.

Proposition 2.2.7 is proved.

\demo{Remark}   We will show later (see 3.2) that $T \otimes \check{T}$
induces
the identity map on the quotient  $\langle M,N\rangle\
{\dsize\buildrel{\text{def}}\over =}\ \underline{M}\otimes
\underline{N}/(M \otimes N)_0$.
\enddemo

We extend ${\Cal L}_u$  and $G$ to functions on  $\Lambda_{{\Bbb C}^*}$  which
are constant on fibers of $r$ and are zero outside of
$r^{-1}('\overline{\Lambda})$
and denote by  $\underline{\Cal L}_u$  and  $\underline{G}$  the corresponding
automorphism of  $\underline{V}$  for any  $V = (\rho, \underline{V})$  in
${}^{[0]}{\Cal C}$.  Then $\underline{\Cal L}_u$  and $\underline{G}$  are
sections of the line bundles  $\check{\underline{\Cal L}}_u$	 and
$\check{\underline{G}}$  over the finite set  $\Pi('\overline{\Lambda})
\subset {\frak S}$.  For any  $V = (\rho_V, \underline{V})$  in
${}^{[0]}{\Cal C}^+\ ({\dsize\buildrel{\text{def}}\over =}\  {\Cal C}^+ \cap
{}^{[0]}{\Cal C})$  and  $W = (\rho_W, \underline{W})$ in
${}^{[0]}{\Cal C}^-\ ({\dsize\buildrel{\text{def}}\over =}\ {\Cal C}^- \cap
{}^{[0]}{\Cal C})$  we denote by  $T \in \text{End}\ \underline{V}$  and
$\check{T} \in \text{End}\ \underline{W}$  the Sugawara operators
corresponding to the functions  $\underline{\Cal L}_u$  and $\underline{G}$.

\proclaim{{\bf 2.2.9} Proposition}  Let $M = (\rho_M,\underline{M}) \in
{}^{[0]}{\Cal C}^+_z$  be a ${\bold U}$-module, $\lambda \in
{\frac{1}{d}} \overline{\Lambda}$ and
let $m \in \underline{M}_{\lambda},$  be such that  $E_im = 0$  for all
$i \in I$.  Then
$$
Tm = q^{[[(\lambda + \rho'), (\lambda + \rho')] - [\rho',\rho']]}
z^{-2d[\lambda,\rho']} q^{-2[\lambda,\rho']}m =
q^{[\lambda,\lambda]} z^{-2d[\lambda,\rho']} m .
$$
\endproclaim

\demo{Proof}  Follows immediately from the definition of $T$.

\subheading{2.2.10}  Let  $M = (\rho_{\underline{M}}, \underline{M})$
be a ${\bold U}$-module in\ ${}^{[0]}{\Cal C}$.  For
any $n \in {\Bbb N}$  we denote
by $M_{(n)} \subset M$  the subspace generated by vectors of the norm
$xm$, $x \in {\bold U}^{\geqslant n}_-$, $m \in \underline{M}$.
For any $N = (\rho_N, \underline{N})$ in ${\Cal C}^-_z$  we define
$N_{(n)} = ({}^\omega N)_{(n)} \subset \underline{N}$  where we identify the
spaces
$\underline{N}$ and ${}^\omega \underline{N}$.  We define
$M_n\ {\dsize\buildrel{\text{def}}\over =}\ \underline{M}/M_{(n)}$,\
$N_n\ {\dsize\buildrel{\text{def}}\over =}\ \underline{N}/N_{(n)}$
and denote by $\pi$ the natural projections $\underline{M} \to M_n$ and
$\underline{N} \to N_n$.

\subheading{2.2.11}  For any  $p \in {\Bbb N}$, we define
$$
T_{(p)}\ {\buildrel{\text{def}}\over =}\
\underline{\Cal L}_{(Z\widetilde{q}^{h^\vee})^{-2}}
\Omega_{\le p} \underline{G},
$$
where  $\Omega_{\le p} \in {\bold U}$  is as in [L] 6.1.1.

\proclaim{Proposition}  a)  For any $M$ in ${}^{[0]}{\Cal C}$  and any pair
$m > n \in {\Bbb N}$  we have  $T_{(m)} M_{(n)} \subset M_{(n)}$.

b)   The induced operator  $T_n$  on $M_n$  does not depend on a choice of $m >
n$.

c)  The system $\{ T_n \in End M_n\}$  is compatible with the natural
projections $M_n \to M_{n-1}$.
\endproclaim

\demo{Proof}  a) and b) follow from [L] 6.1.1, and c) is obvious.
\enddemo
\vskip .2in

\subheading{2.2.12}  For any $M$ in  ${\Cal C}$  we define
$\widehat{M}\ {\dsize\buildrel{\text{def}}\over =}\
{\dsize\lim_{\leftarrow}} M_n$
and denote by  $\widehat{\pi}_n: \underline{\widehat{M}} \to M_n$  the natural
projection.  Let  $\widehat{M}(\infty) \subset \widehat{M}$  be the submodule
as in 2.2.1.  Since  $\widehat{M}(\infty) \in {\Cal C}^+$  we can define a
linear transformation  $T_{\widehat{M}} \in \text{End}\
\underline{\widehat{M}}$.

\proclaim{Proposition}  $\widehat{\pi}_n \circ T_{\widehat{M}} = T_n \circ
\widehat{\pi}_n$  for all $n \in {\Bbb N}$.
\endproclaim

\demo{Proof}  Fix $\widehat{m}\in \widehat{M}(\infty)$ and $n  \in {\Bbb N}$
such that  ${\bold U}_+^{\ge n} \widehat{m} = 0$.  We may assume that $p > n$.
Let  $m \in M$  be such that  $\pi_n(m) = \pi_n(\widehat{m})$.   We have
$\widehat{\pi}_n(T_{\widehat{M}} \widehat{m}) =
\widehat{\pi}_n(T_{(p)}\widehat{m})
= \pi_n(T_{(p)}m) =T_n \pi_n(m) = T_n \widehat{\pi}_n(\widehat{m})$.
Proposition 2.2.12 is proved.
\enddemo
\vskip .3in

%this file qkz2c.tex   [contains 2.3* and 2.4*]

\newpage

\pageno=27

\subheading{2.3.\   The $R$-matrix}

\subheading{2.3.1}  Let  $\Xi$  be the complex-valued function on
${}'\overline{\Lambda} \times '\overline{\Lambda}$  such that
$\Xi(\overline{\lambda}', \overline{\lambda}'') =
q^{-[\overline{\lambda}', \overline{\lambda}'']}$.

For  any  ${\bold U}$-modules  $V' = (\rho', \underline{V}'), V'' = (\rho'',
\underline{V}'')$  in  ${}^{[0]}{\Cal C}$  we denote by $\underline{\Xi}$
the endomorphism of  $\underline{V}' \otimes \underline{V}''$  which preserves
subspaces  $V'_{\lambda'} \otimes V''_{\lambda''}$  for all $\lambda',\lambda''
\in \Lambda_{{\Bbb C}^*}$  and such that
$$
\underline{\Xi}_{|V'_{\lambda'}\otimes V''_{\lambda''}} =
\Xi (\overline{\lambda}',\overline{\lambda}'')Id,
$$
where  $\overline{\lambda}'\ {\dsize{\text{def}}\over =}\ r(\lambda')$  and
$\overline{\lambda}''\ {\dsize{\text{def}}\over =}\ r(\lambda'')$.

\proclaim{{\bf 2.3.2}  Definition}  We say that a pair  $V',V''$,
$V' = (\rho', \underline{V}') \in {}^{s'}{\Cal C}$, $V'' = (\rho'',
\underline{V}'') \in {}^{s''}{\Cal C}$, is admissible if for any
$v' \in \underline{V}'$, $v'' \in \underline{V}''$  we have
$\Theta_\nu(v'' \otimes v') = 0$  for almost all $\nu \in {\Bbb Z}[I]$  where
$\Theta_\nu \in {\bold U}^{(2)}\ {\dsize\buildrel{\text{def}}\over =}\ {\bold
U} \otimes {\bold U}$
are as in [L], 4.1.1.  If  $v' = (\rho',\underline{V}'), V'' = (\rho'',
\underline{V}'')$  is an admissible pair we denote by $\Theta$ an endomorphism
of  $\underline{V}'' \otimes \underline{V}'$  such that
$$
\Theta(v'' \otimes v') = \sum_{\nu\in {\Bbb Z}[I]} \Theta_\nu(v'' \otimes
v')\quad
\text{for\ all}\quad v' \in \underline{V}',\ v'' \in \underline{V}''.
$$
\endproclaim

\demo{Remark}  If  $V', V'' \in {\Cal C}$  are such that either  $V' \in
{\Cal C}^+$  or  $V'' \in {\Cal C}^-$  then the pair $(V', V'')$  is
admissible.
\enddemo

\medpagebreak

\subheading{2.3.3}  If  $V',V''$  is an admissible pair we define linear map
$'s_{V',V''}$  and  $s_{V',V''}$  from $\underline{V}' \otimes \underline{V}''$
 to
$\underline{V}'' \otimes \underline{V}'$ as compositions:
$$
's_{V',V''}\ {\buildrel{\text{def}}\over =}\ \Theta \circ \underline{\Xi}
\circ \sigma,\qquad s_{V',V''}\ {\buildrel{\text{def}}\over =}\
{}'s_{V',V''} \circ (\underline{\Cal L}_{Z^{-1}} \otimes \underline{\Cal
L}_{Z^{-1}})
$$
where  $\sigma:  \underline{V}' \otimes \underline{V}'' \to \underline{V}''
\otimes \underline{V}'$  is the isomorphism such that $\sigma(v' \otimes v'')
= v'' \otimes v'$ for all $v' \in \underline{V}'$, $v'' \in \underline{V}'$.

\medpagebreak

\subheading{2.3.4}     Let  $\varphi^{(2)}$  be the automorphism of
${\bold U} \otimes {\bold U}$, such that
$$
\varphi^{(2)}(x \otimes y) =
\cases
x\otimes y &\quad \text{if}\ x \otimes y \in \{ K_\lambda \otimes K_\mu,
E_i \otimes 1, 1 \otimes E_i, F_i \otimes 1, 1 \otimes F_i\} \\
&\qquad\qquad\qquad  \text{for}\ \lambda,\mu \in \Lambda^\vee, i \in
\overline{I};\\
\widetilde{Z}^{-1} \otimes E_{i_0} &\quad \text{if}\ x \otimes y =
1 \otimes E_{i_0} ;\\
E_{i_0} \otimes \widetilde{Z}^{-1} &\quad \text{if}\ x \otimes y =
E_{i_0} \otimes 1 ;\\
\widetilde{Z} \otimes F_{i_0} &\quad \text{if}\ x \otimes y = 1 \otimes
F_{i_0};\\
F_{i_0} \otimes \widetilde{Z} &\quad \text{if}\ x \otimes y = F_{i_0} \otimes 1
{}.
\endcases
$$

\proclaim{Proposition}  If  $V' = (\rho', \underline{V}'), V'' =
(\rho'', \underline{V}'') \in {}^{[0]}{\Cal C}$  is an admissible pair, then
$$
's_{V',V''} \circ \varphi^{(2)} (\Delta(x)) = \Delta(x) \circ\
{}'s_{V',V''}\quad \text{for\ all}\quad x \in {\bold U}.
$$
\endproclaim

\demo{{\bf 2.3.5}\ Proof}  As in [L] 32.1, we have to show that
$\Theta \cdot \underline{\Xi} \varphi^{(2)}({}^t\Delta(x)) = \Delta(x)
\Theta \underline{\Xi}$  when $x$ runs through a system of generators of
${\bold U}$.
\enddemo

Let $\alpha$  be an automorphism of  ${\bold U}^{(2)}$  as in [L] 32.1.  Then
we have
$\Delta(x)\Theta = \Theta\alpha({}^t\Delta(x))$  for all $x \in {\bold U}$ (see
[L] 32.1).   We have to show that
$$
\varphi^{(2)}({}^t\Delta(x)) = \underline{\Xi}^{-1} \alpha({}^t\Delta(x))
\underline{\Xi}
$$
for a system of generators $x$ in ${\bold U}$.

For any  $\overline{\lambda} \in {}'\overline{\Lambda}$  we define
$$
\underline{V}'_{\overline{\lambda}}\ {\buildrel{\text{def}}\over =}\
\bigoplus_{\lambda \in r^{-1}(\overline{\lambda})}
\underline{V}'_\lambda,\quad
\underline{V}''_{\overline{\lambda}}\ {\buildrel{\text{def}}\over =}\
\bigoplus_{\lambda \in r^{-1}(\overline{\lambda})}
\underline{V}''_\lambda .
$$
Then it is sufficient to prove the equality
$$
\varphi^{(2)}({}^t\Delta(x))(v'' \otimes v') = \underline{\Xi}^{-1}
\alpha({}^t\Delta(x)) \underline{\Xi}(v'' \otimes v')
\tag *
$$
for all\ $v' \in \underline{V}'_{\overline{\lambda}'},
v'' \in \underline{V}''_{\overline{\lambda}''},\quad
\overline{\lambda'}, \overline{\lambda}'' \in \overline{\Lambda}_{{\Bbb C}^*}$
and a system of generators $x$ of ${\bold U}$.  If $x \in {\bold U}_0$, then
the
validity of (*) is obvious.  If  $x = E_i$ or $F_i$  for $i \in \overline{I}$
then the validity of (*) is proven in [L] 32.1.  So it is sufficient to prove
(*) in the case when $x = E_{i_0}$  and  $x = F_{i_0}$.  We prove (*) in the
case when  $x = E_{i_0}$.  The case  $x = F_{i_0}$  is completely analogous.

\subheading{2.3.6}  If we compute the right side of (*) for $x = E_{i_0}$  we
find that it is equal to  $q^{[i'_0,{\overline{\lambda}}']}
E_{i_0} v'' \otimes v' + \widetilde{Z}^{-1}v'' \otimes E_{i_0} v'$.  On the
other hand,
$$
\varphi^{(2)} ({}^t\Delta(E_0)) = \varphi^{(2)}(E_{i_0} \otimes K_{i_0} + 1
\otimes
E_{i_0}) = E_{i_0} \otimes \widetilde{Z}^{-1} K_{i_0} + \widetilde{Z}^{-1}
\otimes E_{i_0} = E_{i_0} \otimes K_\nu^{-1} + \widetilde{Z}^{-1} \otimes
E_{i_0}.
$$
where $\nu = {\dsize\sum_{i\in\overline{I}}} n_i {\frac{(i\cdot i)}{2}} i$.
and we see that the left side of (*) is equal to the right side.  Proposition
2.3.4
is proved.

\medpagebreak

\proclaim{{\bf 2.3.7}\ Corollary}  Let  $V', V''$  be an admissible pair of
${\bold U}$-modules such that  $V' \in {\Cal C}_{z'}, V''~\in~{\Cal C}_{z''},$\
$z',z'' \in {\Bbb C}^*$.  Then  $'s_{V',V''}: \underline{V}' \otimes
\underline{V''} \to \underline{V}'' \otimes \underline{V}'$  is an ${\bold
U}$-module
isomorphism between ${\bold U}$-modules  ${\Cal T}^\varphi_{z^{''-1}}(V')
\otimes
{\Cal T}^\varphi_{z^{'-1}}(V'')$  and  $V'' \otimes V'$.
\endproclaim

\medpagebreak

\subheading{2.3.8}  Let  $V', V''$  be as in 2.3.7.

\proclaim{Corollary}  The linear map  $s_{V',V''}$  defines a ${\bold
U}$-module
isomorphism between ${\bold U}$-modules ${\Cal T}_{z^{''-1}}(V') \otimes
{\Cal T}_{z^{'-1}}(V'')$  and  $V'' \otimes V'$.
\endproclaim

\proclaim{{\bf 2.3.9}\ Lemma} Let  $V', V''$,\ $V \in {}^{[0]}{\Cal C}$  be
such that the pairs  $(V,V'), (V,V'')$  and  $(V,V' \otimes V'')$  are
admissible.  Then
$$
s_{V,V'\otimes V''} = (Id_{V'} \otimes s_{V,V''}) \circ (s_{V,V'} \otimes
Id_{V''}).
$$
\endproclaim

\demo{Proof}  Completely analogous to the proof of Proposition 32.2.4 in [L].
\enddemo
\vskip .3in

\subheading{2.4.\ Quantum algebras over  ${\Bbb C}[t]$}

\subheading{2.4.1}  Let  $A = {\Bbb C}[t], F = {\Bbb C}(t)$, $A_n = A/t {}^nA,
n \in {\Bbb N}$.  For any ${\bold U}$-module $M$ and $n \in {\Bbb N}$  we
define
${}^nM\ {\dsize\buildrel{\text{def}}\over =}\ M \otimes_A A_n$.  Let  ${\bold
U}$
be as in 2.1.4,
$$
\aligned
{\bold U}_A\ & {\buildrel{\text{def}}\over =}\ {\bold U} \otimes_{\Bbb C}
A,\quad
{\bold U}_F\ {\buildrel{\text{def}}\over =}\ {\bold U} \otimes_{\Bbb C} F,\quad
{}^n{\bold U} = {\bold U} \otimes_{\Bbb C} A_n,\\
{\bold U}_A^{(2)}\ & {\buildrel{\text{def}}\over =}\
{\bold U}_A \otimes_A {\bold U}_A\quad
\text{and}\  \quad
{}^n{\bold U}^{(2)} = {\bold U}^{(2)}_A \otimes_A A_n .
\endaligned
$$
We will consider  ${\bold U}$  as a subalgebra in ${\bold U}_A$  and  ${\bold
U}_n$.  The
comultiplications $\Delta: {\bold U} \to {\bold U} \otimes {\bold U}$  defines
$A$-linear comultiplications
${\bold U}_A \to {\bold U}_A^{(2)}$  and  ${}^n{\bold U} \to
{}^n{\bold U}^{(2)}$  which we also denote by
$\Delta$.

\subheading{2.4.2}  Let\quad   $\psi\ {\buildrel{\text{def}}\over =}\
\psi_t$\quad
be the automorphism of\ \   ${\bold U}_F$\ \   as in 2.1.12 and

\noindent
$\Gamma\ {\buildrel{\text{def}}\over =}\ \{ x \in {\bold
U}_A|(1\otimes\psi)(\Delta(x))
\in {\bold U}_A^{(2)}\}$.  As follows from [L], Proposition 3.2.4, we have a
direct sum
decomposition  ${\bold U} = \oplus_\nu {\bold U}_\nu$, where ${\bold U}_\nu\
{\dsize\buildrel{\text{def}}\over =}\
{\bold U}_+({\bold U}_-^\nu)$.  We consider the ${\Bbb C}$-linear map  $\eta:
{\bold U} \to {\bold U}_A$  such
that  $\eta(u) = t^{|\nu|d} u$  for $u \in {\bold U}_\nu$  and extend it to an
$A$-linear map  $\eta_A:  {\bold U} \otimes_{\Bbb C} A \to {\bold U}_A$.  Here
$|\nu| = {\dsize\sum^n_{p=1}} \nu_p {\frac{(i_p\cdot i_p)}{2}}$  if
$\nu = (\nu_1,\cdots,\nu_n)$.  This is called a degree of $\nu$.

\proclaim{Proposition}  a)  $\Gamma$  is a subalgebra of  ${\bold U}_A$,

b)  Im $\eta \subset \Gamma$

c)  The map  $\eta_A: {\bold U} \otimes_{\Bbb C} A \to \Gamma$  is an
isomorphism of $A$-modules.

d)  The algebra  $\Gamma$  is generated by  ${\bold U}_+$  and
$t_i^d F_i$, $i \in I$, where  $t_i = t^{(i\cdot i)/2}$.

e)  $\Delta(\Gamma) \subset \Gamma \otimes \Gamma$.
\endproclaim

\demo{{\bf 2.4.3}\ Proof}  Part a) follows immediately from the definition
of $\Gamma$  since  $(1 \otimes \psi) \circ \Delta$  is an algebra homomorphism
from ${\bold U}_F$  to  ${\bold U}_F \otimes {\bold U}_F$.  Part b) follows
from a)  and the
inclusions  $(1 \otimes \psi) \circ \Delta({\bold U}_+) \subset (1 \otimes
\psi)
({\bold U}_+ \otimes {\bold U}_+) \subset {\bold U}_+ \otimes {\bold U}_+
\subset {\bold U}_A^{(2)}$  and
$$
(1 \otimes \psi) \circ \Delta(t^d_iF_i) = (1 \otimes \psi)(t^d_i(F_i \otimes
\widetilde{K}_{-i}) + t^d_i(1 \otimes F_i)) = t^d_i(F_i \otimes
\widetilde{K}_{-i})
+ 1 \otimes F_i \in {\bold U}_A^{(2)}.
$$
\enddemo

\subheading{2.4.4}  To prove part c) we choose a basis $B$ in ${\bold f}$
consisting of homogeneous elements and containing 1 (see [L] 3.2.4).  For any
$b \in B$  we denote by $|b|$  its degree, i.e., $|b| = |\nu|$ if
$b \in f_\nu$.  As follows from
Proposition 3.2.4 in [L] we can write  $x \in \Gamma$  as a sum
$x = \sum_{b',\mu,b} c_{b',\mu,b} b^{'+} K_\mu b^-$  with
$c_{b',\mu,b} \in A$.  We have to show that  $c_{b',\mu,b} \in t^{|b|}A$.

Since the map  $\Delta: {\bold U} \to {\bold U} \otimes {\bold U}$  is a
monomorphism, we see that the inclusion  $(1 \otimes \psi) \circ \Delta(x) \in
{\bold U}_A \otimes {\bold U}_A$ implies the inclusion
$c_{b',\mu,b}(1 \otimes \psi)\Delta (b^{'+}K_\mu b^-) \in {\bold U}_A
\otimes {\bold U}_A$  for all  $b,b' \in B, \mu \in \Lambda^\vee$.  Therefore
we may assume that  $x = ab^{'+} K_\mu b^-$  for  $b', b \in B$, $\mu \in
\Lambda^\vee$  and  $a \in A$.

\subheading{2.4.5}  Assume that  $a \notin t^{|b|}A$.  We want to show that
$(1 \otimes \psi)$ $(\Delta(x)) \notin {\bold U}_A^{(2)}$.  Choose $n \in {\Bbb
N}$  such
that $a \in t {}^nA - t {}^{n+1}A$  and  define  $\widetilde{x}\
{\dsize\buildrel{\text{def}}\over =}\ t^{|b|-n}x$.  Let
${}^-:\ {\bold U}_A^{(2)} \to {\bold U} \otimes {\bold U}$  be  the natural
projection
(= reduction mod $t$).  As follows from part b), $(1 \otimes \psi)
\Delta(\widetilde{x}) \in {\bold U}_A^{(2)}$.  Let $\overline{x}$ be the
image of this element in  ${\bold U} \otimes {\bold U}$.
Since $|b| - n > 0$  it is sufficient to show that
$\overline{x} \ne 0$.

We have
$$
(1 \otimes \psi) \circ \Delta(\widetilde{x}) = t^{|b|-n} a (1 \otimes \psi)
\Delta(b^{'+}) (\Delta(K_\mu))(1 \otimes \psi)\cdot (\Delta(b^-)).
$$
Therefore it follows from formulas in 2.1.10 and the definition of $\psi$ that
$$
\overline{x} = \overline{a} b^{'+} K_\mu \otimes K_\mu b^-
$$
where  $\overline{a} \in {\Bbb C}$  is the reduction of $at^{-n} \in A$
mod $t$.   By the definition of $n$, we have  $\overline{a} \ne 0$.
Therefore,   $\overline{x} \ne 0$.  Part c) of Proposition 2.4.2 is proved.

Part d) follows immediately from c) and part e) follows from d) and explicit
formulas for $\Delta$ (see 2.1.12). Proposition 2.4.2 is proved.

\proclaim{{\bf 2.4.6}\ Corollary}  $\Gamma =
\{ x \in {\bold U}_A|(\psi \otimes 1)(\Delta(x)) \in {\bold U}_A^{(2)}\}$.
\endproclaim

\demo{Proof}  Let  $\Gamma' = \{ x \in {\bold U}_A|(\psi \otimes 1)\Delta(x)
\in
{\bold U}_A^{(2)}\}$.  It follows from Proposition 2.4.2 that  $\Gamma \subset
\Gamma'$.
An analogous argument shows that  $\Gamma' \subset \Gamma$.  So
$\Gamma = \Gamma'$.
\enddemo

\subheading{2.4.7}   The map $\Delta: {\bold U}_0 \to {\bold U}_0 \otimes
{\bold U}_0 \subset {\bold U} \otimes {\bold U}$  defines an imbedding of
${\bold U}_0$  in  ${\bold U} \otimes {\bold U}$.  Since  ${\bold U}_0
\subset \Gamma$ and  $(1 \otimes \psi) \circ \Delta|_{{\bold U}_0} =
\Delta|_{{\bold U}_0}$  we can consider  ${\bold U}_0$  as a subalgebra of
$\Gamma$.  Therefore the imbeddings  $(1 \otimes \psi) \circ \Delta:
\Gamma \hookrightarrow {\bold U}_A^{(2)}$  and  ${\bold U}_- \otimes
{\bold U}_+ \hookrightarrow {\bold U}_A^{(2)}$  define an $A$-linear
map $\alpha: \widetilde{\bold U}_A^{(2)} \to {\bold U}_A^{(2)}$, where
$\widetilde{\bold U}_A^{(2)}\ {\buildrel{\text{def}} \over =}\
\Gamma \otimes_{{\bold U}_0} ({\bold U}_- \otimes {\bold U}_+)$.
Analogously for any  $n \in {\Bbb N}$
we can define an $A_n$-linear map $\alpha_n: {}^n\widetilde{\bold U}^{(2)}
\to\ {}^n{\bold U}^{(2)}$, where\ ${}^n{\bold U}^{(2)}\
{\buildrel{\text{def}} \over =}\ \widetilde{\bold U}_A^{(2)}/t^n
\widetilde{\bold U}_A^{(2)}$.

\proclaim{Theorem}  The map  $\alpha$  is an isomorphism.
\endproclaim

\demo{{\bf 2.4.8}\ Proof}  We start the proof with the following general
result.
\enddemo

\proclaim{Lemma}  Let $M,N$  be free $A_n$-modules,  let  $\alpha: M \to N$  be
a
morphism such that the induced map  $\overline{\alpha}: M/tM \to N/tN$ is an
isomorphism.  Then $\alpha$ is an isomorphism.
\endproclaim

\demo{Proof of Lemma}  Follows from Nakayama's lemma.
\enddemo

\proclaim{{\bf 2.4.9}\ Proposition}  The maps  $\alpha_n: n \in {\Bbb N}$  are
isomorphisms.
\endproclaim

\demo{Proof}  Consider first the case $n = 1$.  We have\
${}^1{\bold U}^{(2)} = {\bold U} \otimes {\bold U}$,\
${}^1\widetilde{{\bold U}}^{(2)} = \Gamma_1 ({\bold U}_- \otimes {\bold U}_+)$
where the subalgebra  $\Gamma_1 \subset {\bold U} \otimes {\bold U}$  is
generated
by elements   $x^+ \otimes 1, 1 \otimes x{'{}^-}$  and $K_\mu \otimes K_\mu$
for
$x,x' \in {\bold f}$  and $\mu \in \Lambda^\vee$.  Therefore in the case when
$n = 1$ Proposition 2.4.9 follows from the triangular decomposition for
${\bold U}$  (see [L] 3.2).  The general case follows now from Lemma 2.4.8.
Proposition 2.4.9 is proved.
\enddemo

\proclaim{{\bf 2.4.10}\ Corollary}  The map $\alpha$  is a monomorphism.
\endproclaim

\demo{Proof}  Let  $x \in \widetilde{{\bold U}}_A^{(2)}$  be
such that  $\alpha(x) = 0$.  It follows from Proposition 2.4.9 that the
image of  $x$ in  $t^n \widetilde{{\bold U}}_A^{(2)}$
is equal to zero for all $n \in {\Bbb N}$.  Therefore  $x$  is divisible in
$\widetilde{{\bold U}}_A^{(2)}$  by  $t^n$  for all $n \in {\Bbb N}$.
But the $A$-module $\widetilde{{\bold U}}_A^{(2)}$  is free.  Therefore
$x = 0$.  The Corollary is proved.
\enddemo

\medpagebreak

\subheading{2.4.11}  For any $m \in {\Bbb N}$  we denote by ${\bold
U}^{(2)}(m)$  the
subspace of  ${\bold U} \otimes {\bold U}$  spanned by elements of the forms
$x^+K_\mu x{'{}^-} \otimes y^-K_{\mu'}y^{'+}$, where $\mu,\mu' \in
\Lambda^\vee$  and
$x,x',y,y'$  are homogeneous elements such that  $|x| + |y| \le m$  and we
define  ${\bold U}_A^{(2)}(m) = {\bold U}^{(2)}(m) \otimes_{\Bbb C} A \subset
{\bold U}_A^{(2)}$.

\proclaim{Lemma}  For all $m \in {\Bbb N}$  we have
${\bold U}_A^{(2)}(m) \subset Im(\alpha)$.
\endproclaim

\demo{{\bf 2.4.12}\ Proof}  We will prove Lemma 2.4.11 by the induction in $m$.
If  $m = 0$,  then the result follows from the inclusion ${\bold U}^{(2)}(0)
\subset {\bold U}_- \otimes {\bold U}_+$.  Assume that the lemma is true for
$m-1$.  It is sufficient
to show that for any  $\mu,\mu', x,x', y,y'$  as in 2.4.11 there exists
$\gamma \in \Gamma$  and  $\widetilde{u} \in {\bold U}_- \otimes {\bold U}_+$
such that $x^+K_\mu x^{'-} \otimes y^-K_\mu y^{'+} -
(1 \otimes \psi)\Delta(\gamma) \widetilde{u} \in
{\bold U}_A^{(2)} (m-1)$.  But we can take  $\gamma = t^{|y|} x^+y^-$  and
$\widetilde{u} = K_{-\delta} K_\mu x^{'-} \otimes K_{-\delta '}
K_{\mu'} y^{'+}$, where $\Delta(x^+) = x^+ \otimes K_{\delta '} + \cdots,$\
$\Delta(y^-) = K_\delta \otimes y^- + \cdots .$   Lemma 2.4.12 is proved.
\enddemo

\medpagebreak

\subheading{2.4.13}   Now we can finish the proof of Theorem 2.4.7.  It follows
from the triangular decomposition that  ${\bold U}_A^{(2)} = \bigcup_m
{\bold U}_A^{(2)} (m)$.  Therefore Lemma 2.4.11 implies the surjectivity
of $\alpha$.  On the other hand, the injectivity of  $\alpha$  follows
from Corollary 2.4.10.  Theorem  2.4.7 is proved.

\subheading{2.4.14}  We will use a following version of Theorem 2.4.7.  As
follows
from Corollary 2.4.6 we have $(\psi \otimes 1)\Delta(\gamma) \in {\bold
U}_A^{(2)}$ for all
$\gamma \in \Gamma$.  Let  $\beta: \Gamma \otimes_{{\bold U}_0}({\bold U}_+
\otimes {\bold U}_-) \to
{\bold U}_A^{(2)}$  be the morphism such that $\beta(\gamma \otimes u) =
(\psi \otimes 1)  \Delta(\gamma) u$  for all $\gamma \in \Gamma$ and
$u \in {\bold U}_+ \otimes {\bold U}_-$.

\proclaim{theorem}  $\beta$  is an isomorphism of $A$-modules.
\endproclaim

\medpagebreak

\proclaim{{\bf 2.4.15}\ Corollary}  Let  $\beta^{(1)}$  be the linear map
from  ${\bold U} \otimes {\bold U}_+ \otimes {\bold U}_-$  to ${\bold U}
\otimes {\bold U}$  such that
$\beta^{(1)}(x \otimes u) = \Delta(x)u$.  Then  $\beta^{(1)}$  is an
isomorphism of linear spaces.
\endproclaim

\demo{Proof}  Let  $ev^1: A \to {\Bbb C}$  be the morphism of evaluation at
$t = 1$.   Since  $\beta^{(1)} = \beta \otimes_A {\Bbb C}$,  where  $A$  acts
on ${\Bbb C}$  by  $ev^1$.  Theorem 2.4.14 implies the validity of the
 corollary.
\enddemo

\subheading{2.4.16}  As follows from part e) of Proposition 2.4.2 $\Gamma$
has a natural structure of a Hopf $A$-algebra.  Let $\Gamma_0 \subset \Gamma$
be the kernel of the counit  $\epsilon_A: {\bold U}_A \to A$  to $\Gamma$.
We denote by  ${\bold U}_0^{(2)} \subset {\bold U}_A^{(2)}$  the span of
elements of the
form  $(\psi \otimes 1)(\Delta(\gamma))(x)$, $\gamma \in \Gamma_0$,
$x \in {\bold U}_A^{(2)}$.

\proclaim{Proposition}  For any  $x \in {\bold f}_\nu$  there exists
$y,z \in {\bold U}_0^{(2)}$  such that  $(1 \otimes x^+) - y \in t^{|\nu|}
{\bold U}_A^{(2)}$
and  $(x^- \otimes 1) - z \in t^{|\nu|} {\bold U}_A^{(2)}$.
\endproclaim

\demo{{\bf 2.4.17}\ Proof}  We prove the first part of the Proposition.  The
proof of the second part is completely analogous.

Since  ${\bold U}_+$  is a Hopf algebra, it follows from Lemma 1.3.3 that
there exists $a_r \in {\bold U}_+$, $b_r \in {\bold U}_+ \otimes
{\bold U}_+$, $1 \le r \le R$  such that
$\epsilon(a_r) = 0$  and $S^{-1}(x^+) \otimes 1 - 1 \otimes x^+ =
\sum_{1\le r\le R} \Delta(a_r)b_r$, where $S \in Isom({\bold U}_+,
{\bold U}_+^{op})$  is the antipode.  Therefore
$$
\psi(S^{-1}(x^+)) \otimes 1 - 1 \otimes x^+ =  \sum_{1\le r\le R}
(\psi \otimes 1) \circ \Delta(a_r) \cdot (\psi \otimes 1) b_r.
$$
Since  $\psi({\bold U}_+) \subset {\bold U}_+ \otimes A$, we see that
$a_r \in \Gamma$ and  $(\psi \otimes 1)b_r \in {\bold U}_A^{(2)}$  for all
$r$,\ $1 \le r \le R$.  On the other hand, the inclusion  $x \in ({\bold
U}_+)^\nu$
implies that  $S^{-1}(x) \in ({\bold U}_+)^\nu$  and therefore
$\psi(S^{-1}(x^+)) = t^{|\nu|} S^{-1}(x^+)$.  Proposition 2.4.16 is proved.
\enddemo

\subheading{2.4.18}  For any ${\bold U}$-module  $M = (\rho, \underline{M})$
the $A$-module\ $\underline{M}_A\ {\dsize\buildrel {\text{def}}\over =}\
\underline{M} \otimes_{\Bbb C} A$  has a natural  ${\bold U}_A$-module
structure.  Therefore two imbeddings id and $\psi$ of $\Gamma$ into
${\bold U}_A$  define two $\Gamma$-module structures on  $\underline{M}_A$.
We denote the first one as $M[t]_0$ (or simply $M$ if this does not create
confusion) and the second as $M[t]$.  We
define  ${}^nM\ {\dsize\buildrel{\text{def}}\over =}\ M[t]_0 \otimes_A A_n$
and  ${}^nM[t]\ {\dsize\buildrel{\text{def}}\over =}\ M[t] \otimes_A A_n$.

For any  $M' = (\rho',\underline{M}'), M'' = (\rho'', \underline{M}'')$  in
${\Cal C}$  we denote by ${}^n(M'[t] \otimes M''$) the $\Gamma$-module which
is the tensor product of\  ${}^nM'[t]$  and\ ${}^nM''$  over $A_n$.
Corresponding representation is denoted by\ ${}^n(\rho_{M'[t]} \otimes
\rho_{M''})$.

\medpagebreak

\subheading{2.4.19}  For any $\nu \in {\Bbb N}[I]$ we define
$\Theta_\nu^\psi\ {\dsize\buildrel{\text{def}}\over =}\ (1 \otimes
\psi)(\Theta_\nu)
\in {\bold U}_A^{(2)}$  where  $\Theta_\nu \in {\bold U}^{(2)}$  are as in
2.3.2.  Then
$\Theta_\nu^\psi \in t^{d|\nu|} {\bold U}_A^{(2)}$  for all $\nu \in {\Bbb
N}[I]$.
Therefore for any  $M' = (\rho_{M'}, \underline{M}')$, $M'' = (\rho_{M''},
\underline{M}'')$  in  ${}^{[0]}{\Cal C}$  and any   $n \in {\Bbb N}$  we have\
${}^n(\rho_{M''} \otimes \rho_{M'})(\Theta_\nu^\psi) = 0$  for almost all
$\nu \in {\Bbb N}[I]$  and we can define an $A_n$-linear map
${}^ns_{M'[t],M''}$  from  $\underline{M}' \otimes_{\Bbb C} \underline{M}''
\otimes_{\Bbb C} A_n$ to $\underline{M}'' \otimes_{\Bbb C} \underline{M}'
\otimes_{\Bbb C} A_n$  as a finite sum
$$
{}^ns_{M'[t],M''}\ {\buildrel{\text{def}}\over =}\ \sum_{\nu\in {\Bbb N}[I]}
{}^n(\rho_{M''} \otimes \rho_{M'})(\Theta^\psi_\nu) \underline{\Xi}
\sigma \circ (\underline{\Cal L}_{Z^{-1}} \otimes \underline{\Cal L}_{Z^{-1}})
$$
where  $\underline{\Xi}$, $\sigma$ and $\underline{\Cal L}_{Z^{-1}}$
are as in 2.3.3.

\proclaim{{\bf 2.4.20}\ Proposition}  For any  $M', M'' \in {}^{[0]}{\Cal C}$
such that  $M \in {\Cal C}_{z'}   M'' \in {\Cal C}_{z''}$,\ $z',z'' \in
{\Bbb C}^*$  and any  $n \in {\Bbb N}$  the map  ${}^ns_{M'[t],M''}$  is
a $\Gamma$-module morphism from  ${}^n({\Cal T}_{z^{''-1}}(M') [t] \otimes
{\Cal T}_{z^{'-1}}(M''))$  to  ${}^n(M'' \otimes M'[t])$.
\endproclaim

\demo{Proof}  Clear.
\enddemo

\proclaim{{\bf 2.4.21} Proposition}  For any  $M', M'', M$  in ${}^{[0]}{\Cal
C}$
and $n \in {\Bbb N}$  we have
$$
{}^ns_{M[t], M'\otimes M''} = (Id_{M'} \otimes {}^ns_{M[t],M''}) \circ
({}^ns_{M[t],M'} \otimes Id_{M''}).
$$
\endproclaim

\demo{Proof}  Analogous to the proof of Proposition 32.2.4 in [L].
\enddemo

\medpagebreak

\subheading{2.4.22}  For any  $M, N$ in ${\Cal C}$  and  $u \in {\Bbb C}^*$
we denote by $\widetilde{u}$  the automorphism of $A$ such that
$\widetilde{u}(t)
= ut$  and by $\widehat{u} = \widehat{u}_{M,N}$  the ${\Bbb C}$-linear
endomorphism of the space  $\underline{M} \otimes_{\Bbb C} \underline{N}
\otimes_{\Bbb C} A$  such that $\widehat{u}(m \otimes n \otimes a) =
m \otimes n \otimes \widetilde{u}(a)$  for all $m \in \underline{M},
n \in \underline{N}, a \in A$.  It is clear that  $\widehat{u}$  is a
$\widetilde{u}$-linear automorphism of  $\underline{M} \otimes_{\Bbb C}
\underline{N} \otimes A$.

\proclaim{Proposition}  For any  $M, N$ in ${\Cal C}$  and  $u \in {\Bbb C}^*$
we have  $\widehat{u}(M[t] \otimes N)_{(0)} \subset ({\Cal T}_u(M)[t] \otimes
N)_{(0)}$, where the subspaces  $(M[t] \otimes N)_{(0)}$ and  $({\Cal
T}_u(M)[t]
\otimes N)_{(0)}$  of  $\underline{M} \otimes \underline{N} \otimes A$  are as
in 1.1.6.
\endproclaim

\demo{Proof}  Follows immediately from the definitions.
\enddemo

\proclaim{Corollary}  The  $\widetilde{u}$-morphism  $\widehat{u}$  defines
a $\widetilde{u}$-linear isomorphism between  $A$-modules

\noindent
$\langle M[t] \otimes N\rangle\ {\dsize\buildrel\sim\over\longrightarrow}\
\langle {\Cal T}_u(M)[t] \otimes N\rangle$, where  $\langle V\rangle$  denote
$\langle V\rangle_\Gamma$.
\endproclaim

\demo{Proof}  Follows from Proposition 2.4.22.
\enddemo

We denote the induced isomorphisms from $\langle M[t] \otimes N\rangle$  to
$\langle {\Cal T}_u(M)[t] \otimes N\rangle$  also by $\widehat{u}$.
\vskip .3in

%this file is qkz2d.tex	 [contains 2.5*]

\newpage

\subheading{2.5\ $\Gamma$-coinvariants}

\subheading{2.5.1}  For any $\Gamma$-module $M$ we denote by
$\langle M\rangle_\Gamma$ (or simply $\langle M\rangle$) the $A$-module of
$\Gamma$-coinvariants (see 1.1.6).

\subheading{2.5.2}  Let ${\Cal M} = (\rho_{\Cal M}, \underline{\Cal M})$
be a ${\bold U}_+$-module,  ${\Cal N} = (\rho_{\Cal N}, \underline{\Cal N})$
be a  ${\bold U}_-$-module and  $X = (\rho_X, \underline{X})$  a
${\bold U}$-module.  Then $\underline{\Cal M} \otimes \underline{X}
\otimes \underline{\Cal N}$
has a natural structure of a ${\bold U}_0$-module and we denote by
$\langle {\Cal M} \otimes X \otimes {\Cal N}\rangle_{{\bold U}_0}$
the space of ${\bold U}_0$-coinvariants of ${\Cal M} \otimes X \otimes
{\Cal N}$.

Let  $M\ {\dsize\buildrel{\text{def}}\over =}\ {\bold U}
\otimes_{{\bold U}_+}{\Cal M}$
and  $N\ {\dsize\buildrel{\text{def}}\over =}\ {\bold U}
\otimes_{{\bold U}_-}{\Cal N}$
be the induced  ${\bold U}$-modules and let $\underline{j}:
(\underline{\Cal M} \otimes
\underline{\Cal N}) \to \underline{M} \otimes \underline{N} \subset
\underline{M}[t] \otimes \underline{N}$  be the natural imbedding, that is,
$\underline{j}(m \otimes n) = (1 \otimes m) \otimes (1 \otimes n)$ for
$m \in \underline{\Cal M}$, $n \in \underline{\Cal N}$.  Then $\underline{j}$
is a morphism of ${\bold U}_0$-modules and it induces a morphism of
$A$-modules
$$
j: \langle {\Cal M} \otimes {\Cal N}\rangle_{{\bold U}_0} \otimes A
\longrightarrow
\langle M[t] \otimes N\rangle.
$$

\proclaim{Proposition}  The map $j$ is an isomorphism.
\endproclaim

\demo{Proof}  Follows from Theorem 2.4.15.
\enddemo

\proclaim{Corollary}  For any  $n \in {\Bbb N}$  the natural map
$$
\langle (M[t] \otimes N) \otimes _A A_n\rangle \to \langle {}^n(M[t] \otimes
N)\rangle
$$
is an isomorphism.
\endproclaim

\demo{Proof}  Clear.
\enddemo

\proclaim{{\bf 2.5.3}\ Lemma}  The map $\langle  {\Cal M} \otimes {\Cal
N}\rangle_{{\bold U}_0}
\longrightarrow \langle M \otimes N\rangle_{\bold U}$  induced by the imbedding
$\underline{j}$  is an isomorphism.
\endproclaim

\demo{Proof}   Using the same arguments as in the proof of Corollary 2.4.15 we
deduce Lemma 2.5.3 from Proposition 2.5.2.
\enddemo

\subheading{2.5.4}  Let  $V,W$  be ${\bold U}$-modules such that there exist
exact
sequences  $M_1 \to M_0 \to V \to 0$  and  $N_1 \to N_0 \to W \to 0$  of
${\bold U}$-modules such that  $M_i = {\bold U} \otimes_{{\bold U}_+} {\Cal
M}_i$,
$N_i = {\bold U} \otimes_{{\bold U}_-} {\Cal N}_i$,   $i = 0,1$, where
${\Cal M}_i$  (resp. ${\Cal N}_i$) are  ${\bold U}_+$  (resp. ${\bold U}_-$)
modules.

\proclaim{Proposition}  For any  $n \in {\Bbb N}$  the natural
$A_n$-morphism

\noindent
$\langle V[t] \otimes W\rangle \otimes_A A_n \longrightarrow
\langle {}^n(V[t] \otimes W)\rangle$  is an isomorphism.
\endproclaim

\demo{Proof}  The functor  $\otimes_{\Bbb C}$  is exact and therefore
the sequence of
$$
(M_1[t] \otimes N_0) \oplus (M_0[t] \otimes N_1) \to M_0[t] \otimes N_0
\to V[t] \otimes W \to 0   \tag *
$$
of $\Gamma$-modules is exact.  Since  $\langle\ ,\ \rangle$  is a right exact
functor we see that the sequence
$$
\langle M_1[t] \otimes N_0\rangle \oplus \langle M_0[t] \otimes N_1\rangle
\longrightarrow \langle M_0[t] \otimes N_0\rangle \longrightarrow
\langle V[t] \otimes W\rangle \longrightarrow 0
$$
is exact.   Analogously one shows that for any $n \in {\Bbb N}$  the sequence
$$
\langle {}^n(M_1[t] \otimes N_0)\rangle \oplus \langle
{}^n(M_0[t] \otimes N_1)\rangle
\longrightarrow \langle {}^n(M_0[t] \otimes N_0)\rangle \longrightarrow
\langle {}^n(V[t] \otimes W)\rangle \longrightarrow 0
$$
is also exact.  Therefore Proposition 2.5.4 follows from Corollary 2.5.2 and
the
five-homomorphism lemma.

\subheading{2.5.5}  Assume that $V,W$  are ${\bold U}$-modules as in 2.5.4 and
dim
${\Cal M}_k < \infty$, dim ${\Cal N}_k < \infty$  for $k = 0,1$.

\proclaim{Lemma}  $\langle V[t] \otimes W\rangle$  is a finitely generated
$A$-module.
\endproclaim

\demo{Proof} Follows from the exactness of the sequence (*) and
Proposition 2.5.2.
\enddemo

\subheading{2.5.6}  Let  $D$  be the rigid category of finite-dimensional
${\bold U}$-modules $X$ such that $Z$ acts trivially on $X$.  Then ${\Cal D}$
is a
subcategory of  ${}^{[0]}{\Cal C}$.

Given $z \in {\Bbb C}^*$,\ $V \in {}^{[0]}{\Cal C}^+_z, X \in {\Cal D}$  and
$W \in
{}^{[0]}{\Cal C}^-_{z^{-1}}$  we denote by\  ${}^n\delta_{V,X,W}$  the
$A_n$-linear map from  $\langle {}^n((V \otimes X)[t] \otimes W)\rangle$ to
$\langle {}^n((V \otimes {\Cal T}_{z^2\widetilde{q}^{2h^\vee}}(X)[t]
\otimes W)\rangle$  defined as the composition
$$
\aligned
&\langle {}^n((V \otimes X)[t] \otimes W)\rangle\
{\buildrel{s_{V,X} \otimes id_W}\over\llongrightarrow}\
\langle {}^n(({\Cal T}_z (X) \otimes V)[t] \otimes W)\rangle = \\
& = \langle {}^n(({\Cal T}_z (X)[t]) \otimes_{A_n}
{}^nV[t] \otimes_{A_n} {}^nW_A\rangle\
{\buildrel {\beta}\over\longrightarrow} \\
& \longrightarrow \langle {}^n V[t] \otimes_{A_n} {}^nW_A  \otimes_{A_n}
{}^n({\Cal T}_{z\widetilde{q}^{2h^\vee}} (X) [t])\rangle\
{\buildrel {id_{{}^nV[t]} \otimes
s^{-1}_{{\Cal T}_{z\widetilde{q}^{2h^\vee}}(X),W}}\over\llongrightarrow}\\
& \langle {}^n((V[t] \otimes_{A_n}
{}^n({\Cal T}(X)[t]) \otimes_{A_n}
{}^n W_A\rangle =
\langle {}^n((V \otimes {\Cal T}(X))[t]
\otimes W)\rangle ,
\endaligned
$$
where  $\beta\ {\dsize\buildrel{\text{def}}\over =}
\beta^{{}^nV[t] \otimes_{A_n} {}^nW_A}_{{}^n{\Cal T}_z(X)[t]}$  is as
in 1.1.15  and  ${\Cal T}\ {\buildrel{\text{def}}\over =}\
{\Cal T}_{(z\widetilde{q}^{h^\vee})^2}$.

\demo{Remark}  The map  ${}^n\delta_{V,X,W}$  is an analog of the map
$\delta_a^{(3)}$  in 1.2.7.
\enddemo

\subheading{2.5.7}  We have assumed that $W$ lies in  ${}^{[0]}{\Cal
C}^-_{z^{-1}}$.
Then the composition

\noindent
$(T^{-1} \otimes id \otimes \check{T}^{-1}) \circ
{}^n\delta_{V,X,W}$  defines a $\Gamma$-module morphism $\varphi_{V,X,W}$ from
$\langle {}^n(V \otimes X)[t] \otimes W\rangle$ to
$$
\langle {}^n(({\Cal T}(V)  \otimes
{}^n(({\Cal T}(X))[t] \otimes
{\Cal T}(W))\rangle  =
\langle {\Cal T}({}^n((V \otimes X)[t] \otimes W))\rangle.
$$
As follows from 1.1.3 we can identify the $A_n$-module
$\langle {\Cal T}({}^n((V \otimes X)[t] \otimes W))\rangle$
with $\langle {}^n((V \otimes X)[t] \otimes W)\rangle$.  Therefore we can
consider $\varphi_{V,X,W}$  as an endomorphism of the $A_n$-module
$\langle {}^n((V \otimes X)[t] \otimes W\rangle$.

\proclaim{Theorem}  $\varphi_{V,X,W} = id$.
\endproclaim

We will prove this theorem at the end of section 3.2.

%this file is qkz3a.tex	  [contains section 3.1 and section 3.2]

\newpage

\pageno=38

\centerline{\bf \S3.\  The category of smooth representation}
\vskip .3in

\subheading{3.1\  The category ${\Cal O}^+_z$}

\proclaim{{\bf 3.1.1}\ Definition}  A finite dimensional
${\bold U}^f_+$-module ${\Cal N}$ is called a nil-module if there exists
$n \in {\Bbb N}$  such that the  ${\bold U}^{\ge n}_+{\Cal N} = 0$  and
${\Cal N} = \oplus_{\lambda \in \Lambda_{{\Bbb C}^*}} {\Cal N}_\lambda$ (see
2.1.5).
\endproclaim

For any nil-module ${\Cal N}$  and a number $z \in
{\Bbb C}^*$  we extend the action of  ${\bold U}^f_+$  on ${\Cal N}$  to an
action of  ${\bold U}_+$  on ${\Cal N}$  in such a way that  $Z$
acts as  $z Id$  and  define
$$
{\Cal N}^z\ {\buildrel{\text{def}}\over=}\ {{\bold U}}
\otimes_{{\bold U}_+}{\Cal N}.
$$
The map  $n \to 1 \otimes n$  defines a ${\bold U}_0$-covariant imbedding
${\Cal N} \hookrightarrow {\Cal N}^z$  and we will always consider  ${\Cal N}$
as a  ${\bold U}_0$-submodule of ${\Cal N}^z$.

\proclaim{{\bf 3.1.2}\ Definition} ${\Cal N}^z$  is called a generalized
Verma module.
\endproclaim

\demo{Remark} ${\Cal N}^z$  lies in  ${\Cal C}^+_z$.
\enddemo

For any  $a \in \overline{\Lambda}_{{\Bbb C}^*}$  we denote by  ${\Cal V}_a$
the
one-dimensional representation of ${\bold U}^f_+$  such that
$K_\mu$  acts as a multiplication by  $\langle \mu,a\rangle$  for all
$\mu \in \overline{\Lambda}^\vee$  and  ${\bold U}^>_+$  acting as zero.
In this case the  ${\bold U}$-module
${\Cal N}^z$  is denoted by  ${\bold V}^z_a$; it is called a Verma module.

\proclaim{{\bf 3.1.3}\ Proposition}  a)  For any generalized Verma module
${\Cal N}^z$  there
exists a finite filtration by submodules such that the successive quotients are
Verma modules.

b) Let  $V$  be an object in  ${\Cal C}_z^+$.
Then $V$ is a quotient of a generalized Verma module if and only if there
exists $n \ge 1$  such that  dim~$V(n) < \infty$  and  $V(n)$  generates $V$
as a ${\bold U}$-module.

c)  Any object $V$ in ${\Cal C}^+$  is a union of subobjects which are
isomorphic to quotients of generalized Verma modules`.
\endproclaim

\demo{Proof}   Analogous to the proof of Proposition 2.5 in [KL].
\enddemo
\vskip .2in

\subheading{3.1.4}   The following result concerns the action of the Sugarawa
operator $T^a$ on a Verma module  ${\bold V}^z_a$ (see 2.2.4).  For
any  $\ell \in {\Bbb C}^*$  we set
${}_\ell {\bold V}^z_a = \{ v \in {\bold V}^z_a\big\vert
T^a(v) = \ell v\}$.

\proclaim{Proposition}

a)  ${\bold V}^z_a = {\dsize\bigoplus_{\ell\in {\Bbb C}^*}}
{}_\ell {\bold V}^z_a$.

b)  If $|\widetilde{q}^{h^\vee}z| \ne 1$, then  ${}_\ell {\bold V}^z_a$  is a
finite-dimensional vector space, for all $\ell \in {\Bbb C}$.

c)  If  ${}_\ell {\bold V}^z_a \ne 0$,\ then
$\ell = (z^dq)^{2n}$ for some  $n \in {\Bbb N}$.

d)  ${}_1{\bold V}^z_a = {\Cal V}_a$.

e)  ${\bold V}^z_a = {\bold V}^z_a(\infty)$.

f)  The  ${\bold U}$-module  ${\bold V}^z_a$  has unique irreducible quotient
$L^z_a$.

g)  If  $|\widetilde{q}^{h^\vee}z| > 1$ then  ${\bold V}_a^z$  has finite
length for any  $a \in \overline{\Lambda}_{{\Bbb C}^*}$.
\endproclaim

\demo{Proof}  The proof of a) - f) is analogous to the proofs of
Propositions 2.7 - 2.9 in [KL] when one uses Proposition 2.2.7.  In order to
prove g) one should observe that if  ${\bold V}_{a'}^z \to {\bold V}_a^z$  is a
non-trivial homomorphism of ${\bold U}$-modules then
${\Cal L}^a_{(z\widetilde{q}^{h^\vee})^{-2}} (a') G^a(a') =
(z\widetilde{q}^{h^\vee})^{2nd}$  for such $n \in {\Bbb N}$  that
$a - a' = {\dsize\sum^n_{p=1}} i'_p$.  Under our assumptions the LHS is
bounded and the RHS increases when $n \to \infty$.  This implies that there
are finitely many such $a'$  for a given $a$.  Then one can finish the proof
along the lines of [KL], 2.22.
\enddemo
\vskip .2in

\proclaim{Definition}\ a)\ Complex numbers  $z_1,\cdots,z_n \in {\Bbb C}^*$
are multiplicatively independent if for any $\vec{r} = (r_1,\cdots,r_n) \in
{\Bbb Z}^n$, $\vec{r} \ne 0$ we have  $z^r_1 \cdots z_n^{r_n} \ne 1$.

\noindent
b)\ A pair $(s,z)$, $s \in {\frak S}$, $z \in {\Bbb C}^*$
is generic if for any  $a \in \Pi^{-1}(s)$  the complex numbers
$a(i)\ {\dsize\buildrel{\text{def}}\over =}\ \langle i,a\rangle$,
$i \in \overline{I}, z$ and $q$  are multiplicatively independent.

\noindent
c)\  A nil-module ${\Cal N}$  is $z$-generic if for any
$\overline{\lambda} \in \overline{\Lambda}_{{\Bbb C}^*}$  such that
${\Cal N}_{\overline{\lambda}} \ne \{ 0\}$  the pair
$(\overline{\lambda},z)$  is generic.

\proclaim{{\bf 3.1.6}\ Proposition}  a)  If a pair  $(s,z) \in {\frak S}
\times {\Bbb C}^*$  is generic, then the Verma module  ${\bold V}^z_a$
is irreducible for $a \in \Pi^{-1}(s)$.

b)  If a nil-module ${\Cal N}$  is $z$-generic then the
corresponding generalized Verma module ${\Cal N}^z$  is a direct sum of Verma
modules.
\endproclaim

\demo{Proof}  Follows from the same arguments as Propositions 9.9 and 9.10 in
[K].
\enddemo

\proclaim{{\bf 3.1.7}\ Definition}  We denote by  ${\frak Q}_{\ge 0} \subset
{\Bbb C}^*$
the set of numbers such that  $|z| \le 1$.
\vskip .2in

In the remainder of this paper we assume  $z \in {\Bbb C}^*$  is such
that $\widetilde{q}^{h^\vee} z \notin {\frak Q}_{\ge 0}$.

\proclaim{{\bf 3.1.8}\ Proposition}  Let  $V$  be an object in  ${\Cal C}^+_z$.
The following conditions are equivalent:

\item{a)}   There exists a finite composition series of $V$ with subquotients
of the form  $L^z_a$  for various $a \in {\Bbb N}^I$.

\item{b)}  $V$  is a quotient of a generalized Verma module.

\item{c)}  There exists $n \ge 1$  such that  $V(n)$ generates  $V$  as a
${\bold U}$-module  and  $dim V(n) < \infty$.

\item{d)}  $dim V(1) < \infty$.
\endproclaim

\demo{Proof}
Analogous to the proof of Theorems 2.22 and 3.2 in [KL].
\enddemo

\proclaim{{\bf 3.1.9}\ Definition}  ${\Cal O}^+_z$ is the full subcategory
of  ${\Cal C}^+_z$  consisting of modules satisfying the conditions of
Proposition 3.1.8.
\endproclaim

\proclaim{{\bf 3.1.10}\ Corollary}  Any object in  ${\Cal C}^+_z$  can be
represented as an
inductive limit of objects from ${\Cal O}^+_z$.
\endproclaim
\vskip .2in

\subheading{3.1.11} In this subsection ${\Cal A}$ denotes the ring of regular
functions on\quad
$\overline{\Lambda}_{{\Bbb C}^*} \times	{\Bbb C}^*$,\

\noindent
${}_{\Cal A}{\bold U}\ {\dsize\buildrel{\text{def}}\over =}\
{\bold U} \otimes_{\Bbb C} {\Cal A}$,\
${}_{\Cal A}{\bold U}_+\ =  {\bold U}_+ \otimes_{\Bbb C} {\Cal A}$,\ etc.  We
denote by
$\dot{{\Cal A}}$  the
${}_{\Cal A}{{\bold U}}_+$-module which is isomorphic to ${\Cal A}$ as an
${\Cal A}$-module and
such that  ${}_{\Cal A}{\bold U}^>_+$  acts trivially, $K_\mu$  acts as a
multiplication by the function $\langle \mu, \overline{\lambda}\rangle$  if
$\mu \in \Lambda^\vee_0$  and $Z$  acts as a multiplication by  $z$,
where  $(\overline{\lambda},z)$  are natural coordinates on
$\overline{\Lambda}_{{\Bbb C}^*} \times {\Bbb C}^*$.  We denote by ${\Cal V}$
the \
induced  ${}_{\Cal A}{\bold U}$-module
$$
{\Cal V} = {}_{\Cal A} {\bold U} \otimes_{{}_{\Cal A}{{\bold U}}_+} \dot{{\Cal
A}} .
$$
\vskip .2in

\subheading{3.1.12}  For any  $z_0 \in {\Bbb C}^*$  we define by
$ev_{z_0}: {\Cal A} \to {\Bbb C}^*$  the homomorphism of the evaluation at the
point  $(0,z_0) \in \overline{\Lambda}_{{\Bbb C}^*} \times {\Bbb C}^*$.  This
homomorphism defines an  algebra homomorphism
$ev_{z_0}: {}_{\Cal A}{\bold U} \to  {\bold U}$.
Given a ${\bold U}^f_+$-nil-module  ${\Cal N}$  and  $z_0 \in {\Bbb C}^*$  we
can use
$ev_{z_0}$  to define a structure of  ${}_{\Cal A}{{\bold U}}_+$-module on
${\Cal N}$.  We denote by  ${}_{\Cal A}{\Cal N}$ the
${}_{\Cal A}{{\bold U}}_+$-module
which is defined as a tensor product ${}_{\Cal A}{\Cal N} = {\Cal N}
\otimes_{\Bbb C}
{\Cal V}$.  Let  ${}_{\Cal A} \check{\Cal N}$  be the induced module
${}_{\Cal A} \check{\Cal N}\ {\dsize\buildrel{\text{def}}\over =}\
{}_{\Cal A}{\bold U} \otimes_{{}_{\Cal A}{\bold U}_+} {}_{\Cal A}{\Cal N}$.

For any point  $(\overline{\lambda},z) \in \overline{\Lambda}_{{\Bbb C}^*}
\times
{\Bbb C}^*$  we denote by   ${\frak m}_{\overline{\lambda},z} \subset {\Cal A}$
 the
maximal ideal of functions equal to zero at  $(\overline{\lambda},z)$.  Define
${}_{(\overline{\lambda},z)} \check{\Cal N}\ {\dsize\buildrel{\text{def}}\over
=}\
{}_{\Cal A}\check{\Cal N}/{\frak m}_{\overline{\lambda},z}$.  Then
${}_{(\overline{\lambda},z)} \check{\Cal N}$ has a natural structure of a
${\bold U}$-module.

\proclaim{{\bf 3.1.13}\ Proposition}   a) ${}_{\Cal {\Cal A}} \check{\Cal N}$
is a free ${\Cal A}$-module.

b)  For all $(\overline{\lambda},z) \in \overline{\Lambda}_{{\Bbb C}^*} \times
{\Bbb C}^*$  the ${\bold U}$-module  ${}_{(\overline{\lambda},z)}\check{\Cal
N}$\
is a generalized Verma module obtained from the nil-module
${}_{(\overline{\lambda},z)}{\Cal N}\ {\buildrel{\text{def}}\over =}\
{}_{\Cal A} {\Cal N}/{\frak m}_{\overline{\lambda},z_{\Cal A}}{\Cal N}$.

c)  For almost all  $(\overline{\lambda},z) \in \overline{\Lambda}_{{\Bbb C}^*}
\times
{\Bbb C}^*$  the module  ${}_{(\overline{\lambda},z)}\check{\Cal N}$\ is a
direct sum of Verma modules.
\endproclaim

\demo{Proof}  a)  and b)  follow from definitions and  c) follows from
Proposition 3.1.6 b).
\enddemo

\proclaim{{\bf 3.1.14}\ Definition}  We denote by ${\Cal O}^-_z \subset
{\Cal C}^-_{z^{-1}}$  the category of ${\bold U}$-modules $M$ such that
${}^\omega M$  lies in ${\Cal O}^+_z$.
\endproclaim

For any nil-module ${\Cal N}$ and $z \in {\Bbb C}^*$  we
define  ${\Cal N}^z_-\ {\dsize\buildrel{\text{def}}\over =}\
{}^\omega ({\Cal N}^z)$.  Then ${\Cal N}^z_-$  lies in ${\Cal C}^-$ and, as
before, we have a nature imbedding  ${}^\omega {\Cal N} \hookrightarrow
{\Cal N}^z_-$  of  ${\bold U}_0$-modules where we identify
${}^\omega {\Cal N}$  with ${\Cal N}$  as a vector space and the action of
${\bold U}_0$  on  ${}^\omega {\Cal N}$  is given by the map  $(x,n) \to
\omega(x)n$,  $x \in {\bold U}_0, n \in {\Cal N}$.
\vskip .2in

\subheading{3.1.15}  For any Hopf algebra $H$ and two $H$-modules
$M = (\rho_M, \underline{M})$ and $N = (\rho_N, \underline{N})$ we
define  $\langle M,N\rangle_H\ {\dsize\buildrel{\text{def}}\over =}\
\langle M \otimes N\rangle$  (see 1.1.6)
and denote by   $pr_{M,N}$  (or simply $pr$)
the natural projection  $pr: \underline{M} \otimes \underline{N} \to
\langle M,N\rangle_H$.
\vskip .2in

\subheading{3.1.16}  Let  ${\Cal M},{\Cal N}$  be nil-modules $z \in {\Bbb
C}^*$
and  ${\Cal M} \otimes {}^\omega {\Cal N}\ {\dsize\buildrel
i\over\hookrightarrow}\
{\Cal M}^z \otimes {\Cal N}^z_-$  be the natural imbedding.  Then we have a
linear map  $pr \circ i: ({\Cal M} \otimes {}^\omega {\Cal N}) \longrightarrow
\langle {\Cal M}^z, {\Cal N}^z_-\rangle_{{\bold U}}$  which  factorizes
through the map  $\overline{i}: \langle {\Cal M},
{}^\omega {\Cal N}\rangle_{{\bold U}_0} \longrightarrow
\langle {\Cal M}^z,{\Cal N}^z_-\rangle_{{\bold U}}$.

\proclaim{Proposition}  The map $\overline{i}$  is an isomorphism.
\endproclaim

\demo{Proof}  Using Theorem 2.4.15 one can immediately apply the
arguments of the proof of Proposition 9.15 in [KL].
\enddemo
\vskip .3in

\subheading{3.2\  The action of Sugawara operators on coinvariants}

\subheading{3.2.1}  In this section we prove Theorem 2.5.7.  We start with
the special case when $X = \ident$.  So  $V = (\rho_V, \underline{V})$,
$W = (\rho_W, \underline{W})$  be ${\bold U}$-modules such that  $V \in
{\Cal C}^+_z$ and  $W \in {\Cal C}^-_{z^{-1}}$.  Let  $T \otimes \check{T} \in
\text{End}(\underline{V}
\otimes \underline{W})$  be the linear map as in 2.2.8.  As follows from
Proposition 2.2.8 the  $T \otimes \check{T}$  induces an endomorphism of the
space  $\langle V,W\rangle$  which we denote as  $\varphi_{V,W}$.

\proclaim{{\bf 3.2.2}\ Theorem}  $\varphi_{V,W} = id$.
\endproclaim

\demo{Proof of Theorem 3.2.2}  As follows from Proposition 3.1.3 c) it is
sufficient to consider the case when  $V$ and $W$ are generalized Verma
modules,
$V = {\Cal M}^z$, $W = {\Cal N}^w_-$  where ${\Cal M}$  and ${\Cal N}$  are
nil-modules.  It is clear that  $\langle V,W\rangle = \{ 0\}$  if $z \ne w$.
\enddemo

We start with the following special case.

\proclaim{{\bf 3.2.3}\ Proposition}  Theorem 3.2.2 is true in the
case when $V$ and $W$ are Verma modules.
\endproclaim

\demo{Proof of Proposition 3.2.3}  Let  $V = {\bold V}^z_a$,
$W = {}^\omega ({\bold V}^z_b)$.  As follows from Proposition 3.1.16 the map
$\overline{i}:  \langle {\Cal V}_a \otimes
{}^\omega ({\Cal V}_b)\rangle_{{\bold U}_0} \longrightarrow
\langle V,W\rangle_{{\bold U}}$  is an isomorphism.  On the other hand,
it follows from Proposition 2.2.9 that the operator  $T \otimes \check{T}
= T^a \otimes \check{T}^a$  preserves the subspace
${\Bbb C}_a \otimes {}^\omega ({\Bbb C}_b) \subset \underline{M} \otimes
\underline{N}$ 	and acts trivially on this subspace.   Proposition 3.2.3 is
proved.
\enddemo

\subheading{3.2.4}  Consider now the case when $V$ and $W$ are arbitrary
generalized  Verma modules.  Let
${\Cal M}, {\Cal N}$  be nil-modules, $z_0 \in {\Bbb C}^*$  and
${}_{\Cal A}\check{\Cal M}, {}_{\Cal A}\check{\Cal N}$  be modules as in
3.1.12.
As follows from Proposition 2.2.8, the
operator  $T \otimes \check{T}$ defines an endomorphism of the ${\Cal
A}$-module
$\langle {\Cal M}, {\Cal N}^-\rangle\ {\buildrel{\text{def}}\over =}\
\langle {}_{\Cal A}\check{\Cal M} \otimes_{\Cal A} {}_{\Cal A}\check{\Cal
N}^-\rangle$.  We denote
this endomorphism of  $\langle {\Cal M}, {\Cal N}\rangle$ by $\Phi$.  For any
$(\overline{\lambda},z) \in \Lambda_{{\Bbb C}^*} \times
{\Bbb C}^*$  we denote by  $\Phi_{\overline{\lambda},z}$  the natural
morphism

\noindent
$\Phi_{\overline{\lambda},z}: \langle {\Cal M},
{\Cal N}^-\rangle \to \langle {}_{(\overline{\lambda},z)} {\Cal M},
{}_{(\overline{\lambda},z)} {\Cal N}^-\rangle$.

\proclaim{{\bf 3.2.5}\ Lemma} a)  The ${\Cal A}$-module $\langle {\Cal M},
{\Cal N}^-\rangle$
is free as an ${\Cal A}$-module,

b)  $\Phi_{\overline{\lambda},z}$  is surjective
 for all $\overline{\lambda} \in \overline{\Lambda}_{{\Bbb C}^*}$, $z \in {\Bbb
C}^*$,

c)  the kernel of  $\Phi_{\overline{\lambda},z}$  is equal to
${\frak m}_{\overline{\lambda}, z}\langle {\Cal M}, {\Cal N}^-\rangle$.
\endproclaim

\demo{Proof}  Follows from Propositions 3.1.13 and 3.1.16.
\enddemo

Now we can finish the proof of Theorem 3.2.2 in the case when $V$ and
$W$ are generalized Verma modules.  Really since the morphism
$\Phi_{0,1}$  is surjective it is sufficient to prove that
$\Phi = Id$.  On the other hand, since the ${\Cal A}$-module  $\langle {\Cal
M},
{\Cal N}^-\rangle$  is free, it is sufficient to show that the induced
endomorphism  $\Phi_{\overline{\lambda},z}$  on
$\langle {}_{(\overline{\lambda},z)} {\Cal M},
{}_{(\overline{\lambda},z)} {\Cal N}\rangle$  is equal to Id for generic
$(\overline{\lambda},z)$  .  But this follows from Propositions 3.1.13\ c) and
3.2.3.

Theorem 3.2.2 is proved.

\proclaim{{\bf 3.2.6}\ Proposition}  Theorem 2.5.7 is true in the case when
$V$ and $W$ are Verma modules.
\endproclaim

\demo{Proof}  We have  $V = {\bold V}^z_a, W = {}^\omega ({\bold V}^z_b)$  for
some
$a,b \in\ {}'\overline{\Lambda}$.   As follows from Theorem 2.4.7 the natural
imbedding  ${\Cal V}_a \otimes \underline{X} \otimes {}^\omega ({\Cal V}_b)
\hookrightarrow (V \otimes X)[t] \otimes W$  defines an isomorphism
$\langle {\Cal V}_a \otimes \underline{X} \otimes {}^\omega ({\Cal V}_b)
\rangle_{{\bold U}_0}
\otimes_{\Bbb C}  A_n\ {\dsize\buildrel\sim\over\longrightarrow}\
\langle {}^n((V \otimes X)[t] \otimes W)\rangle$ (cf. 2.5.2).  So it is
sufficient to
show that for any  $x \in \underline{X}_{b-a}$  we have
$\varphi_{V,X,W}(1_a \otimes x \otimes 1_b) = 1_a \otimes x \otimes 1_b$, where
$1_a,1_b$  are generators of 1-dimensional spaces  ${\Cal V}_a$  and
${\Cal V}_b$  and we identify  $({\Cal V}_a \otimes \underline{X}
\otimes {}^\omega ({\Cal V}_b))^{{\bold U}_0}$  with its image in
$\langle {\Cal V}_a \otimes \underline{X} \otimes {}^\omega ({\Cal
V}_b)\rangle$.
But this follows immediately from the definitions.  Proposition 3.2.6 is
proved.
\vskip .2in

\subheading{3.2.7}  The same arguments as in 3.2.4-3.2.5 show that
Theorem 2.5.7 follows from Proposition 3.2.6.  Theorem 2.5.7 is proved.
\vskip .3in

%this file is qkz3b.tex	 [contains section 3.3]

\subheading{3.3 Completions}

\subheading{3.3.1}  We  assume until the end of the
section all our infinite-dimensional modules are in  $^{[0]}{\Cal C}$.
It is easy to see that for any $N$ in $^{[0]}{\Cal C}$  and any finite
dimensional
${\bold U}$-module $V$ the tensor product  $N \otimes V$  lies in
$^{[0]}{\Cal C}$.
\vskip .2in

\subheading{3.3.2}  We fix until the end of this section a number
$z$ such that $zq^{h^\vee} \in {\Bbb C}^* - {\frak Q}_{\ge 0}$.  For any $M =
(\rho, \underline{M})$ in
${\Cal C}_z$  we define the spaces  $M_n$  and  $\underline{M}_{(n)}$  as in
2.2.10.

%we denote
%by $M_{(n)} \subset M$  the subspace generated by vectors of the norm
%$xm$, $x \in {\bold U}^{\geqslant n}_-$, $m \in \underline{M}$.
%For any $N = (\rho, \underline{N})$ in ${\Cal C}^-_z$  we define
%$N_{(n)} = (N^\sim)_{(n)} \subset \underline{N}$  where we identify the spaces
%$\underline{N}$ and $\underline{N}^\sim$.  We define
%$M_n\ {\dsize\buildrel{\text{def}}\over =}\ \underline{M}/M_{(n)}$,\
%$N_n\ {\dsize\buildrel{\text{def}}\over =}\ \underline{N}/N_{(n)}$
%and denote by $\pi$ the natural projections $\underline{M} \to M_n$ and
%$\underline{N} \to N_n$.

\proclaim{{\bf 3.3.3}\ Proposition}  If ${\Cal N}^z$  is a generalized Verma
module, $V$  is
a finite-dimensional representation of  ${\bold U}$  and  $M = {\Cal N}^z
\otimes V$  then the natural morphism  ${\Cal N} \otimes V \longrightarrow M_1$
is an isomorphism.
\endproclaim

\demo{Proof}  Follows from 2.4.15.
\enddemo
\vskip .2in

\subheading{3.3.4}  We  denote by
${\Cal E}_z \subset {\Cal C}_z$  be the full subcategory
of modules $M$ such that  dim $M_1 < \infty$.

\proclaim{{\bf 3.3.5}\ Proposition}  a)  For any  $M$ in ${\Cal E}_z$ and any
$n \in {\Bbb N}$
we have  dim $M_n < \infty$.

b)   For any $N$ in ${\Cal O}^+_z$  and any $V$ in ${\Cal D}$  the tensor
product
$N \otimes V$  lies in ${\Cal E}_z$.
\endproclaim

\demo{Proof}  a)  The proof is completely analogous to one in \S\S7.6-7.7 of
[KL].

b) Follows from Proposition 3.3.3.
\enddemo
\vskip .2in

\subheading{3.3.6} For $n \in {\Bbb N}$  we denote by
$L_{n} \subset {\Bbb C}^*$  the
set of eigenvalues of  $T_n$  on  $\overline{M}_{(n-1)}$,
where $\overline{M}_{(n-1)} \subset M_n$  is the image  $M_{(n-1)}$.

\proclaim{Proposition}  For any $M$ in ${\Cal E}_z$  and $n \in {\Bbb N}$
 we have
$$
L_{n+1} \subset {\bigcup_{n \le k \le 3n}} q^k_1 L_1 ,
$$
where  $q_1\ {\buildrel{\text{def}}\over =}\ (z\widetilde{q}^{h^\vee})^{2d}
= (z^dq)^2$.
\endproclaim

\demo{{\bf 3.3.7}\ Proof}  We prove the result by induction in $n$.   If
$n = 0$  then there is nothing to prove.  Assume that we know the
Proposition for $n = n_0$  and we prove it for $n~=~n_0+1$.  Any element
in  $M_{(n-1)}$  is a linear combination of elements of the form $F_im$,
$i \in I$, $m \in M_{(n-2)}$.  Let  $\overline{m}$  be the image of $m$ in
$\overline{M}_{(n-1)}$.  We may assume that  $\overline{m}$  is a generalized
eigenvector of  $T_{n-1}$  on  $\overline{M}_{(n-2)}$  with an eigenvalue equal
to $\lambda$.  It follows then from Proposition 2.2.7 that the image of
$F_im$ in $\overline{M}_{(n-1)}$  is a generalized eigenvector for $T_n$
with the eigenvalue  $\lambda \cdot q_1^{(ii)/2}$.  The inclusion
$\lambda \cdot q_1^{(ii)/2} \in {\bigcup_{n\le k\le 3n}} q^k_1 L_1$
follows now from the inductive assumption and the inclusion  $(ii)/2  \in
\{ 1,2,3\}$.  Proposition 3.3.7 is proved.
\enddemo

%\proclaim{{\bf 3.3.9}\ Proposition}  a)  For any $M$ in ${\Cal C}$  and any
%%pair
%$m > n \in {\Bbb N}$  we have  $T^a_{(m)} M_{(n)} \subset M_{(n)}$.
%
%b)   The induced operator  $T_n$  on $M_n$  does not depend on a choice of $m
%%> n$.
%
%c)  The system $\{ T_n \in End M_n\}$  is compatible with the natural
%projections $M_n \to M_{n-1}$.
%\endproclaim
%
%\demo{Proof}  a) and b) follow from [L] 6.1.1, and c) is obvious.
%\enddemo
\vskip .2in

\subheading{3.3.8} For any $M$ in ${\Cal E}_z$  we define  $\widehat{M}\
{\dsize\buildrel\text{def}\over =}\ {\dsize\lim_\leftarrow} M_n$,
$\widehat{T}\ {\dsize\buildrel\text{def}\over =}\ {\dsize\lim_\leftarrow} T_n
\subset \text{End}(\widehat{M})$  and denote by
$\widehat{\pi}_n: \widehat{M} \to M_n$  the natural projection (see 2.2.10).
Given $M$ in ${\Cal E}_z$    we define (as in [KL], \S29)
for any $\ell \in {\Bbb C}$,  and $n \in {\Bbb N}$  the number
$d_n(\ell)$ to be the dimension of the space ${}_\ell M_n$,
where ${}_\ell M_n\
{\dsize\buildrel\text{def}\over =}\ \cup_m \ker(T_n - \ell)^m$.  It is clear
that for any $\ell \in {\Bbb C}$,  we have
$$
d_1(\ell) \le d_2(\ell) \le \cdots \le
d_n(\ell) \le \cdots .
$$
We define $d(\ell) = {\dsize\lim_{n\to\infty}}
d_n(\ell)$.

\proclaim{{\bf 3.3.9}\ Proposition}  For any $\ell \in {\Bbb C}$,
we have  $d(\ell) < \infty$.
\endproclaim

\demo{Proof}  Since all operators $T_n$  are invertible we may assume that
$\ell \in {\Bbb C}^*$.  It follows from Proposition 3.3.7  that there exists
$n_0 \in {\Bbb N}$  such that for all $n \ge n_0$, $\ell$ is not an eigenvalue
of the restriction of  $T_n$  on  $\overline{M}_{(n-1)}$.
Therefore $d_n(\ell) = d_{n_0}(\ell)$
for all $n \ge n_0$.  Proposition 3.3.9 is proved.
\enddemo

\demo{Remark}  We can rephrase the statement of Proposition 3.3.9 by saying
that
$\widehat{T}$  induces an admissible automorphism of
$\widehat{M}$  (see \S29 in [KL]).
\enddemo
\vskip .2in

\subheading{3.3.10}  For any $M$ in ${\Cal E}_z$  we define a submodule
$\widehat{M}(\infty) \subset \widehat{M}$  as in 2.2.1 and define a subspace
$\widehat{M}^\infty {\dsize\buildrel{\text{def}}\over =}\
{\dsize\bigoplus_{\ell \in {\Bbb C}^*}} {}_\ell \widehat{M}
\subset \widehat{M}$  as in Proposition 29.5  of [KL].  Since
$\widehat{M}(\infty)$  lies in ${}^{[0]}{\Cal C}^+$ one can define the Sugawara
operator
$T \in End\ \widehat{M}(\infty)$.

\proclaim{Proposition} a)  $\widehat{M}^\infty \subset \widehat{M}(\infty)$.

b)  For any $n \in {\Bbb N}$  we have  $\pi_n(\widehat{M}^\infty) = M_n$.
\endproclaim

\demo{Proof}  Part a) follows immediately from Propositions 3.3.9, 3.3.6.
 Part b) follows from Proposition 29.1  in [KL].
\vskip .2in

\subheading{3.3.11}  For any $M$ in ${\Cal E}_z$  the projective system
$\pi_n: M_n \to M_{n-1}$  defines an inductive system
$\pi^*_n:  M^*_{n-1} \hookrightarrow M^*_n$, where
${M}^*_n\ {\dsize\buildrel{\text{def}}\over =}\
{\dsize\bigoplus_{\lambda\in\Lambda_{{\Bbb C}^*}}} \text{Hom}(M_n,{\Bbb
C})_\lambda$
(see 2.1.5).    We denote by $M^*$  the
inductive limit  $M^* = {\dsize\lim_\rightarrow} {M}^*_n$.
It is easy to see that  $M^*$
has a natural structure of  ${\bold U}$-module.  Moreover,  $M^*$
lies in  $^{[0]}{\Cal C}^-$  and therefore we can define the Sugawara operator
$\check{T}~\in~End(M^*)$.

\proclaim{{\bf 3.3.12}\ Proposition}  For any $\zeta \in M^*_n \subset M^*$
and $m \in M$
we have  $\zeta(T_n(\pi_n(m)) = (\check{T}^{-1}\zeta)(m)$  where we consider
$\zeta$ as a linear functional on $M_n$  and  $\check{T}^{-1}\zeta \in M^*$
as a linear functional on $M$.
\endproclaim

\demo{Proof}  Follows from Theorem 3.2.1.
\enddemo

\proclaim{{\bf 3.3.13}\ Proposition}  Let $M$ be an object in ${\Cal E}_z$.

a)   For any  $m \in \widehat{M}(\infty)$  and any
$n \in {\Bbb N}$  we have  $\pi_n(Tm) = T_n \pi_n(m)$.

b)  $\widehat{M}(\infty) = \widehat{M}^\infty$.

c)  The  ${\bold U}$-module  $\widehat{M}(\infty)$  belongs to
${\Cal O}^+_z$.

d)  The natural map  $(\widehat{M}(\infty))_n \to (\widehat{M})_n = M_n$  is an
isomorphism for all $n \in {\Bbb N}$.

e)  $\widehat{M}$  belongs to  ${\Cal E}_z$.
\endproclaim

\demo{Proof}  a)  Follows from 2.2.12.

To prove b) we observe that the inclusion  $\widehat{M}^\infty \subset
\widehat{M}(\infty)$  follows from Proposition 3.3.6 and the inclusion
$\widehat{M}(\infty) \subset \widehat{M}^\infty$  follows from
Proposition 3.1.8 (see the proof of Lemma 26.4 in [KL]).

Part c) follows from Proposition 3.1.8,  part d) follows from
Proposition 3.3.10 b), and part e) is clear.   Proposition 3.3.13 is proved.
\enddemo

\proclaim{{\bf 3.3.14}\ Corollary}  a)  For any module $M$ in ${\Cal E}_z$  the
morphisms
$$
 M \longrightarrow \widehat{M} \longleftarrow \widehat{M}^\infty
$$
induce isomorphisms
$$
M^* \longleftarrow (\widehat{M})^* \longrightarrow (\widehat{M}^\infty)^*.
$$

b)  For any $W$ in ${\Cal O}^-_z$  and any $M$ in
${\Cal E}_z$  the morphisms
$$
\text{Hom}_{{\bold U}}(W, M^\vee) \longleftarrow
\text{Hom}_{{\bold U}}(W, (\widehat{M})^\vee) \longrightarrow
\text{Hom}_{{\bold U}}(W, (\widehat{M}^\infty)^\vee)
$$
are isomorphisms, where as always
$M^\vee\ {\dsize\buildrel{\text{def}}\over =}\ \text{Hom}(M,{\Bbb C})$.
\endproclaim

\demo{Proof}  Part a) is equivalent to part d) of Theorem 3.3.13.
To prove b) we observe that it follows from the  definitions the natural map
Hom$(W,M^*) \to \text{Hom}(W, M^\vee)$  is an isomorphism.  Therefore
part b)  is a restatement of part a).
\enddemo
\vskip .2in

\subheading{3.3.15}  In 3.3.15 - 3.3.18  $\langle M\rangle$  denotes
coinvariants of $M$ with respect to ${\bold U}$.  For any two ${\bold
U}$-modules $W$ and $M$ in ${\Cal C}$
we define the map

\noindent
$\langle M \otimes W\rangle^\vee  \to
\text{Hom}_{{\bold U}}(W,M^\vee)$,  $r \to r^\vee$, where for any
$r \in \langle  M \otimes W\rangle^\vee \subset \text{Hom}(M \otimes W, {\Bbb
C})$
and any  $w \in W$  we define  $r^\vee(w) \in M^\vee$  by the rule
$r^\vee(w)(m)\ {\buildrel{\text{def}}\over =} r(m \otimes w)$.

\proclaim{Proposition}     The map  $\langle M \otimes W\rangle \to
\text{Hom}_{{\bold U}}(W,M^\vee)$  is an isomorphism.
\endproclaim

\demo{Proof}  Clear.
\enddemo
\vskip .2in

\proclaim{{\bf 3.3.16}\ Corollary}  For any  $M$ in ${\Cal E}_z$  and $W$ in
${\Cal O}^-_z$
the morphisms
$$
\langle M \otimes W\rangle \longrightarrow \langle \widehat{M} \otimes W\rangle
\longleftarrow \langle \widehat{M}^\infty \otimes W\rangle
$$
are isomorphisms.
\endproclaim

\subheading{3.3.17} For any finite dimensional  ${\bold U}$-module $X$
we denote by  $\dot{\otimes} X$  the functor from  ${\Cal O}^+_z$  to itself
such that
$$
V \dot{\otimes} X\ {\buildrel \text{def}\over =}\ (V \widehat{\otimes}
X)^\infty
(= (V \widehat{\otimes} X)(\infty))
$$
for all $V$  in ${\Cal O}^+_z$.

For any finite dimensional ${\bold U}$-module $Y$ and $W$ in ${\Cal O}^-_z$
we define  $Y \widehat{\otimes} W\ {\buildrel {\text{def}}\over =}\
{}^\omega (Y \widehat{\otimes} {}^\omega W)$  and $Y \dot{\otimes} W\
{\buildrel {\text{def}}\over =}\ {}^\omega (Y \dot{\otimes} {}^\omega W)$.

\proclaim{{\bf 3.3.18}\ Proposition}  For any $V$ in ${\Cal O}^+_z$, $W$ in
${\Cal O}^-_z$  and a
finite dimensional ${\bold U}$-module $X$ the maps
$$
\langle V \otimes X \otimes W\rangle \longleftarrow \langle (V
\widehat{\otimes} X)
\otimes W\rangle \longrightarrow \langle (V \dot{\otimes} X) \otimes W\rangle
$$
are isomorphisms.
\endproclaim

\demo{Proof}  Follows from Corollary 3.3.16 in the case when $M = V \otimes X$.
\enddemo
\vskip .2in

\subheading{3.3.19}  We say that a ${\bold U}$-module $V$  is locally
$\overline{\bold U}$-finite if for any $v \in V$, $\overline{\bold U}v$
is a finite-dimensional subspace of $V$.

\proclaim{Proposition}  Let  $M$  be a ${\bold U}$-module in ${\Cal E}_z$
which is locally $\overline{\bold U}$-finite.  Then the module
$\widehat{M}^\infty$  is also locally $\overline{\bold U}$-finite.
\endproclaim

\demo{{\bf 3.3.20}\ Proof}  The proposition is an immediate consequence of the
following general and easy result.
\enddemo

\demo{Claim}  Let  $V = (\rho, \underline{V})$  be a locally finite
representation
of ${\bold U}$, $V^{\vee\vee}\ {\dsize\buildrel{\text{def}}\over =}\ (\rho,
\underline{V}^{\vee\vee})$  be the full second dual to $V$ and
$\underline{V}^{**} \subset \underline{V}^{\vee\vee}$  be the subspace of
vectors which are ${\bold U}^+$-finite.  Then  $V^{**}$  is a locally finite
representation of ${\bold U}$.
\enddemo

We will not give a proof of this claim since we will never use Proposition
3.3.19.

\proclaim{{\bf 3.3.21}\ Corollary}  Let  ${\Cal O}^0_z \subset {\Cal O}^+_z$
be the subcategory of locally  $\overline{{\bold U}}$-finite modules.   For any
$V$
in ${\Cal O}^0_z$  and any finite-dimensional representation $X$ of ${\bold U}$
the module  $V \dot{\otimes} X$  lies in ${\Cal O}^0_z$.
\endproclaim

\subheading{3.3.22}  As in [KL] \S27 we denote by  ${\Cal A}_z \subset
{\Cal O}^+_z$  the full subcategory of objects $V$ which admit a filtration
$0 = V_0 \subset V_1 \subset \cdots \subset V_N = V$  such that each quotient
$\overline{V}_n\ {\dsize\buildrel{\text{def}}\over =}\ V_n/V_{n-1}$,\
$1 \le n \le N$, is isomorphic to a Verma module  ${\bold V}^z_{a_i}$,
$a_i \in \overline{\Lambda}_{{\Bbb C}^*}$.  We denote by  $[V]$  the element
of the group ring  ${\Bbb C}[\overline{\Lambda}_{{\Bbb C}^*}]$  defined as a
sum $[V]\ {\dsize\buildrel{\text{def}}\over =}\ {\dsize\sum^N_{n=1}} a_i$.
As follows from the Jordan-H\"older theorem, the element $[V]$  does not
depend on a choice of filtration.

\subheading{3.3.23}  For any finite-dimensional representation $X = (\rho,
\underline{X})$ of ${\bold U}$  we define
$$
\{ X\}\ {\dsize\buildrel{\text{def}}\over =}\
{\dsize\sum_{\overline{\lambda}\in \overline{\Lambda}_{{\Bbb C}^*}}}\
\text{dim}(X_{\overline \lambda}) \cdot {\overline \lambda} \in
{\Bbb C}[\overline{\Lambda}_{{\Bbb C}^*}] ,
$$
where the subspace  $X_{\overline{\lambda}}$  of  $\underline{X}$  is defined
as in 2.1.5.

\proclaim{Proposition}  For any $V$ in ${\Cal A}_z$  and any
${\bold U}$-module $X$ in ${\Cal D}$

a)   $V \dot{\otimes} X \in {\Cal A}_z$  and

b)   $[V \dot{\otimes} X] =   [V] \cdot \{ X\}$.
\endproclaim

\demo{Proof}   The proof of a) is completely analogous to the proof of
Proposition 28.1 in [KL] and the proof of b) is completely analogous to the
proof of Theorem 28.1 in [KL].
\enddemo
\vskip .3in

%this file is qkz3c.tex  [contains section 3.4]

\subheading{3.4\ Comparison of coinvariants}

\subheading{3.4.1}  Given  $V \in {\Cal O}^+_z$, $X, Y \in {\Cal D}$  and
$W \in {\Cal O}^-_z$  we consider\quad $\Gamma$-modules

\noindent
$P(V,X,Y,W)\
{\dsize\buildrel{\text{def}}\over =}\ (V \otimes X)[t] \otimes Y \otimes W$,
\quad $\widehat{P}(V,X,Y,W)\
{\dsize\buildrel{\text{def}}\over =}\ (V \widehat{\otimes} X)[t] \otimes
(Y \widehat{\otimes} W)$\quad    and

\noindent
$Q(V,X,Y,W)\
{\dsize\buildrel{\text{def}}\over =}\ (V \dot{\otimes} X)[t] \otimes
(Y \dot{\otimes} W)$.  We will often write  $\underline{P}, \widehat{P}$
and  $Q$  instead of $P(V,X,Y,W)$,  $\widehat{P}(V,X,Y,W)$  and  $Q(V,X,Y,W)$.
For any  $n \in {\Bbb N}$  we define
${}^nP\ {\dsize{\buildrel{\text{def}}\over =}}\ P \otimes_A A_n$,
${}^n\widehat{P}\ {\dsize{\buildrel{\text{def}}\over =}}\ \widehat{P}
\otimes_A A_n$, and ${}^nQ\ {\dsize{\buildrel{\text{def}}\over =}}\ Q
\otimes_A A_n$.  We denote by $\langle {}^nP\rangle,
\langle {}^n\widehat{P}\rangle$  and  $\langle {}^nQ\rangle$  the
corresponding $A_n$-modules of coinvariants.

The natural imbeddings  $V \otimes X \hookrightarrow V \widehat{\otimes}
X \hookleftarrow V \dot{\otimes} X$  and $Y \otimes W \hookrightarrow
Y \dot{\otimes} W \hookleftarrow Y \dot{\otimes} W$  induce the imbeddings
${}^nP \hookrightarrow {}^n\widehat{P} \hookleftarrow {}^nQ$  and the
corresponding  $A_n$-morphism of $\Gamma$-coinvariants\
$$
\langle {}^nP\rangle\ {\buildrel{{}^n \sigma}\over\longrightarrow}\
\langle {}^n\widehat{P}\rangle\ {\buildrel{{}^n{\eta}}\over\longleftarrow}\
\langle {}^nQ\rangle.
$$

\proclaim{Theorem}  The morphisms  ${}^n\sigma$ and ${}^n\eta$  are
isomorphisms for all
$V \in {\Cal O}^+_z$, $X,Y \in {\Cal D}$  and $W \in {\Cal O}^-_z$.
\endproclaim
\vskip .2in

\demo{{\bf 3.4.2}\ Proof}  Let  ${}^nP_0 \subset {}^nP$  be the
$A_n$-submodule generated by vectors of the form  $\gamma p$,  where
$\gamma \in \Gamma_0$  and  $p \in {}^nP$.  In other words,  ${}^nP_0$
is the kernel of the natural projection  ${}^nP \to \langle {}^nP\rangle$.
We start with the following result.
\enddemo
\vskip .2in

\proclaim{{\bf 3.4.3}\ Lemma}  a)\ $(\underline{V} \otimes \underline{X})_{(n)}
\otimes (\underline{Y} \otimes \underline{W}) \subset {}^nP_0$

b)\ $(\underline{V} \otimes \underline{X}) \otimes
(\underline{Y} \otimes \underline{W})_{(n)} \subset {}^nP_0$.
\endproclaim

\demo{Proof of Lemma}  We prove part a).  The proof of part b) is completely
analogous.  By the definition of the subspace  $(\underline{V} \otimes
\underline{X})_{(n)}$  it is sufficient to show that for any
$\widetilde{v} \in \underline{V} \otimes \underline{X}$,
$\widetilde{w} \in \underline{Y} \otimes \underline{W}$  and
$x \in {\bold f}_\nu$, $|\nu| \ge n$  we have  $x^-\widetilde{v} \otimes
\widetilde{w} \in {}^nP_0$.  But this follows from Proposition 2.4.16.
\enddemo
\vskip .2in

\proclaim{{\bf 3.4.4}\ Corollary}  The map  ${}^n\sigma:  \langle {}^nP\rangle
\to \langle {}^n\widehat{P}\rangle$  is an isomorphism.
\endproclaim
\vskip .2in

\subheading{3.4.5}  Consider the map
$$
{}^n\widetilde{\xi}\ {\dsize\buildrel{\text{def}}\over =}\
({}^n{\sigma})^{-1} \circ {}^n\eta:  \langle {}^nQ\rangle \to
\langle {}^nP\rangle.
$$

\proclaim{Lemma}  The map   ${}^n\widetilde{\xi}$  is surjective.
\endproclaim

\demo{Proof}  Follows from Proposition 3.3.13 d).
\enddemo
\vskip .2in

\proclaim{{\bf 3.4.6}\ Proposition}  Theorem 3.4.1 is true in the case when
$V$ and $W$ are generalized Verma modules.
\endproclaim

\demo{Proof}  In this case we have  $V = {\Cal M}$  and  $W = {\Cal N}^z_-$
where  ${\Cal M}$  and ${\Cal N}$  are nil-modules and it follows from
Proposition 2.5.2 that the $A_n$-module   ${}^nP$  is free of rank equal to
dim$\langle {\Cal M} \otimes X \otimes Y \otimes
{}^\omega {\Cal N}\rangle_{{\bold U}_0}$.
Therefore  dim$_{\Bbb C} \langle {}^nP\rangle = n\cdot \text{dim}\langle
{\Cal M} \otimes X \otimes Y \otimes {}^\omega {\Cal N} \rangle_{{\bold U}_0}$.
On the other hand, it follows from Proposition 3.3.23 and the right exactness
of the functor  $\langle\ ,\ \rangle$  that
dim$_{\Bbb C} \langle {}^nQ\rangle \le n \cdot \text{dim}\langle {\Cal M}
\otimes X \otimes Y \otimes {}^\omega {\Cal N}\rangle_{{\bold U}_0}$.
Therefore Proposition 3.4.6 follows from Lemma 3.4.5.
\vskip .2in

\subheading{3.4.7}  We can now finish the proof of Theorem 3.4.1.  Really for
any  $V \in {\Cal O}^+_z$, $W \in {\Cal O}^-_z$  we can find exact sequences
$M_1 \to M_0 \to V \to 0$  and  $N_1 \to N_0 \to W \to 0$  such that
$M_0,M_1$  are generalized Verma modules in ${\Cal O}^+_z$  and $N_0,N_1$  are
generalized Verma modules in ${\Cal O}^-_z$.  Therefore it follows from
Proposition 3.4.6, the right exactness of the functor $\langle\ ,\ \rangle$
and the five-homomorphisms lemma that  ${}^n\widetilde{\xi}: \langle
{}^nQ\rangle
\to \langle {}^nP\rangle$  is an isomorphism.  Theorem 3.4.1 is proved.
\enddemo

\medpagebreak

\proclaim{{\bf 3.4.8}\ Theorem} Let ${\Cal T} =
{\Cal T}_{(z\widetilde{q}^{h^\vee})^2}$.   The diagram
$$
\matrix \format \c &\quad \c &\quad \c \\
\langle {}^nQ(V,X,Y,W)\rangle & {\buildrel {T^{-1}_{V\dot{\otimes} X} \otimes
id}\over\llongrightarrow} & \langle {}^nQ({\Cal T}(V), {\Cal T}(X),Y,W)\rangle
\\
\vspace{2\jot}
\bigg\downarrow {}^n{\widetilde{\xi}} && \bigg\downarrow {}^n{\widetilde{\xi}}
\\
\vspace{2\jot}
\langle {}^nP(V,X,Y,W)\rangle & {\buildrel {(T^{-1}_{V} \otimes
id)\cdot {}^n\delta_{X,V,Y\otimes W}}\over\llongrightarrow} & \langle
{}^nP({\Cal T}(V), {\Cal T}(X),Y,W)\rangle \\
\endmatrix
$$
is commutative.
\endproclaim

\demo{Proof}  We decompose the diagram into a sequence of simpler diagrams
whose commutativity was proven already.

\newpage

$$
\matrix \format \c &\quad \c &\quad \c \\
\langle {}^nQ(V,X,Y,W)\rangle & {\buildrel {T^{-1}_{V\dot{\otimes} X} \otimes
id}\over\llongrightarrow} & \langle {}^nQ({\Cal T}(V), {\Cal T}(X),Y,W)\rangle
\\
\vspace{1\jot}
\bigg\Vert && \bigg\Vert \\
\vspace{1\jot}
\langle {}^n(V \dot{\otimes} X)[t] \otimes (Y \dot{\otimes} W)\rangle &
{\buildrel {T^{-1}_{V\dot{\otimes} X} \otimes
id}\over\llongrightarrow} & \langle ({\Cal T}(V) \dot{\otimes} {\Cal T}(X))[t]
{\otimes} Y \dot{\otimes} W\rangle \\
\vspace{1\jot}
\bigg\Vert & (1) & \bigg\Vert\\
\vspace{1\jot}
\langle {}^n((V \dot{\otimes} X)[t] \otimes (Y \dot{\otimes} W)\rangle &
{\buildrel {id \otimes \check{T}_{Y\dot{\otimes} W}}\over\llongrightarrow} &
\langle (V \dot{\otimes} X)[t] \otimes {\Cal T}^{-1}(Y \dot{\otimes} W)\rangle
\\
\vspace{1\jot}
\bigg\downarrow & (2) & \bigg\downarrow \\
\vspace{1\jot}
\langle {}^n((V \widehat{\otimes} X)[t] \otimes (Y \dot{\otimes} W)\rangle &
{\buildrel {id \otimes T_{Y\dot{\otimes} W}}\over\llongrightarrow} &
\langle {}^n(V \widehat{\otimes} X)[t] \otimes {\Cal T}^{-1}(Y \dot{\otimes}
W)\rangle \\
\vspace{1\jot}
\bigg\uparrow & (3) & \bigg\uparrow \\
\vspace{1\jot}
\langle {}^n((V \dot{\otimes} X)[t] \otimes (Y \dot{\otimes} W)\rangle &
{\buildrel {id \otimes T_{Y\dot{\otimes} W}}\over\llongrightarrow} &
\langle (V \otimes X)[t] \otimes {\Cal T}^{-1}(Y \dot{\otimes} W)\rangle \\
\vspace{1\jot}
\bigg\Vert & (4) & \bigg\Vert \\
\vspace{1\jot}
\langle {}^n((V {\otimes} X)[t] \otimes (Y \dot{\otimes} W)\rangle &
{\buildrel {(T_V^{-1} \otimes id) \circ {}^n\delta_{X,V,Y\dot{\otimes}
W}}\over\llongrightarrow} &
\langle {}^n((V \otimes X)[t] \otimes {\Cal T}^{-1}(Y \dot{\otimes} W)\rangle
\\
\vspace{1\jot}
\bigg\downarrow & (5) & \bigg\downarrow \\
\vspace{1\jot}
\langle {}^n(V {\otimes} X)[t] \otimes (Y \widehat{\otimes} W)\rangle &
{\buildrel {(T_V^{-1} \otimes id) \circ
{}^n\delta_{X,V,Y\widehat{\otimes} W}}\over\llongrightarrow} &
\langle {}^n((V \otimes X)[t] \otimes {\Cal T}^{-1}(Y \widehat{\otimes}
W))\rangle \\
\vspace{1\jot}
\bigg\uparrow & (6) & \bigg\uparrow \\
\vspace{1\jot}
\langle {}^n((V {\otimes} X)[t] \otimes (Y {\otimes} W)\rangle &
{\buildrel {(T_V^{-1} \otimes id) \circ {}^n\delta_{X,V,Y \otimes
W}}\over\llongrightarrow} &
\langle {}^n((V \otimes X)[t] \otimes {\Cal T}^{-1}(Y \otimes W))\rangle \\
\vspace{1\jot}
\bigg\Vert & (7) & \bigg\Vert  \\
\vspace{1\jot}
{}^n(P(V,X,Y,W) &
{\buildrel {(T_V^{-1} \otimes id) \circ {}^n\delta_{X,V,Y\otimes
W}}\over\llongrightarrow} &
\langle {}^nP({\Cal T}(V), {\Cal T}(X), Y, W)\rangle
\endmatrix
$$
where the commutativity of (1) follows from Theorem 3.2.2, the commutativity
of (4) from Theorem 2.5.7, the vertical map in (2), (3), (5) and (6) are
isomorphisms coming from natural imbeddings (see Theorem 3.4.1) and the right
vertical isomorphism in (7) comes from the natural isomorphism
$\langle M\rangle\ {\buildrel\dsize\sim\over\longrightarrow}\ \langle
{\Cal T}(M)\rangle$ as in 1.1.13.  Theorem 3.4.5 is proved.
\enddemo
\vskip .2in

\newpage

\subheading{3.4.9}  Fix  $z$ as in 3.1.7,  $V$ in ${\Cal O}^+_z$,
$W$ in ${\Cal O}^-_z$, $X,Y$ in ${\Cal D}$  and define  $u\
{\dsize\buildrel{\text{def}}\over =}\ (z\widetilde{q}^{h^\vee})^{-2}$.  For any
$n \in {\Bbb N}$  we define  ${}^nR\ {\dsize\buildrel{\text{def}}\over =}\
\text{Hom}_{A_n}(\langle {}^nQ(V,X,Y,W)\rangle, \langle {}^nP(V,X,Y,W)\rangle$
and

\noindent
${}^nR'\ {\dsize\buildrel{\text{def}}\over =}\
\text{Hom}_{A_n}(\langle {}^nQ({\Cal T}(V), {\Cal T}(X), Y, W)\rangle$,
$\langle {}^nP({\Cal T}(V), {\Cal T}(X), Y, W)\rangle)$  and denote by
${}^n\nabla'$ the $A_n$-linear map from ${}^nR$ to ${}^nR'$  such
that ${}^n\nabla'(a) = (T_V^{-1} \otimes id) \circ {}^n\delta_{V,X,Y\otimes W}
\circ (a \cdot T_{V{\dot{\otimes}} X} \otimes id)$.
\vskip .2in

\subheading{3.4.10}  The $\widetilde{u}$-linear isomorphism  $\widehat{u}$
which we constructed
in 2.4.22 define an $\widetilde{u}$-linear isomorphism $\check{u} =
(\widehat{u})^{-1}$ such that  $\check{u}: {}^nR' \to
{}^nR$  and we define  ${}^n\nabla\ {\dsize\buildrel{\text{def}}\over =}\
\check{u} \circ {}^n\nabla'$.  Then ${}^n\nabla:\ {}^nR \to\ {}^nR$  is a
$\widetilde{u}$-linear automorphism.
We can restate Theorem 3.4.8 as follows.

\proclaim{Theorem}\  ${}^n\nabla({}^n\widetilde{\xi}) =\ {}^n\widetilde{\xi}$.
\endproclaim

 %this file is qkz4a.tex	[contains sections 4.1 and 4.2 and 4.3]

\newpage

\pageno=53

\centerline{\bf \S4.\ Finite-dimensional representations}
\vskip .3in

\subheading{4.1\  The category ${\Cal D}$}

{\bf 4.1.1} As before, we denote by ${\Cal D}$  the category
of unital finite-dimensional ${\bold U}$-modules $M$ such that
$Z$ acts on $M$ as identity.

\demo{Remark} If $\widetilde{\bold U}$  is the quotient of
${\bold U}$  by the ideal generated by $Z-1$ we can identify
${\Cal D}$  with the category of finite dimensional unital
$\widetilde{\bold U}$-modules.
\enddemo

It is clear that ${\Cal D}$  has a natural structure of a strict monoidal
rigid category (see \S1).
\vskip .2in

\subheading{4.1.2}    The algebra
$\widetilde{{\bold U}}$  admits ``loop-like'' generators
$x^\pm_{ik}, h_{ik},  i \in \overline{I},  k \in {\Bbb Z}$  such that

\noindent
$h_{i0} = K_i$,\  $x^\pm_{i0} = E_i$,\quad
$i \in \overline{I}$;  and the generators are subject to certain commutation
relations explained in [D2] (see also [Be1]).  We will use only the following
relations:

$(\alpha)\qquad\qquad\qquad [h_{ik}, h_{j\ell}]\ =\ 0$

$(\beta)\qquad\qquad\qquad  [h_{ik},x^\pm_{j\ell}]\ =\ \pm {\frac{1}{k}}
[ka_{ij}]
\cdot x^\pm_{j,k+\ell}$

$(\gamma)\qquad\qquad\qquad [x^+_{ik},x^-_{j\ell}]\ =\
{\frac{\delta_{ij}(\theta_{i,k+\ell} - \varphi_{i,k+\ell})}{q-q^{-1}}}$,

\noindent
where  $[m] =  {\frac{q^m-q^{-m}}{q-q^{-1}}}$ and the elements $\theta_{i,p}$
and  $\varphi_{i,p}$, $p > 0$, are defined from the relations:
$$
\aligned
exp((q-q^{-1}) \sum_{p>0} h_{ip} \tau^p) &= 1 + \sum_{p>0} K_i^{-1} \theta_{ip}
\tau^p \\
exp((q-q^{-1}) \sum_{p<0} h_{ip} \tau^p) &= 1 + \sum_{p<0} K_i \varphi_{ip}
\tau^p ,
\endaligned
$$
where we consider $\tau$  as a formal parameter.
\vskip .2in

\subheading{4.1.3}  Let  $A^+$ (resp. $A^-$) $\subset \widetilde{\bold U}$ be
the
subalgebra with unity generated by  $x^\pm_{ik},\ h_{ik},$\ \
$i \in \overline{I}$, where $k > 0$  ($k < 0$ resp.).

\proclaim{Proposition} For any $X = (\rho,\underline{X})$ in
${\Cal D}$  we have $\rho(A^+) = \rho(A^-) = \rho(\widetilde{\bold U})$.
\endproclaim

\demo{Proof}  We prove that $\rho(A^+) = \rho(\widetilde{\bold U})$.  The proof
of
the equality $\rho(A^-) = \rho(\widetilde{\bold U})$  is completely analogous.
\enddemo

Let  ${\Cal H}(n) \subset \widetilde{\bold U}$  be the span of
$h_{ik},  k \ge n$, $1 \le i \le r$.
Put  $H = \cap_{n\in{\Bbb N}} \rho({\Cal H}(n))$.   Since
dim$_{\Bbb C} \text{End}\ \underline{X} < \infty$  there exists $N > 0$  such
that
$H = \rho({\Cal H}(N+p))$  for any $p \ge 0$.

Let us prove that for any $\ell < 0$,  $\rho(x^+_{j\ell}) \in \rho(A^+)$.
Fix a vertex $j_0 \in \overline{I}$ such that  $a_{j_0j} \ne 0$.  Then
$$
[\rho(h_{j_0N}), \rho(x^+_{jm})] = {\frac{1}{k}}[ka_{j_0j}]\rho(x^+_{j\ell}),
$$
where $m = \ell - N$.

We have  $\rho(h_{j_0N}) \in H = \rho({\Cal H}(N+p))$, for all
$p \ge 0$  and therefore  $\rho(h_{j_0N})$  lies in $\rho({\Cal H}(N-\ell
+1))$.
Since  $\rho(h_{j_0N})$  lies in $\rho({\Cal H}(N-\ell + 1)$  we
can write  $\rho(h_{j_0N})$  as a linear combination of operators of the
form  $\rho(h_{is})$, $i \in \overline{I}$, $s > N-\ell$.  Therefore
${\frac{1}{k}}[k a_{j,\ell}]\rho(x^+_{j\ell})$  is a linear combination of
operators of the form
$$
[\rho(h_{is}), \rho(x^+_{jm})] = {\frac{1}{k}}[ka_{ij}] x^+_{j,m+s},\quad
s > N-\ell,\ m = \ell - N.
$$
Therefore we see that  ${\frac{1}{k}}[ka_{j_0j}]x^+_{j\ell}$  lies in
$\rho(A^+)$.  Since the number $q$ is not a root of 1 we have  ${\frac{1}{k}}
[ka_{j_0j}] \ne 0$, and therefore $\rho(x^+_{j\ell}) \in \rho(A^+)$.
In a similar way,  one shows that  $\rho(x^-_{j\ell}) \in
\rho(A^+)$  for any $\ell \in {\Bbb Z}$.   $\rho(h_{ik}) \in \rho(A^+)$
for $k \ge 0$  and for any  $i \in \overline{I}$,  we have:  $\rho(\theta_{ik})
\in
\rho(A^+)$  for $k \ge 0$.   The relation
$(\gamma)$  in  4.1.2 implies that $\rho(\varphi_{ik}) \in \rho(A^+)$  for $k
\ge 0$
and all $i \in \overline{I}$.  Since  $\varphi_{ik} = \theta_{ik} = 0$  if  $k
\le 0$
we see that the images of all Drinfeld's generators of  $\widetilde{\bold U}$
belong to  $\rho(A^+)$.  Proposition 4.1.3 is proved.
\vskip .3in

\subheading{4.2\ Endomorphisms of tensor products}

\subheading{4.2.1}  Let  $F = {\Bbb C}[t])$  be the field of fractions of
$A (= {\Bbb C}[t])$, $\overline{A} = {\Bbb C}[[t]]$  be the completion of $A$
and  $\overline{F}$  be the field of fractions of  $\overline{A}$.  It is
clear that   $\underline{X}(t) = \underline{X}[t] \otimes_A F$  and $X((t))
= \underline{X}[t] \otimes_A \overline{F}$  carry the structures of
${\bold U}_F$  and ${\bold U}_{\overline{F}}$-modules respectively.  We
will denote them by  $X(t) = (\rho_{X(t)}, \underline{X}(t))$,
$X((t)) = (\rho_{X((t))}, \underline{X}((t)))$.

For any  $X = (\underline{X}, \rho_X), Y = (\underline{Y}, \rho_Y)$  in
${\Cal D}$,   the $\Gamma$-module structure  $\rho_{X[t]} \otimes \rho_Y$  on
$\underline{X}[t] \otimes_{\Bbb C} \underline{Y}$  defines a $\Gamma$-module
structure  on finite dimensional, respectively, $F$ and $\overline{F}$
vector spaces
$$
\underline{X}(t) \otimes \underline{Y}\ {\buildrel\text{def}\over =}\
(\underline{X}[t] \otimes \underline{Y}) \otimes_AF\quad \text{and}\quad
\underline{X}((t)) \otimes \underline{Y}\ {\buildrel\text{def}\over =}\
(\underline{X}[t] \otimes \underline{Y}) \otimes_A \overline{F}.
$$
We denote the corresponding $\Gamma$-modules by $X(t) \otimes Y  =
(\rho_{X(t)\otimes Y}, \underline{X}(t) \otimes \underline{Y})$\quad  and

\noindent
$X((t)) \otimes Y =  (\rho_{X((t))\otimes Y}, \underline{X}((t)) \otimes
\underline{Y})$.  We denote by  $\rho_{X((t)) \otimes Y}
({\bold U}_{\overline{F}}) \subset
\text{End}(\underline{X}((t)) \otimes \underline{Y})$  the $\overline{F}$-space
of the image of  $\rho_{X((t)) \otimes Y}$.
\vskip .2in

\subheading{4.2.2}  For any $X,Y$  in ${\Cal D}$  we denote (as in 1.1.10)
the ring of endomorphisms of $X$ and $Y$ to be $E_X,E_Y$   and denote
by $E$ the tensor product
$E = E_{X,Y}\ {\dsize\buildrel {\text{def}}\over =}\ E_X \otimes_{\Bbb C} E_Y$.
Let
$$
\overline{\Cal E} {\buildrel {\text{def}}\over =}\ End_{{\bold
U}_{\overline{F}}}(X((t))
\otimes_{\Bbb C} Y), \qquad
{\Cal E}\ {\buildrel {\text{def}}\over =}\ End_{{\bold U}_{F}}(X(t)
\otimes_{\Bbb C} Y).
$$
Then $\overline{\Cal E}$ (resp. ${\Cal E}$) has a natural structure of an
$\overline{F}$
(resp. $F$) module and we have a natural imbedding $E \hookrightarrow
{\Cal E} \hookrightarrow \overline{\Cal E}$.
\vskip .2in

\proclaim{{\bf 4.2.3}\ Theorem}  The natural morphism $E \otimes_{\Bbb C}
\overline{F}
\to \overline{\Cal E}$  is an isomorphism.
\endproclaim

\demo{Proof}  We start with the following obvious result.
\enddemo

\subheading{4.2.4}  Let  $\underline{X},\underline{Y}$  be finite-dimensional
$\overline{F}$-vector spaces,   $B_X \subset End_{\overline F}(\underline{X}),\
B_Y \subset
End_{\overline F}(\underline{Y})$  be ${\overline F}$-subalgebras containing
$id_{\underline X}$
and $id_{\underline Y}$,  $E_X$  and $E_Y$  be centralizers of $B_X$ and
$B_Y$  in  $End_{\overline{F}}(\underline{X})$  and $End_{\overline
F}(\underline{Y})$ respectively.

\proclaim{Lemma}  The centralizer of  $B_X \otimes B_Y$  in
$End_{\overline F}(\underline{X} \otimes \underline{Y}) = End_{\overline
F}(\underline{X})
\otimes_{\overline F} End_{\overline F}(\underline{Y})$  is equal to $E_X
\otimes E_Y$.
\endproclaim

\demo{Proof}  Well known.
\enddemo

Since  $\rho_{X((t))\otimes Y}({\bold U}_{\overline F}) \subset
\rho_{X((t))}({\bold U}_{\overline F}) \otimes_{\overline F}
\rho_Y({\bold U}_{\overline F})$,  Lemma 4.2.4 shows that the following
result implies the validity of Theorem 4.2.3.

\proclaim{{\bf 4.2.5}\ Proposition}  For any  $X = (\rho_X, \underline{X}),
Y = (\rho_Y, \underline{Y})$  in ${\Cal D}$
$$
\rho_{X((t))\otimes Y}(\Gamma) \supset 	\rho_{X((t))}
({\bold U}_{\overline F})\otimes_{\overline F} \rho_Y({\bold U}_{\overline F}).
$$
\endproclaim

\demo{{\bf 4.2.6}\ Proof of Proposition 4.2.5}  Let
$B_X\ {\dsize\buildrel {\text{def}}\over =}\ \rho_X({\bold U}) \subset
End_{\Bbb C}(\underline{X}),$

\noindent
$B_Y\ {\dsize\buildrel {\text{def}}\over =}\
\rho_Y({\bold U}) \subset End_{\Bbb C}(\underline{Y}),$\quad
${\frak p}: End_A(\underline{X}[t]_{\Bbb C}
\otimes_{\Bbb C} \underline{Y}) \to End_{\Bbb C}(\underline{X} \otimes
\underline{Y})$\quad be the reduction mod $m$, where
$m\ {\dsize\buildrel{\text{def}}\over =}\ tA \subset A$  and $\overline{\rho}$
be the composition
$$
\overline{\rho} = {\frak p} \circ \rho_{X[[t]]\otimes Y}:
\Gamma \longrightarrow End_{\Bbb C}(\underline{X} \otimes \underline{Y})  =
End_{\Bbb C}(\underline{X}) \otimes_{\Bbb C} End_{\Bbb C}(\underline{Y}).
$$
Nakayama's lemma shows that Theorem 4.2.3 is a consequence of the following
result.

\proclaim{Claim}  $\overline{\rho}(\Gamma) \supset B_X \otimes B_Y$.
\endproclaim

\demo{{\bf 4.2.7}\ Proof of Claim}     It follows from [Be], Th. 4.7 and
Prop. 5.3 that
for any $k > 0$, $i \in \overline{I}$  we have $h_{ik}, x^\pm_{ik} \in \Gamma$
and moreover for any $\alpha \in \{ h_{ik}, x^\pm_{ik}\}, i \in \overline{I},
k > 0$  we have
$$
(id \otimes \psi)\Delta(\alpha) - \alpha  \otimes 1 \in
{\frak m}({\bold U}_A \otimes_A {\bold U}_A).
$$
Here $\psi$  means the automorphism defined in 2.4.2.
Therefore  $\overline{\rho}(\alpha) = \rho_X(\alpha) \otimes id_Y$  and we
see that $\overline{\rho}(\Gamma) \supset \rho_X(A^+) \otimes id_Y$.
As follows from
Proposition 4.1.3, $\overline{\rho}(\Gamma) \supset B_X \otimes id$.

Analogously we see that for any $i \in \overline{I}, k < 0$  and
$x \in \{ h_{ik}, x^\pm_{ik}\}$  we have  $t^{-k}\alpha \in \Gamma$  and
$(id \otimes \psi) \Delta(t^{-k}\alpha) - 1 \otimes \alpha \in
{\frak m}({\bold U}_A \otimes {\bold U}_A)$. So  $\overline{\rho}(t^{-k}\alpha)
= 1 \otimes \alpha$ and we have  $\overline{\rho}(\Gamma) \supset
id_X \otimes \rho_Y(A^-)$.  As follows from Proposition 4.1.3 we have
$\overline{\rho}(\Gamma)
\supset id_X \otimes B_Y$.   This finishes the proof of the Claim.
Theorem 4.2.3 is proved.

\proclaim{{\bf 4.2.8}\ Corollary}  The natural morphism $E \otimes_{\Bbb C} F
\to {\Cal E}$  is an isomorphism.
\endproclaim

\demo{Proof}  We can consider $\overline{\Cal E}$  and
End$_F(\underline{X}((t)) \otimes
\underline{Y})$  as subspaces in End$_{\overline{F}}(\underline{X}((t))
\otimes \underline{Y})$.  Then  ${\Cal E} = \overline{\Cal E} \cap \text{End}_F
(\underline{X}(t) \otimes \underline{Y})$  and the Corollary follows
immediately from Theorem 4.2.3.
\enddemo

\subheading{4.2.9} Let  $\underline{G}$  be the algebraic ${\Bbb C}$-group such
that  $\underline{G}({\Bbb C})$  is the group of invertible elements in
$E_X \otimes E_Y$.  We can restate Corollary 4.2.7 in the following way.

\proclaim{Corollary}  The natural morphisms  $\underline{G}(F) \to
\text{Aut}(X(t) \otimes Y)$  and  $\underline{G}(\overline{F}) \to
\text{Aut}(X(t)~\otimes~Y)$  are isomorphisms.
\endproclaim

\subheading{4.2.10}  Since the group  $\underline{G}$  is defined over
${\Bbb C}$  any automorphism $\eta$ of the field  $\overline{F}$  preserving
the subfield  $F \subset \overline{F}$  defines group automorphisms of
groups  Aut$(X(t) \otimes Y)$  and  Aut$(X((t)) \otimes Y)$  which we denote
by  $\widehat{\eta}$.
\vskip .3in

\subheading{4.3\  Intertwiners from  $X(t) \otimes Y$  to $Y \otimes X(t)$}

\subheading{4.3.1}  For any $X,Y$  in ${\Cal D}$  we denote by
$$
\overline{\Cal J} \subset Hom_{{\bold U}_{\overline F}}(X((t)) \otimes Y,
Y \otimes X((t))),
$$
$$
{\Cal J} \subset Hom_{{\bold U}_{F}}(X(t) \otimes Y, Y \otimes X(t))
$$
the subsets of morphisms which are isomorphisms.  The group $\underline{G}(F)$
acts
naturally on

\noindent
${\Cal J}: (a,f) \to af$, $a \in {\Cal J}, f \in \underline{G}(F)$.

Let $\underline{Hom}(X(t) \otimes Y, Y \otimes X(t))$  be the $F$-linear
space $Hom_{{\bold U}_{F}}(X(t) \otimes Y, Y \otimes X(t))$  considered
as an algebraic variety.
\vskip .2in

\proclaim{{\bf 4.3.2}\ Lemma}  a)  ${\Cal J}$  is a set of $F$-points
of a Zariski open subset  $\underline{\Cal J}$  in a linear subspace of
$\underline{Hom}(X(t) \otimes Y, Y \otimes X(t))$,  $\overline{\Cal J} =
\underline{\Cal J}({\overline F})$
and the action of $\underline{G}(F)$  on ${\Cal J}$  comes from an algebraic
action of
$\underline{G}$ on $\underline{\Cal J}$.

b)  Either $\underline{\Cal J}$  is empty or  $\underline{G}$  acts  simply
transitively on $\underline{\Cal J}$.
\endproclaim

\demo{Proof}  Clear.
\enddemo
\vskip .2in

\proclaim{{\bf 4.3.3}\ Proposition}  For any $s \in \overline{\Cal J}$  there
exists
an element $f \in \underline{G}({\overline F})$ and such that $s = s_0f$, where
$s_0
\in {\Cal J}$.
\endproclaim

\demo{Proof}  Since the algebraic ${\Bbb C}$-group $\underline{G}$  is a
group of invertible elements in an algebra $\underline{E}$ (see 4.2.2)
it is isomorphic
to a semidirect product  $\underline{L} \ltimes \underline{N}$,
where $\underline{L}$  is a direct product of a number of general linear
groups and $\underline{N}$  is a unipotent group.  Let  $\underline{G}(F)$
be the $F$-group obtained from  $\underline{G}$  by the extension of scalars
from ${\Bbb C}$ to $F$.   Then $\underline{G}(F)$ is also a semidirect
product of $\underline{L}^*_{F}$   and $\underline{N}^*_{F}$ , where
$\underline{L}^*_{F}$   is a direct product of general linear groups
and $\underline{N}^*_{F}$   is unipotent.  As follows from Corollary 4.2.8
and Lemma 4.3.2 b)\quad $\underline{\Cal J}$  is a principal homogeneous
$\underline{G}$-space.  It follows from  Proposition 1.33,  Propositions 3.1
and 3.6 in [S]  that  ${\Cal J} \ne \phi$  and moreover $\overline{\Cal J} =
\underline{G}({\overline F}){\Cal J}$.   Proposition 4.3.3
is proved.
\enddemo
\vskip .3in

%this files is qkz4b.tex [contains sections 4.4]

\subheading{4.4\  The functional equation}

\subheading{4.4.1}  For any  $X = (\rho_X, \underline{X})$,\
$Y = (\rho_Y, \underline{Y})$ in ${\Cal D}$ we denote by $s_{X[[t]],Y} \in
\overline{\Cal J}$  the inverse limit of elements ${}^ns_{X[t],Y}$, where
${}^ns_{X[t], Y}$ are defined in 2.4.19.   It can be treated as a linear
map  $\underline{X}[[t]] \otimes \underline{Y} \to \underline{Y} \otimes
\underline{X}[[t]]$, where  $\underline{X}[[t]] = \underline{X} \otimes
\overline{A}$.

\proclaim{Lemma}    There exists an element  ${\frak f} \in
\underline{G}(\overline{F})$  such that  $s_{X[[t]],Y}\cdot {\frak f}^{-1}
\in {\Cal J}$ and such an element is unique up to a left multiplication
by an element in  $\underline{G}(\overline{F})$.
\endproclaim

\demo{Proof}  Follows from Proposition 4.3.3.
\enddemo

We denote by ${\Cal L}_{X,Y}$
the subset of all elements  ${\frak f} \in \underline{G}(\overline{F})$  such
that  $s_{X[[t]],Y} \cdot {\frak f}^{-1} \in {\Cal J}$.

\proclaim{Corollary}  ${\Cal L}_{X,Y}$  is a left $\underline{G}(F)$  coset
in $\underline{G}(\overline{F})$.
\endproclaim

\subheading{4.4.2}  Let $\eta$  be the continuous automorphism of the field
$\overline{F}$  over ${\Bbb C}$  such that  $\eta(t) =
\widetilde{q}^{2h^\vee} \cdot t$.   Let $\widehat{\eta}$
be the corresponding automorphisms of the groups $\underline{G}(F)$  and
$\underline{G}(\overline{F})$  (see 4.2.10).

\proclaim{Theorem}  $\widehat{\eta}({\Cal L}_{X,Y}) = {\Cal L}_{X,Y}$.
\endproclaim
\vskip .2in

\subheading{4.4.3}  We start the proof of Theorem 4.4.2 with the following
observation.  Let
$$
\aligned
& {\Cal J}^* \subset \text{Hom}_{{\bold U}_{{F}}}(X(t) \otimes Y^*, Y^* \otimes
X(t)),\qquad
\overline{\Cal J}^* \subset \text{Hom}_{{\bold U}_{\overline{F}}}(X((t))
\otimes Y^*, Y^* \otimes X((t))),\\
& {\Cal J}^{**} \subset \text{Hom}_{{\bold U}_{{F}}}(X(t) \otimes Y^{**},
Y^{**} \otimes X(t))\qquad
\overline{\Cal J}^{**} \subset \text{Hom}_{{\bold U}_{\overline{F}}}(X((t))
\otimes Y^{**},
Y^{**} \otimes X((t))
\endaligned
$$
be the subsets of isomorphisms as in 4.3.1 and we denote by $\underline{\Cal
J}^*$
the algebraic subvariety of the linear space
$\underline{\text{Hom}}(X(t) \otimes Y^*, Y^* \otimes X(t))$  corresponding to
${\Cal J}^*$.

There exists an action
$$
\phi: (E_X^{op} \otimes E_Y^{op}) \times \text{Hom}_{{\bold U}_F}(X(t) \otimes
Y^*, Y^*
\otimes X(t)) \to \text{Hom}_{{\bold U}_F}(X(t) \otimes Y^*, Y^* \otimes X(t))
$$
of the algebra  $E_X^{op} \otimes E_Y^{op}$
on the space  Hom$_{{\bold U}_F}(X(t) \otimes Y^*, Y^* \otimes X(t))$  such
that
$$
\phi(a, {\frak f}_X \otimes {\frak f}_Y) = (id \otimes {\frak f}_X) \circ a
\circ (id \otimes {\frak f}^*_Y)
$$
for\  ${\frak f}_X \in E_X, {\frak f}_Y \in E_Y,$
$a \in \text{Hom}_{{\bold U}_F}(X(t) \otimes Y^*, Y^* \otimes X(t))$.  The
action $\phi$
of the algebra  $(E_X \otimes E_Y)^{op}$  on  Hom$(X(t) \otimes Y^*,
Y^* \otimes X(t))$  defines an algebraic action  $(a,g) \longmapsto a*g$ of
$\underline{G}(F)$  on $\underline{\Cal J}^*$  such that
$a * g\ {\dsize\buildrel{\text{def}}\over =}\ \phi(a,g^{-1})$  for
$a \in {\Cal J}^*, g \in \underline{G}(F)$.
\vskip .3in

\subheading{4.4.4}  We consider  $s_{X[[t]],Y^*}$  as an element in
$\overline{\Cal J}^*$  and denote by  ${\Cal M}_{X,Y^*} \subset
\underline{G}(\overline{F})$  the set of all  $g \in
\underline{G}(\overline{F})$
such that $(s_{X[[t]],Y^*})*(g^{-1}) \in {\Cal J}^*$.

\proclaim{Lemma}  ${\Cal M}_{X,Y^*}$  is a left $\underline{G}(F)$  coset in
$\underline{G}(\overline{F})$.
\endproclaim

\demo{Proof}  Clear.
\enddemo
\vskip .2in

\proclaim{{\bf 4.4.5}\ Proposition}  ${\Cal M}_{X,Y^*} = {\Cal L}_{X,Y}$.
\endproclaim

\demo{Proof of Proposition}  Consider the action $\psi$  of the algebra
$E_X \otimes E_Y$  on the space  Hom$_{{\bold U}_F}(Y^* \otimes X(t), X(t)
\otimes Y^*)$
such that  $\psi(b, {\frak f}_X \otimes {\frak f}_Y) = (id \otimes {\frak
f}^*_Y)
\circ b \circ (id \otimes f_X)$  for

\noindent
${\frak f}_X \in E_X,
{\frak f}_Y \in E_Y$,\ $b \in \text{Hom}_{{\bold U}_F}(Y^* \otimes X(t), X(t)
\otimes Y^*)$.
\enddemo

\proclaim{Lemma}  For any  $a \in {\Cal J}^*$, $g \in G$  we have
$(a*g)^{-1} = \psi(a^{-1}, g)$.
\endproclaim

\demo{Proof}  Follows from the definitions.
\enddemo

\subheading{4.4.6}\ The action $\psi$  defines an algebraic action
$(b,g) \longmapsto b \diamondsuit g$  of $\underline{G}(F)$  on

\noindent
Hom$_{{\bold U}_F}(Y^* \otimes X(t), X(t) \otimes Y^*)$  such that $b
\diamondsuit g =
\psi(b,g)$  for $b \in \text{Hom}_{{\bold U}_F}(Y^* \otimes X(t), X(t) \otimes
Y^*)$,\
$g \in \underline{G}(F)$.  As follows from Lemma 4.4.5 we have
$$
{\Cal M}_{X,Y^*} = \{ g \in \underline{G}(\overline{F})|
s^{-1}_{X[[t]],Y^*} \diamondsuit g^{-1} \in \text{Hom}(Y^* \otimes X(t), X(t)
\otimes Y^* \}.
$$

\subheading{4.4.7}  Now we can prove Proposition 4.4.5.  Choose any
${\frak f}_{X,Y} \in {\Cal L}_{X,Y}$.  Then we have $s_{X,Y} = a \cdot
{\frak f}_{X,Y}$, where  $a \in {\Cal J}$.  As follows from Proposition
1.1.18 we have $\varphi^{Y,Y}_{X[[t]]}(s_{X[[t]],Y}) = s^{-1}_{X[[t]],Y^*}$.
So $s^{-1}_{X[[t]],Y^*}	= \varphi^{Y,Y}_{X[[t]]}(a \cdot {\frak f}_{X,Y})$.
By Proposition 1.1.10 the right side is equal to  $\varphi^{Y,Y}_{X[[t]]}(a)
\diamondsuit {\frak f}_{X,Y}$.  Since $\varphi^{Y,Y}_{X[[t]]}$  maps
Hom$_{{\bold U}_F}(X(t) \otimes Y, Y \otimes X(t))$  into

\noindent
Hom$_{{\bold U}_F}(Y^* \otimes X(t),
X(t) \otimes Y^*)$  we see that ${\frak f}_{X,Y} \in {\Cal M}_{X,Y^*}$.
Proposition 4.4.5 is proved.

\proclaim{{\bf 4.4.8}\ Proposition}  ${\Cal L}_{X,Y^*} = {\Cal M}_{X,Y}$  for
all $X,Y$ in ${\Cal D}$.
\endproclaim

\demo{Proof}  Completely analogous to the proof of Proposition 4.4.5.
\enddemo

\proclaim{Corollary}  ${\Cal L}_{X,Y^{**}} = {\Cal L}_{X,Y}$.
\endproclaim

\subheading{4.4.9}  For  $u \in {\Bbb C}^*$  we denote by  $\theta_u$  the
automorphism of the field $F$  such that  $\theta_u(t) = ut$.  We also denote
by  $\theta_u$  the continuous extension of  $\theta_u$  on $\overline{F}$.
The
automorphism  $\theta_u$  define automorphisms of the groups
$\underline{G}(F)$
and  $\underline{G}(\overline{F})$  which we denote by  $\widehat{\theta}_u$.

\proclaim{Lemma}  For any  $X,Y$ in ${\Cal D}$  and  $u \in {\Bbb C}^*$  we
have
$$
{\Cal L}_{{\Cal T}_u(X),Y} = \widehat{\theta}_u({\Cal L}_{X,Y})\quad
\text{and}\quad
{\Cal L}_{X, {\Cal T}_u(Y)} = \widehat{\theta}_u^{-1}({\Cal L}_{X,Y}).
$$
\endproclaim

\demo{Proof}  Follows immediately from the definitions.
\enddemo

\subheading{4.4.10}  Now we can finish the proof of Theorem 4.4.2.  Really
it follows
from Corollary 4.4.8, Lemma 4.4.9 and Proposition 2.1.15.  Theorem 4.4.2 is
proved.

\proclaim{Corollary}  For any  ${\frak f}_{X,Y} \in {\Cal L}_{X,Y}$  we have
$$
\widehat{\eta}({\frak f}_{X,Y}){\frak f}^{-1}_{X,Y} \in \underline{G}(F).
\tag *
$$
\endproclaim

\proclaim{{\bf 4.4.11} Lemma}  There exists  ${\frak f}_{X,Y} \in {\Cal
L}_{X,Y}$
such that  ${\frak f}_{X,Y}$  and $\underline{X}[[t]] \otimes \underline{Y}$
preserve

\noindent
$\underline{X}[[t]] \otimes \underline{Y} \subset \underline{X}((t))
\otimes \underline{Y}$  and  ${\frak f}_{X,Y}(0) = id$.
\endproclaim

\demo{Proof}  The standard (profinite) topology on $\overline{F}$  defines a
topology on the group  $\underline{G}(\overline{F})$.  Since the group
$\underline{G}$  is a semidirect product $\underline{L} \ltimes N$, where
$\underline{L}$  is the direct product of a number of copies of
$\underline{GL}_n$
and  $\underline{N}$  is a unipotent group, the subgroup $\underline{G}(F)$  is
dense in $\underline{G}(\overline{F})$.  Lemma 4.4.11 follows now immediately
from the equality  $s_{X[[t]],Y} \equiv \sigma$ (mod ${\frak m}$), where
$\sigma: \underline{X} \otimes \underline{Y} \to \underline{Y} \otimes
\underline{X}$
is permutation as in 2.3.3.
\enddemo

\subheading{4.4.12}  Let  $A_{mer}$  be the subring of the power series which
converges in a neighborhood of zero and define a meromorphic function in
${\Bbb C}$.

\proclaim{Proposition}  For any  $X, Y$ in ${\Cal D}$  there exists
${\frak f}_0 \in \underline{G} (A_{mer})$  which satisfies the
condition of Lemma 4.4.11.
\endproclaim

\demo{Proof}  Choose any  ${\frak f}_{X,Y} \in {\Cal L}_{X,Y}$  and define
$r\ {\dsize\buildrel{\text{def}}\over =}\ \widehat{\eta}({\frak f}_{X,Y})
{\frak f}^{-1}_{X,Y}$.  As follows from Theorem 4.4.2 we have  $r \in
\underline{G}({\Bbb C}(t) \cap \overline{A})$  and  $r(0) = Id$.  Consider
the equation
$$
\widehat{\eta}({\frak f}){\frak f}^{-1} = r    \tag **
$$
on ${\frak f} \in \underline{G}(\overline{A})$.  It is easy to see that there
exists a unique solution  ${\frak f}_0$  of (**) in
$\underline{G}(\overline{A})$
such that  ${\frak f}_0(0) = Id$  and moreover any solution of (**) has a form
${\frak f} = {\frak f}_0 \cdot c$  for some  $c \in \underline{G}({\Bbb C})$.
One can write ${\frak f}_0 \in \underline{G}(\overline{A})$  as an
infinite convenient product  ${\frak f}_0(t) = {\dsize\prod^\infty_{n=0}}
r^{-1}(\widetilde{q}^{2h^\vee n}t)$.  It is easy to see that this product is
also convergent in  $\underline{G}(A_{mer})$.  Proposition 4.4.12 is
proved.
\enddemo

\proclaim{Corollary} The map $s_{X[[t]],Y}$  from $\underline{X} \otimes_{\Bbb
C}
\underline{Y} \otimes_{\Bbb C} \overline{A}$  to  $\underline{Y} \otimes
\underline{X} \otimes \overline{A}$  maps the subspace

\noindent
$\underline{X} \otimes_{\Bbb C} \underline{Y} \otimes_{\Bbb C} A_{mer}$ to
$\underline{Y} \otimes_{\Bbb C} \underline{X} \otimes_{\Bbb C} A_{mer}$.
\endproclaim
\vskip .2in

\subheading{4.4.13}  As follows from the previous Corollary we can consider
$s = s_{X[[t]],Y}$  as a mermorphic function in $t$ with values in the space
Hom$_{\Bbb C}(\underline{X} \otimes \underline{Y}, \underline{Y} \otimes
\underline{X})$.  The union of the set of poles of $s$ and the points
$t \in {\Bbb C}$ such that  $s(t)$  is not an isomorphism is denoted by
$\Lambda_{X,Y}$.

\proclaim{Lemma}  $\Lambda_{X,Y}$  is a discrete subset of ${\Bbb C}$.
\endproclaim

\demo{Proof}  It is clear that either  $\Lambda_{X,Y} = {\Bbb C}$  or it is
discrete.  Since  $s(0)$  is invertible,

\noindent
$\Lambda_{X,Y} \ne {\Bbb C}$.  The Lemma is proved.
\enddemo
\vskip .2in

\subheading{4.4.14}  Let  $\widetilde{\Cal D}^{(2)}$  be the full subcategory
of the pair of objects in $(X,Y)$  in ${\Cal D}$  such that $1 \notin
\Lambda_{X,Y}$.  As follows from Lemma 4.4.13 for all pairs $(X,Y)$  in ${\Cal
D}$
we have

\noindent
$({\Cal T}_u(X), Y) \in \widetilde{\Cal D}^{(2)}$  for all $u$ in a
dense open set  ${\Bbb C}^* - \Lambda_{X,Y}$.  We can consider
$(\widetilde{\Cal D}^{(2)},s)$  as a weak  braiding on ${\Cal D}$  (see
1.2.1).

\proclaim{Proposition}  The data $({\Cal C}, {\Cal D}, {\Cal O}^\pm_z, s_\pm,
{\Cal T}_\pm, \widetilde{\Cal D}^{(2)}, s)$  is a rigid $KZ$-data for all $z$
such that  $z\widetilde{q}^{h^\vee} \in {\Bbb C}^* - {\frak Q}_{\ge  0}$, where
${\Cal T}_\pm = {\Cal T}_{z^{\pm 1}}$,\ \
$s_{\pm(V,X)}\ {\dsize\buildrel{\text{def}}\over =}\ (s_{V,X})^{\pm 1}$,
and  $s_{V,X}$ is  as in 2.3.3.
\endproclaim

\demo{Proposition}  Follows immediately from the definitions and 4.4.14, 2.3.8.
\enddemo

\demo{{\bf 4.4.15\ Remark.}} Similarly to 4.4.3 one can prove that ``the
square of weak braiding'' $s_{X,Y}(t) s_{Y,X} (t^{-1})$  is an elliptic
function on the curve  ${\Bbb C}^*/\widetilde{q}^{2h^\vee {\Bbb Z}}$  with
values in  Hom$_{\Bbb C}(\underline{X} \otimes \underline{Y}, \underline{Y}
\otimes \underline{X})$.
\enddemo

%this file is qkz5.tex which is section 5

\newpage

\pageno=62

\centerline{\bf \S5.\ The quasi-associativity morphism.}

\subheading{5.1\ The definition}

\subheading{5.1.1}  Fix $z$ such that $z\widetilde{q}^{h^\vee} \in
{\Bbb C}^* - {\frak Q}_{\ge 0}$  and ${\bold U}$-modules
$V$ in ${\Cal O}^+_z$, $W$ in ${\Cal O}^-_z$, $X,Y$ in ${\Cal D}$.   We
assume that  $(z\widetilde{q}^{h^\vee})^{2m} \notin \Lambda_{X,Y}$  for
$m \in {\Bbb Z}$.  Let 	$P = P(V,X,Y,W)$  be the
$\Gamma$-module as in 3.4.1 and  $\langle P\rangle$  the corresponding
$A$-module of coinvariants.

\proclaim{Proposition}  $\langle P\rangle$  is a finitely generated
$A$-module.
\endproclaim

\demo{Proof}  As follows from Proposition 3.1.8, there exist exact sequences
$M_1' \to M_0' \to V \to 0$  and  $N_1' \to N_0' \to W \to 0$  where  $M'_k$
and $N'_k$  are generalized Verma modules in ${\Cal O}^+_z$  and ${\Cal
O}^-_z$,
respectively, $k = 0,1$.  Therefore Proposition 5.1.1
follows from Lemma 2.5.5  and the right exactness of the functor
$\langle\quad \rangle$.
\enddemo

\subheading{5.1.2}  Let  $Q, {}^nP, {}^nQ$  be $\Gamma$-modules as in 3.4.1
and  $\langle Q\rangle, \langle {}^nP\rangle$ and $\langle {}^nQ\rangle$  be
the  corresponding modules of coinvariants.

\proclaim{Proposition}  a)  $\langle Q\rangle$  is a finitely generated
$A$-module.

b)  The natural $A_n$-module morphisms
$$
{}^n\pi_P: \langle P\rangle \otimes_A A_n \to \langle {}^nP\rangle\quad
\text{and}\quad  {}^n\pi_Q: \langle Q\rangle \otimes_A A_n \to
\langle {}^nQ\rangle
$$
are isomorphisms for all $n \in {\Bbb N}$.
\endproclaim

\demo{Proof}  a) Follows from Proposition 5.1.1 and Theorem 3.4.1
and b) follows from Proposition 2.5.4.
\enddemo

\subheading{5.1.3}  For any  $n \in {\Bbb N}$  we denote by
${}^n\xi: \langle P\rangle \otimes_A A_n \to \langle Q\rangle \otimes_A A_n$
the composition  ${}^n\xi = ({}^n\pi_P)^{-1} \circ {}^n \widetilde{\xi} \circ
{}^n\pi_Q$  where  ${}^n\widetilde{\xi}$  is the isomorphism from
$\langle {}^nQ\rangle$  to $\langle {}^nP\rangle$  as in 3.4.5.  It follows
from Theorem 3.4.1 and Proposition 5.1.2 that  ${}^n\xi$  is an isomorphism
for all $n \in {\Bbb N}$.  It is clear that the isomorphisms ${}^n\xi$  are
compatible with the natural projections $\langle P\rangle \otimes_A
A_{n+1} \to \langle P\rangle \otimes_A A_n$ and  $\langle Q\rangle \otimes_A
A_{n+1} \to \langle Q\rangle \otimes_A A_n$.  Therefore the family
${}^n\xi$  defines an isomorphism $\overline{\xi}:
\langle \overline{Q}\rangle \to \langle \overline{P}\rangle$ of
$\overline{A}$-modules, where  $\langle \overline{Q}\rangle\
{\dsize\buildrel{\text{def}}\over =}\
{\dsize\lim_{\leftarrow}} \langle Q\rangle \otimes_A A_n$,
$\langle \overline{P}\rangle\ {\dsize\buildrel{\text{def}}\over =}\
{\dsize\lim_{\leftarrow}} \langle P\rangle \otimes_A A_n$ and
$\overline{A} = {\dsize\lim_{\leftarrow}} A_n$.

Since $\langle P\rangle, \langle Q\rangle$  are finitely generated $A$-modules
we can identify the $\overline{A}$-modules  $\langle \overline{P}\rangle$,
$\langle \overline{Q}\rangle$ with the tensor products $\langle P\rangle
\otimes_A \overline{A}$ and $\langle Q\rangle \otimes_A \overline{A}$
correspondingly.

Let  $A_{mer}$  be the ring as in Section 4.

\subheading{5.1.4}  Let  $\langle P\rangle_{mer}\
{\dsize\buildrel{\text{def}}\over =}\
\langle P\rangle \otimes_A A_{mer},$ $\langle Q\rangle_{mer}\
{\dsize\buildrel{\text{def}}\over =}\
\langle Q\rangle \otimes_A A_{mer}.$   We can identify $\langle
\overline{P}\rangle$
with $\langle P\rangle_{mer} \otimes_{A_{mer}} \overline{A}$  and
$\langle \overline{Q}\rangle$ with $\langle Q\rangle_{mer} \otimes_{A_{mer}}
\overline{A}$.

\proclaim{Theorem}    There exists an isomorphism  $\xi_{mer}: \langle
Q\rangle_{mer}
{\dsize\buildrel{\sim}\over\longrightarrow}\ \langle P\rangle_{mer}$  such that
$\overline{\xi} = \xi_{mer} \otimes 1$.
\endproclaim

\demo{Proof}  For any $n \in {\Bbb N}$  we denote by ${}^nR$  the $A_n$-module
Hom$_{A_n}(\langle {}^nQ\rangle, \langle {}^nP\rangle)$  and  by
${}^n\nabla$  the $\widetilde{u}$-linear automorphism of ${}^nR$  as in
3.4.10.
\endproclaim

Let  $\overline{R}\ {\dsize\buildrel{\text{def}}\over =}\
{\dsize\lim_{\leftarrow}}\
{}^nR$.  We can identify $\overline{R}$   with the  $\overline{A}$-module
\quad
Hom$_{\overline A}(\langle \overline{Q}\rangle, \langle \overline{P}\rangle)$.
Let

\noindent
$R_{mer}\ {\dsize\buildrel{\text{def}}\over =}\  \text{Hom}_{A_{mer}}
(\langle Q\rangle_{mer}, \langle P_{mer}\rangle)$.  One can consider $R_{mer}$
as an
$A_{mer}$-submodule of  $\overline{R}$ and $\overline{R} =
R_{mer} \otimes_{A_{mer}} \overline{A}$.  We denote by  $\overline{\nabla}$
the $\widetilde{u}$-linear automorphism of the module
$\overline{R}\ {\dsize\buildrel{\text{def}}\over =}\ {\dsize\lim_{\leftarrow}}\
{}^nR$  which is the projective limit of ${}^n\nabla$.

\proclaim{5.1.5\ Lemma}\ $\overline{\nabla}(R_{mer}) \subset R_{mer}$.
\endproclaim

\demo{Proof}  Follows from Corollary 4.4.12, Proposition 2.4.22 and the
definition of $\overline{\nabla}$.
\enddemo

\subheading{5.1.6}  Our proof of Theorem 5.1.4  is based on the following
result.

\proclaim{Proposition}  Let  $\varphi: D_R \to M_n({\Bbb C})$  be a holomorphic
function in a disc

\noindent
$|t| \subset D_R = \{ t \in {\Bbb C}\big\vert |t| < R\}$\  such that
$\det\ \varphi(0) \ne 0$,\ $\psi: {\Bbb C} \to GL_n({\Bbb C})$  be
a polynomial function, $p$ a complex number such that  $0 < |p| < 1$  and
$F(t) \in \overline{A} \otimes_{\Bbb C} M_n({\Bbb C})$  be a formal power
series solution of the difference equation
$$
F(pt)\psi(t) = \varphi(t) F(t)   \tag *
$$
such that  $F(0) = Id$.  Then the series  $F(t)$  is convergent in the disc
$|t| < {\frac{R}{\Vert \varphi(0)\Vert + 1}}$.   Moreover if  $\varphi(t)$
has a meromorphic continuation to ${\Bbb C}$  then $F(t)$  has a meromorphic
continuation to ${\Bbb C}$.
\endproclaim

\demo{Proof}  As follows from (*) we have  $\varphi(0) = \psi(0)$.  Let us
write the expansions for $\varphi, \psi$ and $F$
$$
\varphi(t) = \sum^\infty_{j=0} \varphi_jt^j,\quad
\psi(t) = \sum^N_{j=0} \psi_jt^j,\quad
F(t) = \sum^\infty_{j=0} f_jt^j .
$$
Then we can  rewrite (*) in the form
$$
p^jf_j \varphi(0) - \varphi(0)f_j  =
\sum^j_{k=1} \varphi_k f_{j-k} - \sum^N_{k=1} f_{j-k} \psi_k.
\tag **
$$
Since $\varphi$ is convergent in $D_R$  we have $\Vert \varphi_j\Vert \le
r^{-j}$  for all  $r < R$.  Let  $K_j$  be an endomorphism of
$M_n({\Bbb C})$  such that $K_jM\ {\buildrel{\text{def}}\over =}\
p^jM\varphi(0) - \varphi(0)M$ for $M \subset M_n({\Bbb C})$
and $K > 0$  be a constant such that
$\Vert K_j^{-1}\Vert < K$  for $j \gg 0$.  Let us prove that
$\Vert f_j\Vert < C(r/K')^{-j}$   suitable constants $K'$ and $C$.
This will imply the validity of the first part of Proposition 5.1.6.
\enddemo

\subheading{5.1.7}   Applying the norm to both sides of (**) we get an
inequality
(for $j \gg 0)$)
$$
\Vert f_j\Vert \le K'( \sum^j_{k=1} r^{-k}\Vert f_{j-k}\Vert).
$$
This implies that the sequence $\Vert f_j\Vert$ is dominated by the
sequence $g_j$  which satisfies the equalities $g_j = K'({\dsize\sum^j_{k=1}}
r^{-k} g_{j-k})$  for all $j \gg 0$.  The sequence  $g_j$  satisfies the
equation  $rg_{j+1}/K' - g_j/K' = g_j$ and therefore $g_j$  is a geometric
progression.
So $f_j$  is dominated by a geometric progression and therefore  $F(t)$  is
convergent in a neighborhood of 0.  If  $\varphi$  has a meromorphic
continuation to ${\Bbb C}$, it follows then from (*) that $F$
has a meromorphic continuation to ${\Bbb C}$.  Proposition 5.1.6 is proved.

\proclaim{{\bf 5.1.8}\ Corollary}  Let $L,M$  be finitely generated
$A$-modules, $\varphi \in (End\ L) \otimes_A A_{mer}$,\ $\psi \in (End\ M)$
and  $\overline{\xi} \in \text{Hom}(L,M) \otimes_A \overline{A}$  be such that
$$
\varphi(0) \in Aut\ L/tL,\quad \psi(0) \in Aut\ M/t\overline{M},\quad
\overline{\xi}(0) \in Isom(L/tL, M/tM)
$$
and the formal power series  $\overline{\xi}$  satisfy the equation
$$
\overline{\xi}(pt) = \varphi(t) \overline{\xi}(t) \psi(t).
$$
then $\overline{\xi} \in \text{Hom}(L,M) \otimes_A A_{mer}$.
\endproclaim

\vskip .2in

\subheading{5.1.9}  We can now prove Theorem 5.1.4.  We put
$p = (z\widetilde{q}^{h^\vee})^{-1}$.  As follows from
Proposition 3.1.8 it is sufficient to prove Theorem 5.1.4 in the case when
$V$ and $W$ are generalized Verma modules.  In this case it follows from
Proposition 2.5.2 that $P$ and $Q$ are free $A$-modules and the result
follows from Theorem 3.4.8 and Corollary 5.1.8.  Theorem 5.1.4 is proved.
\vskip .3in

\newpage

\subheading{5.2.1}  Let  $z, V, W, X$ and $Y$  satisfy the assumptions of
5.1.1.  It follows from the Theorem 5.1.4 that there exists a small disc
$|t| < \epsilon$  such that  $\xi_{mer}$  defines an isomorphism of
$\langle Q\rangle_{mer}$  and $\langle P\rangle_{mer}$  at any point of this
disc.  Since  $|z\widetilde{q}^h|> 1$  we can find an even integer $m < 0$
such
that  $|(z\widetilde{q}^{h^\vee})^m| < \epsilon$  so we get an isomorphism
of  $\langle Q\rangle$  and $\langle P\rangle$  at
$t_0 = (z\widetilde{q}^{h^\vee})^m$.    Iterating  $\widetilde{u}$-linear
isomorphism  $\overline{\nabla}$  from Lemma 5.1.5 and using Corollary
2.4.22 and the condition  $(z\widetilde{q}^{h^\vee})^n \notin
\Lambda_{X,Y}$, for  $n \in {\Bbb Z}$  we obtain the isomorphism of
$\langle Q\rangle$  and $\langle P\rangle$  at  $t_1 = 1$.

\proclaim{{\bf 5.2.2}\ Proposition}  We have constructed a functorial
isomorphism
$$
\nabla:  \langle V \dot{\otimes} X,\ Y \dot{\otimes} W\rangle_{\bold U}
\simeq \langle V,X,Y,W\rangle_{\bold U}
$$
where  $\langle \cdot\ \rangle_{\bold U}$  denotes the vector space of
${\bold U}$-coinvariants.
\endproclaim

\medpagebreak

\proclaim{{\bf 5.2.3.}\ Corollary}  Let  $z,V,X$ and $Y$  be as in 5.1.1.
Then the isomorphism $\nabla$  from 5.2.2 gives rise to the functorial
quasi-associativity constraint
$$
a_{V,X,Y}:  (V \dot{\otimes} X) \dot{\otimes} Y \simeq V \dot{\otimes}
(X \otimes Y).
$$
It satisfies the pentagon axiom (see [M]) with respect to $X$ and $Y$.
\endproclaim

\demo{Proof}  Follows from 5.2.2 in the same way as the isomorphism
18.2 (b) follows from the theorem 17.29 in [KL].
\enddemo

\demo{{\bf 5.2.4}\ Remark}  The following is the tautological reformulation
of 5.2.3.  We say that  $X$  and  $Y$  in ${\Cal D}$  are in generic position
with respect to the elliptic curve  $E =
{\Bbb C}^*/(\widetilde{q}^{h^\vee}z)^{2{\Bbb Z}}$  if the image of
$\Lambda_{X,Y}$  under the natural projection  ${\Bbb C}^* \to E$  does not
contain the unity of the group $E$.  If $X$ and $Y$ are in generic position
with respect to the elliptic curve $E$ then 5.2.3 holds.

%this file is qkzref.tex

\newpage

\pageno=66

\centerline{\bf References}
\vskip .3in

\baselineskip=12pt

\item{[Be]}  Beck, J.:  Braid group action and quantum affine algebras.
Commun. Math. Phys., 165: 3 (1994), 555-568.
\vskip .2in

\item{[D1]}  Drinfeld, V.G.:  On quasitriangular quasi-Hopf algebras closely
related to $Gal(\overline{Q}/Q)$, Algebra i Anal. 2 (1990), 149-181.
\vskip .2in

\item{[D2]}  Drinfeld, V.:  A new realization of Yangians and quantized
affine algebras.  Preprint of the Phys.-Tech. Inst. of Low Temperatures,
1986.
\vskip .2in

\item{[D3]}  Drinfeld, V.: Quantum groups.  Proc. ICM (Berkeley, vol. 2,
pp. 789-819.
\vskip .2in

\item{[DM]} Deligne, P., Milne, J.S.:  Tannakian Categories, Lecture Notes
in Math. 900, Springer-Verlag, 1981, 104-228.
\vskip .2in

\item{[FR]} Frenkel, I., Reshetikhin, N.:  Quantum affine algebras and
holonomic difference equations.  Commun. Math. Phys., 146  (1992), 1-60.
\vskip .2in

\item{[K]}  Kac, V.:  Infinite-dimensional Lie algebras.  Cambridge University
Press (1991).
\vskip .2in

\item{[KL]}  Kazhdan, D.,  Lusztig, G.:  Tensor structures arising from affine
Lie algebras, I-IV.  American J. Math., 1993-1994.
\vskip .2in

\item{[L]}  Lusztig, G.:  Introduction to quantum groups.  Birkh\"auser, 1993.
\vskip .2in

\item{[M]} Mac Lane, S.:  Categories for the working mathematician.
Springer-Verlag
(1971).
\vskip .2in

\item{[S]}  Serre, J.-P.:  Cohomologie Galoisienne, Coll\`ege de France, 1963.
\vskip .2in

\item{[TUY]} Tsuchiya, A., Ueno, K., and Yamada, Y., Conformal field theory on
universal family of stable courses with gauge symmetries, Advanced Stud. Pure
Math., vol. 19, Academic Press, Boston, MA, 1989, pp. 459-565.

\bye